\documentclass[12pt]{article}

\usepackage{amsmath}
\usepackage{amssymb}
\usepackage{amsthm}
\usepackage{geometry}
\usepackage{titlesec}
\usepackage[numbers,sort&compress]{natbib}
\usepackage{tabularx}
\usepackage{booktabs}
\usepackage{mathtools}
\usepackage{enumitem}
\usepackage{cases}
\usepackage{threeparttable}
\usepackage{float}

\numberwithin{equation}{section}
\newtheorem{theorem}{Theorem}[section]
\newtheorem{lemma}[theorem]{Lemma}
\newtheorem{proposition}[theorem]{Proposition}
\newtheorem{corollary}[theorem]{Corollary}
\theoremstyle{definition}
\newtheorem{definition}[theorem]{Definition}
\newtheorem{remark}[theorem]{Remark}

\theoremstyle{plain}

\titleformat*{\section}{\large\bfseries}
\titleformat*{\subsection}{\normalsize\bfseries}

\newcommand{\bC}{\ensuremath{\mathbb{C}}}

\newcommand{\bN}{\ensuremath{\mathbb{N}}}
\newcommand{\bZ}{\ensuremath{\mathbb{Z}}}
\newcommand{\bZz}{\ensuremath{\mathbb{Z}_2}}
\newcommand{\subno}[1]{\noindent\textnormal{(#1)}}
\newcommand{\afCase}[1]{\textbf{Case~#1:}}

\newcommand{\End}{\mathrm{End }}

\newcommand{\Span}{\mathrm{span }}

\newcommand{\Rad}{\mathrm{Rad }}
\newcommand{\Ann}{\mathrm{Ann }}
\newcommand{\ep}{\ensuremath{{\bar{0}}}}
\newcommand{\op}{\ensuremath{{\bar{1}}}}

\newcommand{\fg}{\mathfrak{g}}
\newcommand{\Der}{\mathrm{Der }}
\newcommand{\Id}{\mathrm{Id }}

\newcommand{\fs}{\mathfrak{s}}

\newcommand{\cA}{\ensuremath{\mathcal{A}}}
\newcommand{\cQ}{\ensuremath{\mathcal{Q}}}
\newcommand{\cL}{\ensuremath{\mathcal{L}}}

\newcommand{\bone}{\mathbf{1 }}
\newcommand{\wt}{\mathrm{wt }}
\newcommand{\normord}[1]{%
  \mathopen{\vcenter{\hbox{$\genfrac{}{}{0pt}{}{\circ}{\circ}$}}}%
  #1%
  \mathclose{\vcenter{\hbox{$\genfrac{}{}{0pt}{}{\circ}{\circ}$}}}%
}

\newcommand{\fJ}{\mathfrak{J}}

\begin{document}
\begin{center}
    \bfseries\large
		Canonical weight-two algebras and strongly generated vertex operator superalgebras
		\footnote{Supported by the National Natural Science Foundation of China (No. 12471027) and Natural Science Foundation of Shanghai (No. 24ZR1471900).
		\\ \indent \ \ $^*$Corresponding author: Xiaoqing Yue (xiaoqingyue@tongji.edu.cn)}

    \mdseries\normalsize
    \bigskip
		Jinrong Wang, Xiaoqing Yue$^*$
		
    \footnotesize
    \smallskip
		School of Mathematical Sciences, Tongji University, Shanghai 200092, China,\\
        Key Laboratory of Intelligent Computing and  Applications(Tongji University), Ministry of Education
		
    \smallskip
		E-mails: 2130926@tongji.edu.cn, xiaoqingyue@tongji.edu.cn 
\end{center}
    {
    \footnotesize
    \noindent\textbf{Abstract}. This paper studies vertex operator superalgebras strongly generated by their weight-1, weight-$\frac{3}{2}$ and weight-2 homogeneous subspaces. The presence of a nontrivial weight-1 subspace introduces derivative and quadratic current terms in weight-2. To address this, we use low-weight quasi-primary projections to identify two canonical commutative associative algebras, $\mathcal{A}$ and $A$, within the weight-2 subspace. We demonstrate that the algebraic structures of the generating subspaces are encoded by a quadruple that generalizes Lie--Rinehart algebras and their weak representations. This quadruple yields necessary and sufficient compatibility conditions for constructing an associated Lie superalgebra $\mathcal{L}(\mathfrak{g},A,U)$. Based on this, we construct a family of universal vertex operator superalgebras $V_{\mathcal{L}}(\ell,0)$ and establish a reconstruction theorem for their simple quotients under certain conditions. If $A$ is finite-dimensional, we prove that its orthogonal idempotent decomposition into local Artinian summands naturally induces a tensor product isomorphism of the corresponding vertex operator superalgebras. In addition, our general construction encompasses many vertex operator superalgebras, such as superconformal algebras and generalized $\mathrm{BMS}_3$ (super)algebras.

    
    }

\section{Introduction}\label{sect:intro}
\hspace{1.5em}Vertex (operator) superalgebras, which are mathematical counterparts of chiral algebras in conformal field theory, have deep connections with infinite-dimensional Lie (super)algebras, finite simple groups and mathematical physics \cite{B,FB,FHL,K98,Li,T,Xu}. The notion of strong generation was introduced by Kac \cite{K98} to study singular operator product expansions. In brief, a vertex (super)algebra $V$ is said to be strongly generated by a subset $T$ if $V$ is linearly spanned by the vacuum vector $\bone$ and vectors $u_{-n_1}^{(1)}u_{-n_2}^{(2)}\cdots u_{-n_r}^{(r)}\bone$ for $r\ge 1$, $u^{(i)}\in T$, $n_i\ge1$.

It is well known that an affine vertex algebra can be characterized as a $\bZ$-graded vertex algebra $V=\bigoplus_{n\in\bZ}V_{(n)}$ which is strongly generated by its weight-one subspace $V_{(1)}$ (cf. \cite{FZ,Lian}). On the other hand, vertex operator algebras strongly generated by their weight-two subspaces have been studied in \cite{Lam,P}. More recently, Li and Yu studied vertex operator superalgebras strongly generated by their weight-$\frac{3}{2}$ and weight-two subspaces in \cite{LY}. Among their main results, they proved that if such a superalgebra is simple, the weight-two subspace possesses a canonical commutative associative algebra structure equipped with a non-degenerate symmetric invariant bilinear form. Meanwhile, the weight-$\frac{3}{2}$ subspace is naturally a module over this algebra, satisfying a set of compatibility conditions.

However, many prominent examples in representation theory and mathematical physics possess nontrivial weight-one currents alongside generators of weights $\frac{3}{2}$ and $2$. For instance, the $N=2$ and small $N=4$ superconformal algebras \cite{K98,KPP,KW}, which play important roles in conformal field theory, exhibit precisely this low-weight structure. Moreover, various supersymmetric and Weyl extensions of the $\mathrm{BMS}_3$ algebra also fit into a low-weight framework in which the generating weights belong to $\{1,\frac{3}{2},2\}$ \cite{AST,O,Man}. The $\mathrm{BMS}_3$ algebra and its extensions have attracted considerable attention in connection with three-dimensional flat gravity, asymptotic symmetries and flat holography. On the mathematical side, the $\mathrm{BMS}_3$ algebra is known as the $W(2,2)$ algebra in \cite{ZD}, arising naturally in the classification of simple vertex operator algebras generated by a two-dimensional weight-two subspace.

Motivated by these considerations, we investigate conformal vertex superalgebras $V=\bigoplus_{q\in\frac12\bZ}V_{(q)}$ strongly generated by $V_{(1)}+V_{(\frac32)}+V_{(2)}$. The presence of weight-one generators introduces substantial difficulties. 

First, the weight-two space contains derivatives and normally ordered products of weight-one vectors. Therefore, the commutative associative algebra structure appearing in \cite{Lam,LY} does not extend directly to the entire weight-two space. In particular, the zeroth products of weight-two vectors may have nonzero components in the quasi-primary subspace $\Omega$ of $V_{(3)}$. Imposing $\Omega=0$ would force $V_{(1)}=0$ whenever the central charge $c_V\neq-2$ (see Lemma~\ref{lem:c=-2}). 

Second, strong generation in these weights does not by itself imply that the singular products of the generators close linearly. For instance, the zeroth product of two weight-$\frac32$ vectors may contain quadratic current terms and such terms occur in certain minimal $W$-algebras (see Remark~\ref{rmk:linear.closure}).

To address these issues, we first identify two canonical subspaces $\cA$ and $A$ of $V_{(2)}$. The space $\cA$ consists of the quasi-primary vectors $a\in V_{(2)}$ satisfying $a_0V_{(2)}\subset L(-1)V_{(2)}$, and $A$ is obtained by imposing
the additional condition $a_0V_{(\frac{3}{2})}\subset L(-1)V_{(\frac{3}{2})}$. We show in Theorem~\ref{thm:canonical.algebras} that they are commutative associative algebras with identity $\omega$ and a symmetric invariant bilinear form $\beta$. We then consider the case in which $V_{(1)}+V_{(\frac{3}{2})}+A$ strongly generates $V$ and
\[
    u_0v\in A+L(-1)V_{(1)},\ \ \forall\ u,v\in V_{(\frac{3}{2})}.
\]
The second condition restricts the quasi-primary component of $u_0v$ to $A$ and ensures that the odd--odd singular products close linearly. When $V_{(1)}=0$, strong generation gives $A=V_{(2)}$, so the second condition is automatic. In this case, our relations in Proposition~\ref{prop:mode.relations} recover those in \cite{LY}. 

Denote $\mathfrak{g}=V_{(1)}$ and $U=V_{(\frac{3}{2})}$. Under the above assumptions, to describe the resulting algebraic structures, we associate to $V$ a quadruple
\[
    \cQ(V)=\big((\fg,[\cdot,\cdot],\kappa),(A,\circ,\beta),(\rho,\cdot,\theta),(U,\tau,\eta,\sigma)\big),
\]
which consists of the Lie algebra $(\fg,[\cdot,\cdot])$, the commutative associative algebra $(A,\circ)$, the mutual actions between $\fg$ and $A$, a compatible $(\fg,A)$-module structure on $U$, and six bilinear maps determined by the singular products of the generators (see Tables~\ref{tab:product.actions}--\ref{tab:compatibility}). Notably, the mutual actions between $\mathfrak{g}$ and $A$ are closely related to Lie--Rinehart algebras \cite{Hue,Rine}, except that the derivation action $\rho:\fg\to\Der(A)$ satisfies a symmetry condition rather than $A$-linearity (see Remark~\ref{rmk:Lie-Rine}). Moreover, when $(\fg,A)$ is a Lie--Rinehart algebra, the compatible actions on $U$ define a weak representation of Lie-Rinehart algebras in the sense of \cite{PSTZ} (see Remark~\ref{rmk:weak.rep}). Conversely, starting from an abstract $Z$-valued quadruple $\mathcal{Q}_Z$ for an even vector space $Z$, we prove in Theorem~\ref{thm:construction.L} that the compatibility conditions listed in Table~\ref{tab:compatibility} are necessary and sufficient for the bracket relations to make
\[
    \cL(\fg,A,U)=\big((\fg\oplus A\oplus U)\otimes\bC[t,t^{-1}]\big)\oplus Z
\]
into a Lie superalgebra.

For each linear functional $\ell\in Z^*$, we construct an induced vacuum module for the Lie superalgebra $\cL(\fg,A,U)$ with each $c\in Z$ acts as the scalar $\ell(c)$ and endow it with a vertex superalgebra structure, denoted by $V_\cL(\ell,0)$. When $\fg$, $A$ and $U$ are finite-dimensional, $V_{\cL}(\ell,0)$ is a vertex operator superalgebra of CFT-type with central charge $2\beta_\ell(e,e)$, where $e$ denotes the identity of $A$ and $\beta_\ell$ is the composition of $\beta$ and $\ell$ (see Theorem~\ref{thm:VOSA}). Assuming further that $\fg$ is perfect, we show that $L(1)\fg=0$ and study the unique simple quotient $L_\cL(\ell,0)$ of $V_{\cL}(\ell,0)$ by means of its normalized symmetric invariant bilinear form. Since $A$ is in general a proper subspace of $V_{(2)}$, the restriction of the form $\beta$ to $A$ may be degenerate even when $V$ is simple. Consequently, the non-degeneracy of $\beta$ is replaced by the condition $\Rad\, \beta\cap \Ann_A(\fg)=0$ (see Proposition~\ref{prop:simple->}). Conversely, if $\kappa_\ell$ and $\tau_\ell$ are non-degenerate and $\Rad\, \beta_\ell\cap \Ann_A(\fg)=0$, then the quotient map $V_\cL(\ell,0)\rightarrow L_\cL(\ell,0)$ is injective on $\fg\oplus A\oplus U$, where $X_\ell$ is the composition of $X$ and $\ell$, $X\in\{\kappa,\tau,\beta\}$ (see Proposition~\ref{prop:simple<-}). Finally, we obtain the reconstruction theorem $L_\cL(\ell_0,0)\cong V$ for simple vertex operator superalgebras, where we take $Z=\bC$ and $\ell_0=\Id\in Z^*$ (see Theorem~\ref{thm:reconstruction}).

In particular, we illustrate the role of the canonical commutative associative algebra $A$ in the finite-dimensional setting. In this case, $A$ is Artinian and thus admits a decomposition into local Artinian algebras via orthogonal idempotents $A=A_1\oplus\cdots\oplus A_r$ for some $r\ge1$. We prove that this decomposition at the algebra level naturally induces a direct sum decomposition of the associated Lie superalgebra $\cL(\fg,A,U)/Z$ and a tensor product isomorphism $V_{\cL}(\ell,0)\cong V_{\cL_1}(\ell,0)\otimes\cdots\otimes V_{\cL_r}(\ell,0)$ of vertex operator superalgebras (see Theorem~\ref{thm:tensor_decomp}). This reduces the study of $V_{\cL}(\ell,0)$ to the case where $A$ is a local Artinian algebra.

Our general construction encompasses many important vertex operator superalgebras. These include the vertex operator (super)algebras associated with the affine-Virasoro algebra, the twisted Heisenberg--Virasoro algebra, the $N=2$ and small $N=4$ superconformal algebras, and truncated Virasoro and Neveu--Schwarz algebras. Furthermore, we construct the generalized Weyl--BMS algebras and present a novel supersymmetric example in which the pair $(\fg,A)$ strictly generalizes the Lie--Rinehart algebra structure.

This paper is organized as follows. We first review basic notations and facts about vertex (operator) superalgebras in Section~\ref{sect:prelim}. In Section~\ref{sect:Construction}, we determine the low-weight quasi-primary decompositions and construct the canonical commutative associative algebras $\cA$ and $A$. In Section~\ref{sect:data}, we study the conformal vertex superalgebras $V$ which are strongly generated by their subspaces $V_{(1)}+V_{(\frac{3}{2})}+A$. We show that the associated algebraic structures can be encoded in a quadruple $\cQ(V)$. The related definitions and compatibility conditions are listed in Tables~\ref{tab:product.actions}--\ref{tab:compatibility}. Under these conditions, we construct a Lie superalgebra $\cL(\fg,A,U)$. In Section~\ref{sect:VOSA}, we associate vertex (operator) superalgebras $V_\cL(\ell,0)$ to the Lie superalgebra $\cL(\fg,A,U)$. Under certain conditions, the reconstruction theorem for simple quotients is obtained. Section~\ref{sect:examples} establishes a tensor product decomposition and presents several examples.

Throughout this paper, we use notations $\bC$, $\bZ$, $\bZ_+$ and $\bN$ to represent the set of complex numbers, integers, positive integers and non-negative integers, respectively. All vector spaces are defined over $\bC$. For $m\in\bZ$ and $r\in\bZ_+$, we use the convention $\binom{m}{0}=1$ and $\binom{m}{r}=\frac{m(m-1)\cdots(m-r+1)}{r!}$.

\section{Preliminaries}\label{sect:prelim}
\hspace{1.5em}In this section, we recall some classical notions and results. Let $V$ be a superspace that is a \bZz-graded linear space with a direct sum $V=V_\ep \oplus V_\op$. If $x\in V_\theta$, $\theta\in\bZz=\{\ep,\op\}$, then we say $x$ is \emph{homogeneous} and of \emph{parity $\theta$}. For $V=V_\ep\oplus V_\op$, a homogeneous element of $V_\ep$ is said to be \emph{even} while an element of $V_\op$ is said to be \emph{odd}. The parity of a homogeneous element $x$ is denoted by $|x|$. Throughout what follows, if $|x|$ occurs in an expression, then it is assumed that $x$ is homogeneous and the expression extends to the other elements by linearity.

\begin{definition}
    A \emph{Lie superalgebra} $\fg=\fg_\ep\oplus\fg_\op$ is a vector superspace equipped with a bilinear map $[\cdot,\cdot]:\fg\times\fg\rightarrow\fg$ called super bracket, satisfying the following axioms $(x,y,z\in\fg)$:
    \begin{align*}
        &[\fg_i,\fg_j]\subset\fg_{i+j},\ \ \forall\ i,j\in\bZz,\\
        &[x,y]=-(-1)^{|x||y|}[y,x],\\
        &[x,[y,z]]=[[x,y],z]+(-1)^{|x||y|}[y,[x,z]].
    \end{align*}
\end{definition}

Now we recall some definitions and results related to vertex operator superalgebras (cf. \cite{B,FB,K98,Li,T}).
\begin{definition}
    A \emph{conformal vertex superalgebra} is a quadruple $(V,Y,\bone,\omega)$, which consists of a $\frac{1}{2}\bZ$-graded vector space $V=\bigoplus_{q\in\frac{1}{2}\bZ}V_{(q)}$, viewed as a superspace $V=V^\ep\oplus V^\op$, where 
    \begin{equation*}
        V^\ep=\bigoplus_{n\in\bZ}V_{(n)}\ \ \text{and }\ \ V^\op=\bigoplus_{n\in\bZ} V_{(n+\frac{1}{2})},
    \end{equation*}
    and equipped with a linear map
    \begin{align*}
        Y(\cdot,z):V&\longrightarrow(\End V)[[z,z^{-1}]]\\
        v&\longmapsto Y(v,z)=\sum_{n\in\bZ}v_n z^{-n-1},\ \ v_n\in \End V,
    \end{align*}
    and two distinguished even vectors: \emph{vacuum vector} $\bone\in V_{(0)}$ and \emph{conformal vector} $\omega\in V_{(2)}$, satisfying the following conditions for $u,v\in V$:

    \subno{1} the \emph{lower truncation condition}:
    \begin{equation*}
        u_nv=0\ \ \text{for}\ n\ \text{sufficiently large};
    \end{equation*}

    \subno{2} the \emph{vacuum and creation properties}:
    \begin{equation}\label{vac.cre}
        \bone_mv=\delta_{m,-1}v,\ \ v_{-1}\bone=v\ \ \text{and}\ \ v_n\bone=0,\ \ m\in\bZ,\,n\in\bN;
    \end{equation}

    \subno{3} the \emph{Jacobi identity}:
    \begin{align}
        z_0^{-1}\delta\left(\frac{z_1-z_2}{z_0}\right)&Y(u,z_1)Y(v,z_2)-(-1)^{|u||v|}z_0^{-1}\delta\left(\frac{z_2-z_1}{-z_0}\right)Y(v,z_2)Y(u,z_1)\notag\\
        =z_2^{-1}&\delta\left(\frac{z_1-z_0}{z_2}\right)Y\big(Y(u,z_0)v,z_2\big),\label{Jac.1}
    \end{align}
    where $\delta(z)=\sum_{n\in\bZ}z^n$ and $(z_i-z_j)^n$ is expanded as a formal power series in $z_j$;

    \subno{4} the \emph{Virasoro algebra relations}:
    \begin{equation}\label{Vir}
        [L(m),L(n)]=(m-n)L(m+n)+\frac{m^3-m}{12}\delta_{m+n,0}c_V,\ \ m,n\in\bZ,
    \end{equation}
    where $L(n)=\omega_{n+1}$ for $n\in\bZ$ and $c_V\in\bC$ is the \emph{central charge} of $V$, satisfying the $L(-1)$-\emph{derivative property}:
    \begin{align}
        &Y\big(L(-1)v,z\big)=\frac{d}{dz}Y(v,z),\ \ v\in V;
   \shortintertext{and}
        &V_{(q)}=\{v\in V \mid  L(0)v=qv\},\ \ q\in\frac{1}{2}\bZ.\label{L(0)}
    \end{align}
    For $v\in V_{(q)}$ with $q\in\frac{1}{2}\bZ$, we say $v$ is \emph{conformally} or $\frac{1}{2}\bZ$-homogeneous of (\emph{conformal}) \emph{weight} $q$ and write $\wt v=q$.
\end{definition}

Throughout this paper, Lie algebras and conformal vertex algebras are regarded as Lie superalgebras and conformal vertex superalgebras with zero odd part, respectively.

\begin{definition}
    A \emph{vertex operator superalgebra} is a conformal vertex superalgebra satisfying the two grading restriction conditions:
    \begin{equation*}
        \dim V_{(q)}<\infty\ \ \text{for}\ q\in\frac{1}{2}\bZ\ \ \text{and }\ V_{(q)}=0\ \ \text{for}\ q\ \text{sufficiently negative}.
    \end{equation*}
    A vertex operator superalgebra $V=\bigoplus_{q\in\frac12\bZ}V_{(q)}$ is said to be of \emph{CFT-type} if $V_{(0)}=\bC\bone$ and $V_{(q)}=0$ for $q<0$.
\end{definition}

\begin{definition}
    A \emph{homomorphism} from a vertex operator superalgebra $V_1$ to another vertex operator superalgebra $V_2$ is a linear map $f:V_1\rightarrow V_2$ such that
    \begin{equation*}
        f(\bone_1)=\bone_2,\ \ f(\omega_1)=\omega_2,\ \  f\big(Y_1(u,z)v\big)=Y_2\big(f(u),z\big)f(v),\ \ \forall\ u,v\in V_1.
    \end{equation*}
    A bijective homomorphism is called an \emph{isomorphism}. 
\end{definition}

Let $V$ be a conformal vertex superalgebra. In terms of components, Jacobi identity \eqref{Jac.1} can be rewritten as
\begin{align*}
    \sum_{i=0}^\infty(-1)^i&\binom{m}{i}\big(u_{m+k-i}v_{n+i}-(-1)^{|u||v|}(-1)^mv_{m+n-i}u_{k+i}\big)\\
    &=\sum_{i=0}^\infty\binom{k}{i}(u_{m+i}v)_{n+k-i},
\end{align*}
for $u,v\in V$, $m,n,k\in\bZ$. As immediate consequences, for $u,v\in V$, $m,n\in\bZ$, we have:

\subno{1} the \emph{Borcherds commutator formula}
\begin{equation}\label{comm}
    [u_m,v_n]=u_mv_n-(-1)^{|u||v|}v_nu_m=\sum_{i=0}^\infty\binom{m}{i}(u_iv)_{m+n-i};
\end{equation}

\subno{2} the \emph{Borcherds iterate formula}
\begin{equation}\label{iter}
    (u_mv)_n=\sum_{i=0}^\infty(-1)^i\binom{m}{i}\big(u_{m-i}v_{n+i}-(-1)^{|u||v|}(-1)^mv_{m+n-i}u_i\big);
\end{equation}

\subno{3} the \emph{skew-symmetry}
\begin{equation*}
    Y(u,z)v=(-1)^{|u||v|}e^{zL(-1)}Y(v,-z)u,
\end{equation*}
which amounts to
\begin{equation}\label{skew.symm}
    u_nv=(-1)^{|u||v|}(-1)^{n+1}\left(v_nu+\sum_{i\ge1}\frac{1}{i!}(-1)^iL(-1)^iv_{n+i}u\right).
\end{equation}

Recalling that $ L(n)=\omega_{n+1}$ for any $n\in\bZ$, we have 
\begin{gather}
    [L(-1),v_n]=\big(L(-1)v\big)_n=-nv_{n-1},\label{der}\\
    u_nV_{(q)}\subset V_{(p+q-n-1)},\label{wt}
\end{gather}
for $v\in V$, $u\in V_{(p)}$, $n\in\bZ$ and $p,q\in\frac{1}{2}\bZ$.

A nonzero vector $v$ of $V$ is called a \emph{quasi-primary vector of weight} $h\in\frac{1}{2}\bZ$ (see \cite{FHL}) if $v\in V_{(h)}$ and $L(1)v=0$. If $v$ is a quasi-primary vector of weight $h$, then
\begin{equation}\label{quasi.prim}
    [L(1),v_n]=(2h-n-2)v_{n+1},\ \ \forall\ n\in\bZ.
\end{equation}

The notion of strongly generated vertex (operator) superalgebras was introduced by Kac in his study of operator product expansions and Lie conformal algebras \cite{K98}.

\begin{definition}\textup{(\cite{K98})}
    A conformal vertex superalgebra $V$ is said to be \emph{strongly generated} by a subset $T$ if $V$ is linearly spanned by $\bone$ and the vectors of the form $a^{(1)}_{-n_1}\cdots a^{(r)}_{-n_r}\bone$ with $r\ge 1$, $a^{(i)}\in T$ and $n_i\in\bZ_+$.
\end{definition}

Let $T$ be a generating subspace of $V$. For $u,v\in T$, the products $u_nv$ with $n\geq0$ will be called the \emph{singular products} of $u$ and $v$. By the truncation property, only finitely many of these products are nonzero. We say that the singular products among the elements of $T$ \emph{close linearly} on $T$ if they can be expressed as linear combinations of the vacuum vector, elements of $T$, and their $L(-1)$-derivatives.

The following two results are due to \cite{Lam} and \cite{LY}, respectively.
\begin{proposition}\label{prop:Lam}\textup{(\cite{Lam})}
Let $V=\bigoplus_{n\in\bZ}V_{(n)}$ be a conformal vertex algebra strongly generated by $V_{(2)}$. Then $a_0b=\frac12L(-1)a_1b$ and
\begin{equation*}
    [a_m,b_n]=\frac12(m-n)(a_1b)_{m+n-1}+\binom{m}{3}\langle a,b\rangle\delta_{m+n-2,0},
\end{equation*}
for all $a,b\in V_{(2)}$ and $m,n\in\bZ$, where the
bilinear form $\langle\cdot,\cdot\rangle$ on $V_{(2)}$ is defined by $a_3b=\langle a,b\rangle\bone$.
\end{proposition}

\begin{proposition}\label{prop:LY}\textup{(\cite{LY})}
    Let $V=\bigoplus_{q\in\frac{1}{2}\bZ}V_{(q)}$ be a conformal vertex superalgebra such that $V_{(\frac{3}{2})}+V_{(2)}$ strongly generates $V$ and such that $\langle V_{(2)}\rangle$ is strongly generated by $V_{(2)}$. Then
    \begin{align*}
        &[a_m,b_n]=\frac{1}{2}(m-n)(a_1b)_{m+n-1}+\binom{m}{3}\langle a,b\rangle\delta_{m+n-2,0},\\
        &[a_m,u_n]=\frac{1}{3}(m-2n)(a_1u)_{m+n-1},\\
        &[u_m,v_n]_+=(u_0v)_{m+n}+\binom{m}{2}\langle u,v\rangle\delta_{m+n-1,0},
    \end{align*}
    for $a,b\in V_{(2)}$, $u,v\in V_{(\frac{3}{2})}$, $m,n\in\bZ$, where $[u_m,v_n]_+=u_mv_n+v_nu_m$. Here $\langle\cdot,\cdot\rangle$ denotes the bilinear forms on $V_{(2)}$ and $V_{(\frac{3}{2})}$ defined by $a_3b=\langle a,b\rangle\bone$ and $u_2v=\langle u,v\rangle\bone$, respectively.
\end{proposition}

\section{Canonical weight-two algebras}\label{sect:Construction}
\hspace{1.5em}We define canonical commutative associative algebras inside the weight-two subspace of a conformal vertex superalgebra strongly generated by the low-weight subspaces. The low-weight decompositions separate quasi-primary vectors from $L(-1)$-derivatives and lead to the subspaces $\cA$ and $A$ of $V_{(2)}$. The main result of this section is Theorem~\ref{thm:canonical.algebras}, which establishes the algebra structures and invariant bilinear forms on $\cA$ and $A$. 

The following result is straightforward.
\begin{lemma}\label{lem:CFT}
    Let $V=\bigoplus_{q\in\frac{1}{2}\bZ}V_{(q)}$ be a conformal vertex superalgebra strongly generated by $V_{(1)}+V_{(\frac{3}{2})}+V_{(2)}$. Then
    \begin{align*}
        V_{(q)}=0\ \text{for}\ q\in-\frac{1}{2}\bZ_+,\ \ V_{(0)}=\bC\bone,\ \ V_{(\frac{1}{2})}=0,
    \end{align*}
    and $V_{(1)}$ carries a Lie algebra structure with Lie bracket $[x,y]=x_0y$ for $x,y\in V_{(1)}$ \textup{(cf. \cite[Theorem~3.7]{Lian})}. If $V$ is further a vertex operator superalgebra, then it is of CFT-type.
\end{lemma}

Assume that $V=\bigoplus_{q\in\frac{1}{2}\bZ}V_{(q)}$ is a conformal vertex superalgebra strongly generated by $V_{(1)}+V_{(\frac{3}{2})}+V_{(2)}$. By Lemma~\ref{lem:CFT}, we can obtain
\begin{equation}\label{CFT}
    V=V_{(0)}\oplus V_{(1)}\oplus V_{(\frac{3}{2})}\oplus V_{(2)}\oplus V_{(\frac{5}{2})}\oplus V_{(3)}\oplus\cdots
\end{equation}
with $V_{(0)}=\bC\bone$, and there exists a Lie algebra $\fg$ such that $V_{(1)}=\fg$. 

Throughout this section, $V$ is a conformal vertex superalgebra strongly generated by $V_{(1)}+V_{(\frac32)}+V_{(2)}$ unless it is stated otherwise. We use the notation
\begin{equation*}
    \fg:=V_{(1)},\ \ U:=V_{(\frac{3}{2})}.
\end{equation*}
The space $\fg$ is even and $U$ is odd. In particular, $L(1)U=V_{(\frac{1}{2})}=0$, whereas $L(1)\fg\subset\bC\bone$ need not vanish. We will also use the facts that $L(-1)\bone=\omega_0\bone=0$ and the injectivity of $L(-1)$ on every positive-weight subspace.

\subsection{Low-weight subspace decompositions}
\hspace{1.5em}We separate the quasi-primary vectors from the $L(-1)$-derivatives in the subspaces of weights $2$, $\frac52$ and $3$, respectively.
\begin{proposition}\label{prop:decom}
    Let $V=\bigoplus_{q\in\frac12\bZ}V_{(q)}$ be a conformal
    vertex superalgebra strongly generated by
    $\fg+U+V_{(2)}$. Set
    \begin{align*}
        \Gamma=\{v\in V_{(2)}\mid L(1)v=0\},\ \Theta=\{v\in V_{(\frac52)}\mid L(1)v=0\},\ \Omega=\{v\in V_{(3)}\mid L(1)v=0\}.
    \end{align*}
    Then
    \begin{align*}
        V_{(2)}=L(-1)\fg\oplus\Gamma,\ \ V_{(\frac52)}=L(-1)U\oplus\Theta,\ \ V_{(3)}=L(-1)V_{(2)}\oplus\Omega.
    \end{align*}
    The projections onto the quasi-primary summands are given by
    \begin{align}
        \pi_\Gamma(v)&=v-\frac12L(-1)L(1)v,&&v\in V_{(2)},\label{proj.Gamma}\\
        \pi_\Theta(u)&=u-\frac13L(-1)L(1)u,&&u\in V_{(\frac52)},\label{proj.Theta}\\
        \pi_\Omega(w)
        &=w-L(-1)\Big(\frac14-\frac1{24}L(-1)L(1)\Big)L(1)w,&&w\in V_{(3)}.\label{proj.Omega}
    \end{align}
\end{proposition}

\begin{proof}
    We prove the assertions for $V_{(3)}$. The others can be verified similarly.

    Let $v\in V_{(3)}$ and set
    \begin{equation*}
        c=\Big(\frac14-\frac1{24}L(-1)L(1)\Big)L(1)v\in V_{(2)},\ \ g=v-L(-1)c.
    \end{equation*}
    We first show that $L(1)g=0$. By \eqref{vac.cre}, \eqref{Vir} and \eqref{L(0)}, we obtain
    \begin{align*}
        L(1)L(-1)c&=\frac14L(1)L(-1)L(1)v-\frac1{24}L(1)L(-1)^2L(1)^2v\\
        &=\frac12L(0)L(1)v+\frac14L(-1)L(1)^2v-\frac1{12}L(0)L(-1)L(1)^2v\\
        &\quad-\frac1{24}L(-1)L(1)L(-1)L(1)^2v\\
        &=L(1)v+\frac14L(-1)L(1)^2v-\frac16L(-1)L(1)^2v-\frac1{12}L(-1)L(0)L(1)^2v\\
        &\quad-\frac1{24}L(-1)^2L(1)^3v\\
        &=L(1)v,
    \end{align*}
    where $L(1)^3v\in V_{(0)}=\bC\bone$. Hence
    \begin{equation*}
        L(1)g=L(1)v-L(1)L(-1)c=0,
    \end{equation*}
    which proves that $V_{(3)}=L(-1)V_{(2)}+\Omega$.

    Now assume that $v\in L(-1)V_{(2)}\cap\Omega$. Then there exists $c'\in V_{(2)}$ such that $v=L(-1)c'$ and $L(1)v=0$. Thus $L(1)L(-1)c'=0$. Applying $L(-1)L(1)$ to this identity gives
    \begin{align*}
        0&=L(-1)L(1)L(1)L(-1)c'\\
        &=2L(-1)L(1)L(0)c'+L(-1)L(1)L(-1)L(1)c'\\
        &=4L(-1)L(1)c'+2L(-1)L(0)L(1)c'+L(-1)^2L(1)^2c'\\
        &=6L(-1)L(1)c'\\
        &=-12L(0)c'+6L(1)L(-1)c'\\
        &=-24c'.
    \end{align*}
    Therefore $c'=0$ and hence $v=0$. This proves that $V_{(3)}=L(-1)V_{(2)}\oplus\Omega$. The expression for $g$ gives the projection formula
    \eqref{proj.Omega}.
\end{proof}

The projections also describe the normally ordered contributions to these spaces.
\begin{lemma}\label{lem:theta,omega}
    With the notation in Proposition~\ref{prop:decom}, we have
    \begin{gather*}
        \Theta=\Span\{\pi_\Theta(x_{-1}u)\mid x\in\fg,\ u\in U\},\\
        \Omega=\Span\{\Omega_{uv},\Omega_{xa},\Omega_{xyz}
        \mid x,y,z\in\fg,\ u,v\in U,\ a\in V_{(2)}\},
    \end{gather*}
    where $\Omega_{uv}=\pi_\Omega(u_{-1}v)$, $\Omega_{xa}=\pi_\Omega(x_{-1}a)$ and $\Omega_{xyz}=\pi_\Omega(x_{-1}y_{-1}z)$.
\end{lemma}
\begin{proof}
    Let $x,y,z\in\fg$ and $u,v\in U$. Strong generation and weight counting show that $V_{(\frac52)}$ is spanned by $L(-1)U$ and the vectors $x_{-1}u$. Applying $\pi_\Theta$ gives the first assertion. Similarly, modulo $L(-1)V_{(2)}$, the weight-three space is spanned by vectors of the forms $u_{-1}v$, $x_{-1}a$, $x_{-1}y_{-1}z$ and $x_{-2}y$. The last type is covered by the second, since
    \begin{equation*}
        x_{-2}y=L(-1)(x_{-1}y)-x_{-1}\big(L(-1)y\big),\ \ L(-1)y\in V_{(2)}.
    \end{equation*}
    Applying $\pi_\Omega$ proves the second assertion.
\end{proof}

\begin{remark}\label{rmk:normal.order}
    For $u,v\in V$, the special case $u_{-1}v$ can be expressed as the normal order product of $Y(u,z)$ and $Y(v,z)$ as $Y(u_{-1}v,z)$:
    \begin{equation*}
        Y(u_{-1}v,z)=\normord{Y(u,z)Y(v,z)}.
    \end{equation*}
    The \emph{normal order} $\normord{\ \ }$ of vertex operators is defined as follows:
\begin{gather*}
    \normord{Y(u,z_1)Y(v,z_2)}=Y^-(u,z_1)Y(v,z_2)+(-1)^{|u||v|}Y(v,z_2)Y^+(u,z_1),
    \shortintertext{where}
    Y(u,z)=Y^-(u,z)+Y^+(u,z)=\sum_{n<0}u_nz^{-n-1}+\sum_{n\ge 0}u_nz^{-n-1}.
\end{gather*}
Thus, the spanning elements of $\Theta$ and $\Omega$ arise from normally ordered products (modulo the $L(-1)$-derivative terms). Such vectors may occur in the singular operator products of generators, as happens in nonlinear $W$-algebras (cf. \cite{DK,KRW,KW,KPP}). 
\end{remark}
To obtain linear closure, we shall restrict the relevant singular products, rather than require $\Theta$ and $\Omega$ themselves to vanish. The following lemma illustrates why the latter requirement would be too restrictive when $\fg\ne0$.

\begin{lemma}\label{lem:c=-2}
    Let $V=\bigoplus_{q\in\frac{1}{2}\bZ}V_{(q)}$ be a conformal vertex superalgebra strongly generated by $\fg+U+V_{(2)}$. If the central charge $c_V\neq-2$ and $\Omega=0$, then $\fg=0$.
\end{lemma}
\begin{proof}
    For any $x\in\fg$, assume that $L(1)x=k\bone$ for some $k\in\bC$. By Proposition~\ref{prop:decom},
    \begin{equation*}
        f=x_{-1}\omega-L(-1)\Big(\frac14-\frac1{24}L(-1)L(1)\Big)L(1)(x_{-1}\omega)
        \in\Omega.
    \end{equation*}
    Note that $\omega$ is quasi-primary:
    \begin{equation}\label{w.quasi.prim}
        L(1)\omega=L(1)L(-2)\bone=3L(-1)\bone=0.
    \end{equation}
    Skew-symmetry gives $x_0\omega=0$, and the Borcherds commutator formula gives
    \begin{equation*}
        [L(1),x_{-1}]=[\omega_2,x_{-1}]=x_0+k\Id_V.
    \end{equation*}
    Thus $L(1)(x_{-1}\omega)=k\omega$ and $f=x_{-1}\omega-\frac{k}{4}L(-1)\omega$. If $\Omega=0$, then
    \begin{equation*}
        x_{-1}\omega=\frac{k}{4}L(-1)\omega.
    \end{equation*}
    By skew-symmetry, we have $x_{-1}\omega=L(-2)x-\frac12L(-1)^2x$. This together with \eqref{Vir} shows that
    \begin{equation*}
        L(2)L(-2)x=(4+\tfrac{c_V}{2})x,\ \ L(2)L(-1)^2x=6x,
    \end{equation*}
    where we use $L(-1)L(1)x=0$. Hence $L(2)(x_{-1}\omega)=(1+\frac{c_V}{2})x$. On the other hand, by \eqref{Vir} and \eqref{w.quasi.prim},
    \begin{equation*}
        L(2)L(-1)\omega=3L(1)\omega+L(-1)L(2)\omega=0.
    \end{equation*}
    Thus $(1+\frac{c_V}{2})x=0$, and $c_V\ne-2$ implies $x=0$. This completes the proof.
\end{proof}

\subsection{The commutative associative algebras $\cA$ and $A$}
\hspace{1.5em}Motivated by the preceding decompositions, we first select quasi-primary weight-two vectors whose 0-th products with $V_{(2)}$ have no component in $\Omega$. We then select those elements whose 0-th products with $U$ have no component in $\Theta$ in a similar way.

Let $\cA$ be the linear subspace of $\Gamma$ defined by
\begin{equation}\label{mathcal.A}
    \cA=\{a\in\Gamma \mid a_0v\in L(-1)V_{(2)}, \ \forall\ v\in V_{(2)}\}.
\end{equation}
Equivalently, by Proposition~\ref{prop:decom}, $a\in\cA$ means $a\in\Gamma$ and $\pi_\Omega(a_0v)=0$ for every $v\in V_{(2)}$. For any $a,b\in\cA$, the definition of $\cA$ gives $a_0b\in L(-1)V_{(2)}$. By the injectivity of $L(-1)$ on $V_{(2)}$, we can deduce that there exists a unique $c\in V_{(2)}$ such that $a_0b=L(-1)c$. So we define the operation $\circ$ on $\cA$ by $a\circ b=c$. For the weight-$\frac32$ subspace, denote
\begin{align}
    A&=\{a\in\cA \mid a_0v\in L(-1)U,\ \forall\ v\in U\}\notag\\
        &=\{a\in\Gamma \mid  a_0V_{(2)}\subset L(-1)V_{(2)},\,a_0U\subset L(-1)U\}.\label{A}
\end{align}
In other words, $A$ consists of the elements $a\in\Gamma$ such that $\pi_\Omega(a_0v)=0$ and $\pi_\Theta(a_0u)=0$ for all $v\in V_{(2)}$, $u\in U$.

\begin{theorem}\label{thm:canonical.algebras}
The spaces $\cA$ and $A$ defined in \eqref{mathcal.A} and \eqref{A} have the following properties\textup{:}

\subno{i} For $a,b\in\cA$, we have
\begin{equation*}
    a\circ b=\frac12a_1b,\ \ a_2b=0.
\end{equation*}
The product determined by $a_0b=L(-1)(a\circ b)$ makes $\cA$ a commutative associative algebra with identity $\omega$. 

\subno{ii} The space $A$ is a subalgebra of $\cA$ containing $\omega$.

\subno{iii} The bilinear form $\beta$ defined by $a_3b=\beta(a,b)\bone$ is symmetric and invariant on $\cA$, and its restriction to $A$ has the same properties.
\end{theorem}

We prove Theorem~\ref{thm:canonical.algebras} through the following lemmas.

\begin{lemma}\label{lem:alg.A.o}
    The product $\circ$ is a bilinear operation on $\cA$, and $\omega\in\cA$ is a left identity. Moreover, for $a,b\in\cA$,
    \begin{equation}\label{a2b}
        a\circ b=\frac12a_1b,\ \ a_0b=\frac12L(-1)a_1b,\ \ a_2b=0.
    \end{equation}
\end{lemma}
\begin{proof}
    Let $a,b\in \cA$. Then by the definition of $\cA$, there exists $c\in V_{(2)}$ such that $a_0b=L(-1)c$. From \eqref{comm} and \eqref{der}, we have
    \begin{equation}\label{a_0b_0}
        [a_0,b_0]=(a_0b)_0=\big(L(-1)c\big)_0=0.
    \end{equation}
    As $a\in\Gamma$, skew-symmetry gives
    \begin{align*}
        a_0\omega=-\omega_0a+L(-1)\omega_1a-\frac{1}{2}L(-1)^2\omega_2a=-L(-1)a+2L(-1)a=L(-1)a.
    \end{align*}
    Similarly, we can obtain $b_0\omega=L(-1)b$ for $b\in\Gamma$. Applying \eqref{a_0b_0} to $\omega$ yields
    \begin{align*}
        0=[a_0,b_0]\omega=a_0L(-1)b-b_0L(-1)a=L(-1)(a_0b-b_0a).
    \end{align*}
    The injectivity of $L(-1)$ on $V_{(3)}$ implies
    \begin{equation}\label{zero.product.symm}
        a_0b=b_0a.
    \end{equation}
    Applying $L(1)$ and $L(1)^2$ to both sides of \eqref{zero.product.symm} and using \eqref{quasi.prim}, we get
    \begin{equation}\label{a1b,a2b}
        a_1b=b_1a\ \ \text{and}\ \ a_2b=b_2a,
    \end{equation}
    respectively. On the other hand, by skew-symmetry \eqref{skew.symm}, and noting that $L(-1)^i(b_ia)=0$ for $i\ge 3$, we have
    \begin{gather}
    a_0b=-b_0a-\sum_{i\ge 1}\frac{1}{i!}(-1)^iL(-1)^i(b_ia)=-b_0a+L(-1)b_1a-\frac{1}{2}L(-1)^2b_2a,\label{symm.a0b}\\
        a_2b=-b_2a-\sum_{i\ge1}\frac{1}{i!}(-1)^iL(-1)^i(b_{2+i}a)=-b_2a.\label{symm.a2b}
    \end{gather}
    Combining \eqref{zero.product.symm}--\eqref{symm.a2b}, we can derive that
    \begin{equation*}
        a_2b=0\ \ \text{and}\ \ a_0b=\frac{1}{2}L(-1)a_1b.
    \end{equation*}
    Furthermore, $c=\frac{1}{2}a_1b$. Thus $L(1)c=\frac{1}{2}L(1)a_1b=\frac{1}{2}[L(1),a_1]b=\frac{1}{2}a_2b=0$, and $c\in\Gamma$.

    We next verify that $c_0v\in L(-1)V_{(2)}$ for every $v\in V_{(2)}$. Letting $v\in V_{(2)}$ be an arbitrary element, by \eqref{iter} and \eqref{der}, we have
    \begin{equation}\label{c_0x}
        c_0v=-\big(L(-1)c\big)_1v=-(a_0b)_1v=b_1(a_0v)-a_0(b_1v).
    \end{equation}
    It follows from $b_1v\in V_{(2)}$ and $a\in\cA$ that $a_0(b_1v)\in L(-1)V_{(2)}$. Moreover, since $a\in \cA$ and $v\in V_{(2)}$, there exists $u\in V_{(2)}$ such that $a_0v=L(-1)u$. Then
    \begin{equation*}
        b_1(a_0v)=b_1\big(L(-1)u\big)=[b_1,L(-1)]u+L(-1)b_1u=b_0u+L(-1)b_1u\in L(-1)V_{(2)}.
    \end{equation*}
    It follows from \eqref{c_0x} that $c_0v\in L(-1)V_{(2)}$. As a consequence, $\cA$ is closed under $\circ$.

    The bilinearity of $\circ$ follows from the bilinearity of the map $(a,b)\mapsto a_0b$ and the injectivity of $L(-1)$ on $V_{(2)}$. So $(\cA,\circ)$ is an algebra. Finally, since $\omega\in\Gamma$ and $\omega_0v=L(-1)v$ for every $v\in V_{(2)}$, we have $\omega\in\cA$ and $\omega\circ a=a$. This completes the proof.
\end{proof}

The projection formula in Proposition~\ref{prop:decom} also shows that, for $a,b\in\cA$,
    \begin{equation}\label{c^ab}
        a\circ b=\Big(\frac{1}{4}-\frac{1}{24}L(-1)L(1)\Big)L(1)(a_0b).
    \end{equation}

\begin{lemma}\label{lem:cA.o}
    The algebra $(\cA,\circ)$ defined above is a commutative associative algebra with a symmetric invariant bilinear form $\beta(\cdot,\cdot)$ defined by
    \begin{align*}
        a_3 b=\beta(a,b)\bone,\ \ \forall\ a,b\in \cA.
    \end{align*}
\end{lemma}
\begin{proof}
    Let $a,b\in \cA$. It is direct from \eqref{a2b} and \eqref{a1b,a2b} that
    \begin{equation}\label{comm.o}
        a\circ b=b\circ a=\frac{1}{2}a_1b\ \ \text{and}\ \ a_1b=2(a\circ b)\in \cA.
    \end{equation}
    By \eqref{comm}, \eqref{a2b} and \eqref{comm.o}, we get
    \begin{equation}\label{[a1,b1]}
        [a_1,b_1]=(a_0b)_2+(a_1b)_1=\big(L(-1)(a\circ b)\big)_2+2(a\circ b)_1=-2(a\circ b)_1+2(a\circ b)_1=0.
    \end{equation}
    It implies that $a_1(b_1d)=b_1(a_1d)$ for any $d\in \cA$. Then we can obtain
    \begin{equation}\label{abd}
        a\circ (b\circ d)=\frac{1}{4}a_1b_1d=\frac{1}{4}b_1a_1d=b\circ(a\circ d).
    \end{equation}
    By the commutativity of $\circ$ in \eqref{comm.o} and \eqref{abd}, we deduce that
    \begin{align*}
        a\circ (b\circ d)=b\circ(a\circ d)=b\circ(d\circ a)=d\circ(b\circ a)=d\circ (a\circ b)=(a\circ b)\circ d.
    \end{align*}
    Therefore, $(\cA,\circ)$ is a commutative associative algebra.

    Note that as $a_3b\in V_{(0)}=\bC\bone$, the bilinear form $\beta(\cdot,\cdot)$ on $\cA$ is well defined. By skew-symmetry \eqref{skew.symm} and the fact that $b_{3+i}a\in V_{(-i)}=0$ for $i\ge 1$, we have
    \begin{equation*}
        a_3b=b_3a+\sum_{i\ge1}\frac{1}{i!}(-1)^iL(-1)^i(b_{3+i}a)=b_3a.
    \end{equation*}
    Then $\beta(a,b)=\beta(b,a)$. Applying the Borcherds commutator formula \eqref{comm} to $[a_m,b_n]$ with $m=1$ and $n=3$, and using \eqref{der}, we get
    \begin{align*}
        [a_1,b_3]=(a_0b)_4+(a_1b)_3&=\frac{1}{2}\big(L(-1)(a_1b)\big)_4+(a_1b)_3\\
        &=-2(a_1b)_3+(a_1b)_3\\
        &=-(a_1b)_3.
    \end{align*}
    Hence, for all $d\in \cA$,
    \begin{equation*}
        (a_1b)_3d=-[a_1,b_3]d=b_3(a_1d)-a_1(b_3d)=b_3(a_1d),
    \end{equation*}
    where $b_3d\in V_{(0)}=\bC\bone$ and $a_1\bone=0$. It follows that $\beta( b,a_1d)=\beta(a_1b,d)$, namely
    \begin{align*}
        \beta( b, a\circ d)=\beta( a\circ b,d).
    \end{align*}
    Consequently, we have shown that $\beta(\cdot,\cdot)$ is symmetric and invariant.
\end{proof}

We now verify that $A$ is a subalgebra of $\cA$.

\begin{lemma}\label{lem:A.subalg}
    The space $A$ is a subalgebra of $(\cA,\circ)$ containing $\omega$. In particular, $(A,\circ)$ is a commutative associative algebra with identity $\omega$, and the restriction of $\beta(\cdot,\cdot)$ to $A$ is symmetric and invariant.
\end{lemma}
\begin{proof}
    Let $a,b\in A$ and set $c=a\circ b$. By Lemmas~\ref{lem:alg.A.o} and~\ref{lem:cA.o}, we have $c\in\cA$. For any $v\in U$, a similar calculation as in \eqref{c_0x} gives
    \begin{equation*}
        c_0v=b_1(a_0v)-a_0(b_1v).
    \end{equation*}
    Since $b_1v\in U$ and $a\in A$, the second term on the right-hand side belongs to $L(-1)U$. Write $a_0v=L(-1)w$ with $w\in U$. Then
    \begin{equation*}
        b_1(a_0v)=b_1L(-1)w=b_0w+L(-1)b_1w
        \in L(-1)U,
    \end{equation*}
    because $b\in A$. It shows that $c\in A$. Furthermore, $\omega\in\cA$ and $\omega_0v=L(-1)v$ for all $v\in U$, so $\omega\in A$. The remaining assertions follow by restriction from $\cA$.
\end{proof}

Combining Lemmas~\ref{lem:alg.A.o}--\ref{lem:A.subalg}, we can obtain Theorem~\ref{thm:canonical.algebras}. By a slight abuse of notation, $\circ$ and $\beta(\cdot,\cdot)$ also denote their restrictions to $A$.

\section{Associated algebraic structures}\label{sect:data}
\hspace{1.5em}This section identifies the algebraic structures and compatibility conditions guaranteeing the singular products of the low-weight generators.

We first study the mutual actions of $\fg$ and $\cA$, together with the bilinear forms recording the scalar terms in their singular products. We then restrict the even data to $A$ and determine compatible module structures on $U$. The odd--odd products are described by the bilinear maps. Under the additional closure condition on these products,
the resulting data contributes to all commutation relations among the generators. Finally, we formulate these data and prove that their compatibility conditions are necessary and sufficient for the brackets to define the Lie superalgebra $\cL(\fg,A,U)$. The desired vertex operator superalgebra will be constructed in Section~\ref{sect:VOSA}.

For a conformal vertex superalgebra strongly generated by $V_{(1)}+V_{(\frac32)}+V_{(2)}$, we retain the notation $\fg=V_{(1)}$, $U=V_{(\frac32)}$, and algebras $\cA$ and $A$ from \eqref{mathcal.A} and \eqref{A}, respectively.

\subsection{Algebraic structures associated with $\fg+\cA$}\label{sect:even.generators}
\hspace{1.5em}Let $V=\bigoplus_{n\in\bZ}V_{(n)}$ be a conformal vertex algebra strongly generated by $\fg+\cA$, where $\fg=V_{(1)}$ and $\cA\subset V_{(2)}$ such that $(\cA,\circ)$ is a commutative associative algebra defined in Theorem~\ref{thm:canonical.algebras}. 

It is well known that the weight-one subspace $V_{(1)}$ forms a Lie algebra $\fg$ equipped with a symmetric invariant bilinear form:
\begin{equation}\label{Lie}
    [x,y]=x_0y,\ \ \kappa(x,y)\bone=x_1y,\ \ \forall\ x,y\in\fg.
\end{equation}

For the commutative associative algebra $(\cA,\circ)$, a linear map $D:\cA\rightarrow \cA$ is called a \emph{derivation}, if it satisfies 
\begin{equation*}
    D(a\circ b)=D(a)\circ b+a\circ D(b),\ \ \forall\ a,b\in\cA.
\end{equation*}
The vector space of all such derivations is denoted by $\Der(\cA)$. In addition, $\Der(\cA)$ carries a Lie algebra structure
\begin{align*}
    [D_1,D_2](a)=D_1\big(D_2(a)\big)-D_2\big(D_1(a)\big),\ \ \forall\ D_1,D_2\in\Der(\cA),\, a\in \cA.
\end{align*}

\begin{lemma}\label{lem:g-mod.A}
    The map
    \begin{equation*}
        \rho: \fg\rightarrow \Der(\cA);\ \ \rho(x)a=x_0a,\ \ x\in\fg,\,a\in \cA,
    \end{equation*}
    is a Lie algebra homomorphism.
\end{lemma}
\begin{proof}
    Let $x\in\fg$ and $a\in\cA$. We first show that $x_0a\in \cA$, that is, $x_0a$ is quasi-primary and $(x_0a)_0v\in L(-1)V_{(2)}$ for all $v\in V_{(2)}$. Noting that $a$ is quasi-primary and using \eqref{vac.cre}, \eqref{comm} and \eqref{der}, we have
    \begin{align}\label{x_0a}
        L(1)(x_0a)=[L(1),x_0]a=\big(L(-1)x\big)_2a+2\big(L(0)x\big)_1a+\big(L(1)x\big)_0a=0,
    \end{align}
    namely $x_0a$ is quasi-primary. In addition, for $v\in V_{(2)}$, by \eqref{comm}, we can obtain
    \begin{equation}\label{x0a0v}
        (x_0a)_0v=[x_0,a_0]v=x_0(a_0v)-a_0(x_0v).
    \end{equation}
    On the right-hand side of the equation \eqref{x0a0v}, it is clear that $a_0(x_0v)\in L(-1)V_{(2)}$. For the first term, there exists $u\in V_{(2)}$ such that $a_0v=L(-1)u$. Then by \eqref{der},
    \begin{align*}
        x_0\big(L(-1)u\big)=[x_0,L(-1)]u+L(-1)x_0u=L(-1)x_0u\in L(-1)V_{(2)}.
    \end{align*}
    So $(x_0a)_0v\in L(-1)V_{(2)}$. This together with \eqref{x_0a} proves $x_0a\in \cA$.

    We now verify that $\rho$ is a Lie algebra homomorphism. For any $x,y\in\fg$ and $a\in \cA$, by \eqref{Lie}, we get
    \begin{gather*}
        \rho([x,y])a=\rho(x_0y)a=(x_0y)_0a,\\
        [\rho(x),\rho(y)]a=\rho(x)\rho(y)a-\rho(y)\rho(x)a=x_0(y_0a)-y_0(x_0a).
    \end{gather*}
    It then follows from \eqref{comm} that $\rho([x,y])=[\rho(x),\rho(y)]$.

    Moreover, for any $x\in\fg$, $a,b\in \cA$, by \eqref{comm} and \eqref{comm.o}, we can derive 
    \begin{align}
        (x_0a)\circ b+a\circ(x_0b)&=\frac{1}{2}(x_0a)_1b+\frac{1}{2}a_1(x_0b)=\frac{1}{2}[x_0,a_1]b+\frac{1}{2}a_1(x_0b)\notag\\
        &=\frac{1}{2}x_0(a_1b)=x_0(a\circ b).\label{x_0.der}
    \end{align}
    It implies that $\rho(x)=x_0$ acts as a derivation on $(\cA,\circ)$. This completes the proof.
\end{proof}

\begin{lemma}\label{lem:A-mod.g}
    Define a left action and a right action of the commutative associative algebra $(\cA,\circ)$ on $\fg$ by
    \begin{equation*}
        a\cdot x:=a_1x,\ \ x\cdot a:=x_1a,\ \ \text{for}\ a\in \cA,\, x\in\fg.
    \end{equation*}
    Then $\fg$ is a symmetric $\cA$-bimodule. More explicitly, the identity $\omega$ of $\cA$ acts as the identity on $\fg$ and the following equations hold\textup{:}
    \begin{align}
        a\cdot x &= x\cdot a,\label{a.x.comm}\\
        (a\circ b)\cdot x &= a\cdot (b\cdot x),\label{a.x.asso}\\
        x\cdot (a\circ b) &= (x\cdot a)\cdot b,\label{a.x.asso.r}
    \end{align}
    for all $a,b\in \cA$ and $x\in\fg$.
\end{lemma}
\begin{proof}
    Let $a,b\in \cA$ and $x\in\fg$. Since $\wt(a_1x)=\wt(x_1a)=1$, we have $a\cdot x\in\fg$ and $x\cdot a\in\fg$. By the skew-symmetry \eqref{skew.symm} and the fact that $L(-1)^i(x_ia)=0$ for $i\ge 2$, we obtain $a_1x=x_1a$, which implies equation \eqref{a.x.comm}. Then by \eqref{L(0)},
    \begin{equation*}
        x\cdot\omega=\omega\cdot x=\omega_1x=L(0)x=x,
    \end{equation*}
    which verifies that $\omega$ acts as the identity. Since $(\cA,\circ)$ is commutative, if the left associativity \eqref{a.x.asso} holds, the right associativity \eqref{a.x.asso.r} follows naturally. Thus, it suffices to prove \eqref{a.x.asso}.

    Since $x_2a\in V_{(0)}=\bC\bone$, the skew-symmetry \eqref{skew.symm} and the truncation $L(-1)^i(x_ia)=0$ ($i\ge 2$) yield
    \begin{equation}\label{a_0x}
        a_0x=-x_0a+L(-1)x_1a.
    \end{equation}
    Noticing that $x_0a\in \cA$ yields $b_2(x_0a)=0$ by \eqref{a2b}. Applying $b_2$ to \eqref{a_0x} and using $b_2x_1a\in V_{(0)}=\bC\bone$, we get
    \begin{equation}\label{b_2a_0x}
        b_2(a_0x)=-b_2(x_0a)+b_2L(-1)x_1a=[b_2,L(-1)]x_1a+L(-1)b_2x_1a=2b_1x_1a.
    \end{equation}
     Evaluating the commutator $[a_0,b_2]$ on $x$, by \eqref{comm} and \eqref{a2b}, we deduce that
    \begin{equation}\label{a0b2.1}
        [a_0,b_2]x=(a_0b)_2x=\big(L(-1)(a\circ b)\big)_2x=-2(a\circ b)_1x.
    \end{equation}
    Further expanding the left-hand side and using \eqref{b_2a_0x} implies
    \begin{equation}\label{a0b2.2}
        [a_0,b_2]x=a_0(b_2x)-b_2(a_0x)=-2b_1a_1x,
    \end{equation}
    where $b_2x\in\bC\bone$ and $a_0\bone=0$. Moreover, by \eqref{[a1,b1]}, we have $a_1(b_1x)=b_1(a_1x)$. Equating \eqref{a0b2.1} and \eqref{a0b2.2} gives
    \begin{equation*}
        (a\circ b)_1x=b_1x_1a=b_1a_1x=a_1b_1x.
    \end{equation*}
    Therefore, $(a\circ b)\cdot x=a\cdot (b\cdot x)$, which completes the proof.
\end{proof}

For convenience, we also write $ax=xa$ for $a\cdot x=x\cdot a$. This notation refers to the $\cA$-module structure on $\fg$. The action of $\fg$ on $\cA$ is always written as $\rho(x)a=x_0a$.

\begin{lemma}\label{lem:Lie-Rine}
    Equipped with the $\fg$-module and symmetric $\cA$-bimodule structures defined in Lemmas~\ref{lem:g-mod.A} and \ref{lem:A-mod.g}, respectively, $\cA$ and $\fg$ satisfy the following compatibility conditions\textup{:}
    \begin{gather}
        [x,a\cdot y]=a\cdot [x,y]+\rho(x)a\cdot y,\label{compt.1}\\
        \rho(a\cdot x)b-a\circ\rho(x)b=\rho(b\cdot x)a-b\circ \rho(x)a,\label{compt.2}
    \end{gather}
    where $a,b\in \cA$ and $x,y\in \fg$.
\end{lemma}
\begin{proof}
    By Lemmas~\ref{lem:g-mod.A} and \ref{lem:A-mod.g}, $\cA$ and $\fg$ carry mutual module structures. Let $a,b\in \cA$ and $x,y\in\fg$. The Borcherds commutator formula \eqref{comm} gives $[x_0,a_1]=(x_0a)_1$, which yields
    \begin{align*}
        \rho(x)a\cdot y=(x_0a)_1y=[x_0,a_1]y=x_0(a_1y)-a_1(x_0y)=[x,a\cdot y]-a\cdot[x,y].
    \end{align*}
    This establishes \eqref{compt.1}.

    Since $b_1x\in\fg$ and $a\circ b\in \cA$, equations \eqref{a.x.comm} and \eqref{a_0x} imply $x_1(a\circ b)=(a\circ b)_1x$ and $a_0(b_1x)=-(b_1x)_0a+L(-1)a_1(b_1x)$, respectively. Applying $[x_1,a_0]=x_1a_0-a_0x_1$ to $b$, the right-hand side gives
    \begin{align}
        x_1(a_0b)-a_0(x_1b)&=x_1\big(L(-1)(a\circ b)\big)-a_0(b_1x)\notag\\
        &= [x_1,L(-1)](a\circ b)+L(-1)x_1(a\circ b)+(b_1x)_0a-L(-1)a_1(b_1x)\notag\\
        &=x_0(a\circ b)+L(-1)x_1(a\circ b)+(b_1x)_0a-L(-1)a_1(b_1x)\notag\\
        &=x_0(a\circ b)+(b_1x)_0a+L(-1)\big((a\circ b)_1x-a_1(b_1x)\big)\notag\\
        &=(x_0a)\circ b+a\circ (x_0b)+(b_1x)_0a,\label{[x_1,a_0].1}
    \end{align}
    where the latter steps are from the associativity \eqref{a.x.asso} and the derivation property \eqref{x_0.der}. On the other hand, applying the Borcherds commutator formula to $[x_1, a_0]$ and using \eqref{comm.o} and \eqref{a.x.comm} yields
     \begin{equation}\label{[x_1,a_0].2}
         [x_1,a_0]b=(x_0a)_1b+(x_1a)_0b = 2(x_0a)\circ b+(a_1x)_0b.
     \end{equation}
     Equating \eqref{[x_1,a_0].1} and \eqref{[x_1,a_0].2}, we obtain
     \begin{equation*}
         a\circ (x_0b)+(b_1x)_0a=(x_0a)\circ b+(a_1x)_0b.
     \end{equation*}
     Rearranging the terms directly leads to
     \begin{equation*}
         (a_1x)_0b-a\circ(x_0b)=(b_1x)_0a-b\circ(x_0a).
     \end{equation*}
    This verifies the second compatibility condition \eqref{compt.2}.
\end{proof}

\begin{remark}\label{rmk:Lie-Rine}\textbf{(Connection with Lie--Rinehart algebras).}
The structures in Lemma~\ref{lem:Lie-Rine} are closely related to \emph{Lie--Rinehart algebras} \cite{Hue,Rine}. Indeed, the symmetric $\mathcal A$-bimodule structure makes $\fg$ an $\cA$-module, and the Lie algebra homomorphism $\rho:\fg\rightarrow\Der(\cA)$ satisfies \eqref{compt.1}. Hence $(\fg,\cA)$ is a Lie--Rinehart algebra with \emph{anchor} $\rho$ if and only if $\rho$ is $\cA$-linear, that is,
\begin{equation}\label{A.linear}
    \rho(a\cdot x)b=a\circ\big(\rho(x)b\big),\ \ \forall\ a,b\in\mathcal A,\ x\in\mathfrak g.
\end{equation}
Setting
\begin{equation*}
    D_x(a,b):=\rho(a\cdot x)b-a\circ\big(\rho(x)b\big),
\end{equation*}
condition \eqref{compt.2} implies that $D_x(a,b)=D_x(b,a)$, which means that $D_x(\cdot,\cdot)$ is a symmetric bilinear map from $\cA\times \cA$ to $\cA$, while $\mathcal A$-linearity \eqref{A.linear} requires $D_x(a,b)=0$. In particular, the case $\rho=0$ yields a
Lie--Rinehart algebra with zero anchor.
\end{remark}

To complete the description of the relation between $\fg$ and $\cA$, we introduce the bilinear map $\theta:\fg\times \cA\rightarrow\bC$ defined by
\begin{equation*}
    x_2a=\theta(x,a)\bone,\ \ \forall\ x\in\fg,\, a\in \cA.
\end{equation*}
It follows directly from the skew-symmetry $x_2a=-a_2x$ that
\begin{equation}\label{a_2x}
    a_2x=-\theta(x,a)\bone.
\end{equation}

\begin{lemma}\label{lem:theta}
    For any $x,y\in\fg$ and $a,b\in \cA$, the bilinear map $\theta$ satisfies
    \begin{gather}
        \theta(ax,b)=\theta(x,a\circ b)=\theta(bx,a),\label{theta.1}\\
        \theta\big(x,\rho(y)a\big)=\theta([x,y],a)=-\theta\big(y,\rho(x)a\big),\label{theta.2}
    \end{gather}
    and the compatibility conditions\textup{:}
    \begin{gather}
        \beta\big(\rho(x)a,b\big)=0=\beta\big(a,\rho(x)b\big),\label{theta.3}\\
        \kappa(x,ay)-\kappa(ax,y)=\theta([x,y],a).\label{theta.4}
    \end{gather}
\end{lemma}
\begin{proof}
    Let $x,y\in\fg$ and $a,b\in \cA$. Using the Borcherds commutator formula, we obtain $[a_1,b_2]x=(a_0b)_3x+(a_1b)_2x$. Further by \eqref{comm.o} and \eqref{a_2x}, we get
    \begin{gather*}
        a_1(b_2x)-b_2(a_1x)=-a_1\theta(x,b)\bone-b_2(a\cdot x)=\theta(ax,b)\bone,\\
        (a_0b)_3x+(a_1b)_2x=\big(L(-1)(a\circ b)\big)_3x+2(a\circ b)_2x=-(a\circ b)_2x=\theta(x,a\circ b)\bone.
    \end{gather*}
    Equating the two equations, we have $\theta(ax,b)=\theta(x,a\circ b)$. Since $\cA$ is commutative, it follows that
    \begin{equation*}
        \theta(ax,b)=\theta(x,a\circ b)=\theta(x,b\circ a)=\theta(bx,a).
    \end{equation*}
    Similarly, by $[x_2,y_0]a=(x_0y)_2a+2(x_1y)_1a+(x_2y)_0a=(x_0y)_2a$, we can get $\theta\big(x,\rho(y)a\big)=\theta([x,y],a)$, and thus
    \begin{equation*}
        \theta\big(x,\rho(y)a\big)=\theta([x,y],a)=-\theta([y,x],a)=-\theta\big(y,\rho(x)a\big).
    \end{equation*}
    Therefore, we obtain \eqref{theta.1} and \eqref{theta.2}.

    Since $b_0x=-\rho(x)b+L(-1)(bx)$ and $a_3x=0$, we can deduce that
    \begin{gather*}
        [a_3,b_0]x=a_3(b_0x)=-a_3\big(\rho(x)b\big)+[a_3,L(-1)](bx)=-\beta\big(a,\rho(x)b\big)\bone-3\theta(bx,a)\bone,\\
        (a_0b)_3x+3(a_1b)_2x+3(a_2b)_1x+(a_3b)_0x=(a_0b)_3x+3(a_1b)_2x=-3\theta(x,a\circ b)\bone.
    \end{gather*}
    Comparing the two expressions gives
    \begin{equation*}
        3\theta(x,a\circ b)-3\theta(bx,a)=\beta\big(a,\rho(x)b\big).
    \end{equation*}
    This together with \eqref{theta.1} shows that $\beta\big(a,\rho(x)b\big)=0$. Due to the fact that $\beta(\cdot,\cdot)$ is invariant, the equation \eqref{theta.3} holds.

    Finally, from $[x_1,y_1]a=(x_0y)_2a+(x_1y)_1a$, we arrive at \eqref{theta.4}, and the lemma follows.
\end{proof}

Combining the structures obtained above, we associate $\fg+\cA$ to the datum
\begin{equation}\label{datum.even}
    \cQ_{ev}(V)=\big((\fg,[\cdot,\cdot],\kappa),(\cA,\circ,\beta),(\rho,\cdot,\theta)\big).
\end{equation}
The first two components consist of the Lie algebra $\fg$ and the commutative associative algebra $\cA$, together with their symmetric invariant bilinear forms. The third consists of the action $\rho$ of $\fg$ on $\cA$, the symmetric $\cA$-bimodule structure on $\fg$, and the bilinear map $\theta$. Based on these, we have the following result.

\begin{proposition}\label{prop:even.mode.relations}
    Let $V=\bigoplus_{n\in\bZ}V_{(n)}$ be a conformal vertex algebra strongly generated by $\fg+\cA$, where $\fg=V_{(1)}$ and $\cA\subset V_{(2)}$ is defined by \eqref{mathcal.A}. Then 
    \begin{align}
        [x_m,y_n]&=[x,y]_{m+n}+m\kappa(x,y)\delta_{m+n,0},\label{even.mode.gg}\\
        [a_m,b_n]&=(m-n)(a\circ b)_{m+n-1}
            +\binom{m}{3}\beta(a,b)\delta_{m+n-2,0},\label{even.mode.AA}\\
        [x_m,a_n]&=\big(\rho(x)a\big)_{m+n}+m(ax)_{m+n-1}
            +\binom{m}{2}\theta(x,a)\delta_{m+n-1,0},\label{even.mode.gA}
    \end{align}
    for $x,y\in\fg$, $a,b\in\cA$ and $m,n\in\bZ$. Here the operations are those of the datum \eqref{datum.even}, satisfying \eqref{compt.1}, \eqref{compt.2} and \eqref{theta.1}--\eqref{theta.4}.
\end{proposition}
\begin{proof}
     Let $x,y\in\fg$ and $a,b\in\cA$. Since $x_iy\in V_{(1-i)}=0$ for $i\ge 2$, applying the Borcherds commutator formula \eqref{comm} to $[x_m,y_n]$ immediately yields
    \begin{equation*}
        [x_m,y_n]=(x_0y)_{m+n}+m(x_1y)_{m+n-1}=[x,y]_{m+n}+m\kappa(x,y)\delta_{m+n,0}.
    \end{equation*}
    Noticing that $a_ib\in V_{(3-i)}=0$ for $i\ge4$ and $a_2b=0$ from \eqref{a2b}, then by Lemma~\ref{lem:cA.o}, we have
    \begin{align*}
        [a_m,b_n]&=(a_0b)_{m+n}+m(a_1b)_{m+n-1}+\binom{m}{3}(a_3b)_{m+n-3}\\
        &=\big(L(-1)(a\circ b)\big)_{m+n}+2m(a\circ b)_{m+n-1}+\binom{m}{3}\beta(a,b)\bone_{m+n-3}\\
        &=(m-n)(a\circ b)_{m+n-1}+\binom{m}{3}\beta(a,b)\delta_{m+n-2,0}.
    \end{align*}
    Since $x_ia\in V_{(2-i)}=0$ for $i\ge3$, then by Lemmas~\ref{lem:g-mod.A}--\ref{lem:theta}, we can obtain
    \begin{align*}
        [x_m,a_n]&=(x_0a)_{m+n}+m(x_1a)_{m+n-1}+\binom{m}{2}(x_2a)_{m+n-2}\\
        &=\big(\rho(x)a\big)_{m+n}+m(xa)_{m+n-1}+\binom{m}{2}\theta(x,a)\bone_{m+n-2}\\
        &=\big(\rho(x)a\big)_{m+n}+m(ax)_{m+n-1}+\binom{m}{2}\theta(x,a)\delta_{m+n-1,0}.
    \end{align*}
    The compatibility relations have been established above. This completes the proof.
\end{proof}

When $\fg=0$, we have $\cA=V_{(2)}$. Substituting $a\circ b=\frac12a_1b$ into \eqref{even.mode.AA}, we recover the commutation relation in Proposition~\ref{prop:Lam}.

\subsection{Algebraic structures associated with $\fg+A+U$}\label{sect:odd.generators}

\hspace{1.5em}We now consider a conformal vertex superalgebra $V=\bigoplus_{q\in\frac{1}{2}\bZ}V_{(q)}$ satisfying the condition:
\begin{equation}\label{assumption}
    V\ \text{is strongly generated by}\ \fg+A+U,
\end{equation}
where $A$ is defined in \eqref{A}, that is $A=\{a\in\Gamma \mid  a_0V_{(2)}\subset L(-1)V_{(2)},\, a_0U\subset L(-1)U\}$. The structural identities established in Subsection~\ref{sect:even.generators} also hold in the present setting, since their proofs use only the standing low-weight assumptions and the definition of $\cA$.

We first restrict the even datum to $A$ and then equip $U$ with compatible $\fg$-module and $A$-module structures. The products $u_2v$ and $u_1v$ define the maps $\tau$ and $\eta$. Finally, we separate the quasi-primary component $\sigma(u,v)$ of $u_0v$. The assumption $\sigma(U,U)\subset A$ is needed to ensure the linear closure.  

\begin{lemma}\label{lem:g-stable.A}
The subalgebra $A$ is stable under the action of $\fg$ on $\cA$. Consequently, $\rho$ restricts to a Lie algebra homomorphism $\rho:\fg\rightarrow\Der(A)$.
\end{lemma}
\begin{proof}
    Let $x\in\fg$ and $a\in A$. Lemma~\ref{lem:g-mod.A} gives $x_0a\in\cA$. For $v\in U$, write $a_0v=L(-1)w$ with $w\in U$. Then
    \begin{align*}
        (x_0a)_0v=x_0(a_0v)-a_0(x_0v)=L(-1)x_0w-a_0(x_0v)\in L(-1)U,
    \end{align*}
since $x_0v,x_0w\in U$.
\end{proof}

Combining Lemmas~\ref{lem:A.subalg} and~\ref{lem:g-stable.A}, we can restrict \eqref{datum.even} to
\begin{equation}\label{datum.even.A}
    \big((\fg,[\cdot,\cdot],\kappa),(A,\circ,\beta),(\rho,\cdot,\theta)\big),
\end{equation}
where, by a slight abuse of notation, $\rho$, $\cdot$ and $\theta$ also denote their respective restrictions to $A$.

Now we show that $U$ is endowed with $\fg$-module and $A$-module structures satisfying certain compatibility conditions. For convenience, a space equipped with such compatible structures will be referred to as a $(\fg,A)$-module.

\begin{lemma}
    Define a left action of the Lie algebra $\fg$ on $U$ by
    \begin{equation*}
        xv:=x_0v,\ \ \text{for}\ x\in\fg,\, v\in U.
    \end{equation*}
    Then $U$ is a $\fg$-module, namely $[x,y]v=x(yv)-y(xv)$ for any $x,y\in\fg$, $v\in U$.
\end{lemma}
\begin{proof}
    Let $x,y\in\fg$ and $v\in U$. By \eqref{wt}, we have $x_0v\in U$. Applying $[x_0,y_0]=(x_0y)_0$ to $v$ shows that
    \begin{equation*}
        x(yv)-y(xv)=x_0(y_0v)-y_0(x_0v)=(x_0y)_0v=[x,y]_0v=[x,y]v.
    \end{equation*}
    Hence $U$ is a $\fg$-module.
\end{proof}

By the definition \eqref{A}, $a_0v\in L(-1)U$ for $a\in A$ and $v\in U$. Thus Proposition~\ref{prop:decom} gives
\begin{equation}\label{a_1v}
    a_0v=\frac{1}{3}L(-1)L(1)(a_0v)=\frac{1}{3}L(-1)\big([L(1),a_0]v+a_0L(1)v\big)=\frac{2}{3}L(-1)a_1v,
\end{equation}
where we use the fact that $L(1)A=L(1)U=0$ and \eqref{quasi.prim}.

\begin{lemma}
    Define a left action of the commutative associative algebra $(A,\circ)$ on $U$ by
    \begin{equation*}
        av:=\frac{2}{3}a_1v,\ \ \text{for}\ a\in A,\, v\in U.
    \end{equation*}
    Then the identity $\omega$ of $A$ acts as the identity on $U$ and
    \begin{equation*}
        a(bv)=(a\circ b)v=b(av),\ \ \forall\ a,b\in A,\, v\in U.
    \end{equation*}
\end{lemma}
\begin{proof}
    Let $a,b\in A$ and $v\in U$. It is clear that $a_1v\in U$. Consider the equality $[a_0,b_2]v=(a_0b)_2v$. As $b_2v\in V_{(\frac{1}{2})}$, the left-hand side gives
    \begin{align*}
        a_0(b_2v)-b_2(a_0v)=-\frac{2}{3}b_2L(-1)a_1v=-\frac{2}{3}\big([b_2,L(-1)]a_1v+L(-1)b_2a_1v\big)=-\frac{4}{3}b_1a_1v,
    \end{align*}
    where we use \eqref{a_1v} and the fact that $b_2a_1v\in V_{(\frac{1}{2})}=0$. On the other hand, we have
    \begin{align*}
        (a_0b)_2v=\big(L(-1)(a\circ b)\big)_2v=-2(a\circ b)_1v=-3(a\circ b)v.
    \end{align*}
    Equating the above two equations, we obtain
    \begin{align*}
        (a\circ b)v=\frac{4}{9}b_1a_1v=\frac{2}{3}b_1(av)=b(av)=a(bv),
    \end{align*}
    where the last equality follows by interchanging $a$ and $b$ and using $a\circ b=b\circ a$. Moreover, by the action of $A$ on $U$,
    \begin{equation*}
        \omega v=\frac{2}{3}\omega_1v=\frac{2}{3}L(0)v=v.
    \end{equation*}
    This completes the proof.
\end{proof}

\begin{proposition}
    Let $V=\bigoplus_{q\in\frac{1}{2}\bZ}V_{(q)}$ be a conformal vertex superalgebra under the assumption \eqref{assumption}. Then the $\fg$-module and $A$-module structures on $U$ satisfy the following compatibility condition\textup{:}
    \begin{align}\label{compt.xav}
        (xa)v=x(av)-\frac{3}{2}\big(\rho(x)a\big)v=a(xv)-\frac{1}{2}\big(\rho(x)a\big)v,
    \end{align}
    for all $x\in\fg$, $a\in A$ and $v\in U$; and thus make $U$ a $(\fg,A)$-module.
\end{proposition}
\begin{proof}
    Let $x\in\fg$, $a\in A$ and $v\in U$. Applying the Borcherds commutator formula to $[x_1,a_0]v$, we obtain
    \begin{equation*}
        [x_1,a_0]v=(x_0a)_1v+(x_1a)_0v=\frac{3}{2}\big(\rho(x)a\big)v+(xa)v.
    \end{equation*}
    For the left-hand side, since $\wt(x_1v)=\wt(x_1a_1v)=\frac{1}{2}$ implies $x_1v=x_1a_1v=0$, we use \eqref{a_1v} to get
    \begin{align*}
        x_1(a_0v)-a_0(x_1v)=\frac{2}{3}x_1\big(L(-1)a_1v\big)=\frac{2}{3}[x_1,L(-1)]a_1v=\frac{2}{3}x_0(a_1v)=x(av).
    \end{align*}
    Equating the above two equations yields
    \begin{equation}\label{xav}
        x(av)-\frac{3}{2}\big(\rho(x)a\big)v=(xa)v.
    \end{equation}

    Similarly, by $[x_0,a_1]v=(x_0a)_1v=\frac{3}{2}\big(\rho(x)a\big)v$, we have
    \begin{equation*}
        x_0(a_1v)-a_1(x_0v)=\frac{3}{2}x_0(av)-a_1(xv)=\frac{3}{2}x(av)-\frac{3}{2}a(xv),
    \end{equation*}
    which implies that
    \begin{equation}\label{compt.Lie-Rine}
        x(av)-a(xv)=\big(\rho(x)a\big)v.
    \end{equation}
    Substituting it into \eqref{xav} gives
    \begin{equation*}
        a(xv)-\frac{1}{2}\big(\rho(x)a\big)v=(xa)v.
    \end{equation*}
    Combining these results, we can obtain the equation \eqref{compt.xav}.
\end{proof}

\begin{remark}\label{rmk:weak.rep}\textbf{(Connection with weak representations of Lie--Rinehart algebras).}
Suppose that $(\fg,A)$ is a Lie--Rinehart algebra with anchor $\rho$. Then every $(\fg,A)$-module $U$ is a \emph{weak representation} of this Lie--Rinehart algebra in the sense of \cite{PSTZ}. Indeed, the two actions satisfy $x(au)-a(xu)=\big(\rho(x)a\big)u$ by \eqref{compt.Lie-Rine}, for $x\in\fg$, $a\in A$ and $u\in U$. Our compatibility condition \eqref{compt.xav} further measures the failure of the $\fg$-action to be $A$-linear:
\begin{equation*}
    (ax)u-a(xu)=-\frac12\big(\rho(x)a\big)u.
\end{equation*}
Consequently, the given actions define a Lie--Rinehart module structure on $U$ if and only if 
\begin{equation*}
    \big(\rho(x)a\big)u=0.
\end{equation*}
In particular, this holds when $\rho=0$.
\end{remark}

Having described the actions of $\fg$ and $A$ on $U$, we now turn to the singular products between the elements of $U$. For $u,v\in U$, weight considerations give 
\begin{equation}\label{trun.uv}
    u_0v\in V_{(2)},\ u_1v\in V_{(1)}=\fg,\ u_2v\in V_{(0)}=\bC\bone,\ \text{and}\ u_nv=0\ \text{for}\ n\ge3.
\end{equation}
We equip $U$ with a symmetric $\bC$-valued bilinear form $\tau$:
\begin{equation}\label{tau}
    \tau: U\times U\rightarrow\bC;\ \ u_2v=\tau(u,v)\bone,\ \ \tau(u,v)=\tau(v,u),
\end{equation}
and a skew-symmetric $\fg$-valued bilinear form $\eta$:
\begin{equation}\label{eta}
    \eta: U\times U\rightarrow \fg;\ \ \eta(u,v)=u_1v,\ \ \eta(u,v)=-\eta(v,u),
\end{equation}
where $u,v\in U$. The symmetry of $\tau$ and the skew-symmetry of $\eta$
follow from skew-symmetry in $V$ and $L(-1)\bone=0$.

\begin{lemma}\label{lem:uv.part.1}
    The bilinear forms $\tau$ and $\eta$ on $U$ satisfy the following conditions\textup{:}
    \begin{align}
        &[x,\eta(u,v)]=\eta(xu,v)+\eta(u,xv),\label{compt.uv.3}\\
        &a\cdot\eta(u,v)=\eta(au,v)=\eta(u,av),\label{compt.uv.4}\\
        &\kappa\big(x,\eta(u,v)\big)=\tau(xu,v),\label{compt.uv.5}\\
        &\theta\big(\eta(u,v),a\big)=0,\label{compt.uv.6}
    \end{align}
    for $x\in\fg$, $a\in A$ and $u,v\in U$.
\end{lemma}
\begin{proof}
    Let $x\in\fg$, $a\in A$ and $u,v\in U$. By the Borcherds commutator formula, we have $x_0(u_1v)-u_1(x_0v)=(x_0u)_1v$. Since $x_0u=xu$, $x_0v=xv$ and $u_1v=\eta(u,v)\in\fg$, this identity becomes $[x,\eta(u,v)]-\eta(u,xv)=\eta(xu,v)$, which proves \eqref{compt.uv.3}.

    Considering the equality $[a_1,u_1]v=(a_0u)_2v+(a_1u)_1v$, by \eqref{der} and \eqref{a_1v}, we have
    \begin{gather*}
        a_1(u_1v)-u_1(a_1v)=a_1\eta(u,v)-\frac{3}{2}u_1(av)=a\eta(u,v)-\frac{3}{2}\eta(u,av),\\
        (a_0u)_2v+(a_1u)_1v=\frac{2}{3}\big(L(-1)a_1u\big)_2v+(a_1u)_1v=-\frac{1}{3}\eta(a_1u,v)=-\frac{1}{2}\eta(au,v).
    \end{gather*}
    Equating the above two equations gives $a\cdot \eta(u,v)=\frac{3}{2}\eta(u,av)-\frac{1}{2}\eta(au,v)$. Then interchanging $u$ and $v$, by skew-symmetry of $\eta$, we obtain $a\cdot\eta(u,v)=\frac{3}{2}\eta(au,v)-\frac{1}{2}\eta(u,av)$. Comparing these two identities yields $\eta(au,v)=\eta(u,av)$. Substituting this equation back into either identity, we can get \eqref{compt.uv.4}.

    Finally, since $x_1u=x_1v=0$, the Borcherds commutator formula yields
    \begin{equation*}
        \kappa\big(x,\eta(u,v)\big)\bone=x_1(u_1v)=[x_1,u_1]v=(x_0u)_2v=\tau(xu,v)\bone,
    \end{equation*}
    and \eqref{compt.uv.5} holds. Moreover, since $a_2u=a_2v=0$, equations  \eqref{comm}, \eqref{der} and \eqref{a_1v} give
    \begin{align*}
        -\theta\big(\eta(u,v),a\big)\bone=a_2(u_1v)=[a_2,u_1]v=(a_0u)_3v+2(a_1u)_2v=0.
    \end{align*}
    Thus \eqref{compt.uv.6} follows, completing the proof.
\end{proof}

The preceding results establish the linear closure of the singular products among the elements of $\fg+A$, as well as those between $\fg+A$ and $U$. By \eqref{trun.uv} and Lemma~\ref{lem:uv.part.1}, the only remaining closure condition is that $u_0v\in V_{(2)}$. Recall that Proposition~\ref{prop:decom} provides a canonical way to separate its weight-one derivative part from its quasi-primary part. We have the following proposition.

\begin{proposition}\label{prop:sigma}
    For $u,v\in U$, define
    \begin{equation}\label{def.sigma}
        \sigma(u,v)=u_0v-\frac12L(-1)(u_1v).
    \end{equation}
    Then $\sigma:U\times U\rightarrow\Gamma$ is a symmetric bilinear map, and
    \begin{equation}\label{decom.uv}
        u_0v=\sigma(u,v)+\frac12L(-1)\eta(u,v),
    \end{equation}
    where $\eta(u,v)=u_1v$ is the $\fg$-valued bilinear form defined in \eqref{eta}.
\end{proposition}
\begin{proof}
    It follows from \eqref{quasi.prim} that
    \begin{equation*}
        L(1)(u_0v)=[L(1),u_0]v+u_0\big(L(1)v\big)=u_1v,\ \ \forall\ u,v\in U.
    \end{equation*}
    Therefore, by \eqref{proj.Gamma},
    \begin{equation*}
        \pi_\Gamma(u_0v)=u_0v-\frac12L(-1)L(1)(u_0v)=\sigma(u,v)\in\Gamma.
    \end{equation*}
    The bilinearity is immediate. As $u,v$ are odd, skew-symmetry and \eqref{trun.uv} imply
    \begin{equation*}
        u_0v=v_0u-L(-1)(v_1u)+\frac12L(-1)^2(v_2u)=v_0u+L(-1)(u_1v),
    \end{equation*}
    where $v_2u\in\bC\bone$ and $v_1u=-u_1v$. Therefore,
    \begin{equation}\label{symm.sigma}
        \sigma(u,v)=v_0u+\frac12L(-1)(u_1v)=v_0u-\frac12L(-1)(v_1u)=\sigma(v,u).
    \end{equation}
    Equation \eqref{decom.uv} follows from \eqref{eta} and \eqref{def.sigma}.
\end{proof}

To ensure that the odd--odd singular products also close linearly on $\fg+A+U$, we impose the additional condition
\begin{equation}\label{linear.closure}
    \sigma(u,v)\in A,\ \ \forall\ u,v\in U.
\end{equation}
By \eqref{def.sigma}, this is equivalent to
\begin{equation*}
    u_0v\in A\oplus L(-1)\fg,\ \ \forall\ u,v\in U.
\end{equation*}

\begin{remark}\label{rmk:linear.closure}\textbf{(Connection with minimal $W$-algebras)}
Without the condition \eqref{linear.closure}, the quasi-primary component $\sigma(u,v)$ need not belong to $A$ and may involve quadratic current terms of the form $x_{-1}y$ with $x,y\in\fg$. However, several structural and compatibility identities remain valid. In particular, $\sigma$ is still symmetric. Terms of this form occur naturally in minimal $W$-algebras obtained by quantum Hamiltonian reduction \cite{KRW,KW}. 

More explicitly, let $\fs=\fs_\ep\oplus\fs_\op$ be a basic Lie superalgebra equipped with a fixed nondegenerate even invariant supersymmetric bilinear form. Let $x\in\fs_\ep$ be $ad$-diagonalizable with eigenvalues in $\frac12\bZ$, and write $\fs=\bigoplus_{j\in\frac12\bZ}\fs_j$ for the corresponding grading. Choose a nonzero $f\in\fs_{-1}\cap\fs_\ep$ such that the centralizer $\fs^f:=\{a\in\fs\mid[a,f]=0\}$ is contained in $\bigoplus_{j\leq0}\fs_j$. The datum $(\fs,x,f)$ is called \emph{minimal} if
\begin{equation*}
    \fs_j=0\ \ \text{for }|j|>1,\ \  \fs_{-1}=\bC f.
\end{equation*}
For $k\in\bC\setminus\{-h^\vee\}$, where $h^\vee$ is the dual Coxeter number of $\fs$, the quantum Hamiltonian reduction of $V^k(\fs)$ associated with $(x,f)$ is called the \emph{universal minimal $W$-algebra}. This algebra is strongly generated by a Virasoro field $L$ and primary fields $J^a$ and $G^u$. The weights of $J^a$ and $G^u$ are $1$ and $\frac32$, respectively \cite{KW}. The strong generation and weight counting allow quadratic current terms only in the zeroth products $G^u_0G^v$, which have weight $2$. 

Extending the present framework to algebras that need not satisfy
\eqref{linear.closure} will require substantial further work.
\end{remark}

\begin{lemma}\label{lem:sigma.compat}
Under the assumption \eqref{linear.closure}, for $x\in\fg$, $a\in A$ and $u,v,w\in U$, the symmetric bilinear map $\sigma:U\times U\rightarrow A$ satisfies the following compatibility conditions.

    \subno{i} Compatibility with $\fg$\textup{:}
    \begin{align}
    \rho(x)\sigma(u,v)&=\sigma(xu,v)+\sigma(u,xv),\label{sigma.g}\\
    2\sigma(u,v)\cdot x&=\eta(xu,v)-\eta(u,xv),\label{sigma.x}\\
    \theta\big(x,\sigma(u,v)\big)&=0.\label{sigma.theta}
    \end{align}
    \subno{ii} Compatibility with $A$ and the bilinear forms\textup{:}
    \begin{align}
    2a\circ\sigma(u,v)&=\sigma(au,v)+\sigma(u,av),\label{sigma.A}\\
    \rho\big(\eta(u,v)\big)a&=\sigma(au,v)-\sigma(u,av),\label{sigma.rho}\\
    \beta\big(a,\sigma(u,v)\big)&=\frac32\tau(au,v).\label{sigma.beta}
    \end{align}
    \subno{iii} Compatibility among three odd vectors\textup{:}
    \begin{align}
    3\big(\sigma(u,v)w-\sigma(u,w)v\big)&=\eta(u,v)w-\eta(u,w)v-2\eta(v,w)u.\label{sigma.uvw}
\end{align}
\end{lemma}
\begin{proof}
Let $x\in\fg$, $a\in A$ and $u,v,w\in U$. By \eqref{def.sigma}, $[x_0,u_i]=(x_0u)_i$ for $i=0,1$, and $[x_0,L(-1)]=0$, we have
\begin{align*}
    \rho(x)\sigma(u,v)&=x_0\big(u_0v-\frac12L(-1)u_1v\big)\\
    &=[x_0,u_0]v+u_0(x_0v)-\frac{1}{2}L(-1)\big([x_0,u_1]v+u_1(x_0v)\big)\\
    &=(xu)_0v+u_0(xv)-\frac{1}{2}L(-1)(xu)_1v-\frac{1}{2}L(-1)u_1(xv)\\
    &=\sigma(xu,v)+\sigma(u,xv).
\end{align*}
This proves \eqref{sigma.g}. Similarly, using $x_1u=x_1v=0$, \eqref{comm}, \eqref{der} and \eqref{compt.uv.3}, we obtain
\begin{align*}
    \sigma(u,v)\cdot x&=x_1\sigma(u,v)=x_1(u_0v)-\frac12[x_1,L(-1)](u_1v)\\
    &=[x_1,u_0]v-\frac{1}{2}x_0(u_1v)\\
    &=(x_0u)_1v+(x_1u)_0v-\frac{1}{2}[x,u_1v]\\
    &=\eta(xu,v)-\frac{1}{2}\big(\eta(xu,v)+\eta(u,xv)\big)\\
    &=\frac12\big(\eta(xu,v)-\eta(u,xv)\big).
\end{align*}
Thus \eqref{sigma.x} holds. Moreover, using $x_2u=x_2v=0$ together with \eqref{comm}, \eqref{der} and \eqref{compt.uv.5}, we can deduce that
\begin{align*}
    \theta\big(x,\sigma(u,v)\big)\bone&=x_2\sigma(u,v)=[x_2,u_0]v-\frac{1}{2}[x_2,L(-1)](u_1v)\\
    &=(xu)_2v-x_1\eta(u,v)\\
    &=\tau(xu,v)\bone-\kappa\big(x,\eta(u,v)\big)\bone=0,
\end{align*}
which proves \eqref{sigma.theta}.

Since $a_0u=\frac{2}{3}L(-1)(a_1u)$ and $a_1u=\frac32au$, then
\begin{equation*}
    [a_1,u_0]=(a_0u)_1+(a_1u)_0=\big(L(-1)(au)\big)_1+\frac{3}{2}(au)_0=\frac{1}{2}(au)_0.
\end{equation*}
Applying it to $v$ gives $a_1(u_0v)=\frac12(au)_0v+\frac32u_0(av)$. By skew-symmetry,
\begin{align*}
    a_0\eta(u,v)=-\eta(u,v)_0a+L(-1)\big(\eta(u,v)_1a\big)=-\rho\big(\eta(u,v)\big)a+L(-1)\big(a\cdot\eta(u,v)\big).
\end{align*}
Then by \eqref{decom.uv}, we can get
\begin{align*}
    2a\circ\sigma(u,v)&=a_1\sigma(u,v)=a_1\big(u_0v-\frac12L(-1)\eta(u,v)\big)\\
    &=\frac12(au)_0v+\frac32u_0(av)-\frac12a_0\eta(u,v)-\frac{1}{2}L(-1)a_1\eta(u,v)\\
    &=\frac12\sigma(au,v)+\frac32\sigma(u,av)
      +\frac12\rho\big(\eta(u,v)\big)a\\
    &\quad+L(-1)\big(\frac14\eta(au,v)+\frac34\eta(u,av)-a\cdot\eta(u,v)\big)\\
    &=\frac12\sigma(au,v)+\frac32\sigma(u,av)+\frac12\rho\big(\eta(u,v)\big)a,
\end{align*}
where the last equality follows from \eqref{compt.uv.4}. Interchanging $u$ and $v$, using the symmetry of $\sigma$ and the skew-symmetry of $\eta$, we have
\begin{equation*}
    2a\circ\sigma(u,v)=\frac32\sigma(au,v)+\frac12\sigma(u,av)-\frac12\rho\big(\eta(u,v)\big)a.
\end{equation*}
Adding and subtracting these two identities proves \eqref{sigma.A} and \eqref{sigma.rho}, respectively.

For \eqref{sigma.beta}, consider the equality $[a_3,u_0]v=\sum_{i=0}^3\binom{3}{i}(a_iu)_{3-i}v$. On the left-hand side,
\begin{align*}
    [a_3,u_0]v&=a_3(u_0v)-u_0(a_3v)=a_3(u_0v)=a_3\sigma(u,v)+\frac{1}{2}a_3L(-1)\eta(u,v)\\
    &=a_3\sigma(u,v)+\frac{3}{2}a_2\eta(u,v)=\beta\big(a,\sigma(u,v)\big)\bone-\frac{3}{2}\theta\big(\eta(u,v),a\big)\bone=\beta\big(a,\sigma(u,v)\big)\bone,
\end{align*}
where the last equality is from \eqref{compt.uv.6}. On the other hand,
\begin{align*}
    \sum_{i=0}^3\binom{3}{i}(a_iu)_{3-i}v&=(a_0u)_3v+3(a_1u)_2v=\frac{2}{3}\big(L(-1)(a_1u)\big)_3v+3(a_1u)_2v\\
    &=(a_1u)_2v=\frac32\tau(au,v)\bone.
\end{align*}
Equating both sides, we obtain \eqref{sigma.beta}.

By skew-symmetry, we have $v_1b=b_1v=\frac32bv$ for $b\in A$ and $v_0y=-y_0v=-yv$ for $y\in\fg$. Consider the equality $[v_1,u_0]_+w=(v_0u)_1w+(v_1u)_0w$. By \eqref{decom.uv}, its left-hand side is
\begin{align*}
    v_1(u_0w)+u_0(v_1w)&=v_1\sigma(u,w)+\frac{1}{2}[v_1,L(-1)]\eta(u,w)+u_0\eta(v,w)\\
    &=\frac{3}{2}\sigma(u,w)v+\frac{1}{2}v_0\eta(u,w)-\eta(v,w)_0u\\
    &=\frac32\sigma(u,w)v-\frac12\eta(u,w)v-\eta(v,w)u,
\end{align*}
whereas the right-hand side is
\begin{align*}
    (v_0u)_1w+(v_1u)_0w&=\sigma(v,u)_1w+\frac{1}{2}\big(L(-1)\eta(v,u)\big)_1w+\eta(v,u)_0w\\
    &=\frac{3}{2}\sigma(u,v)w-\frac{1}{2}\eta(v,u)_0w+\eta(v,u)_0w\\
    &=\frac32\sigma(u,v)w-\frac12\eta(u,v)w,
\end{align*}
where $\sigma(u,v)=\sigma(v,u)$ by \eqref{symm.sigma} and $\eta(u,v)=-\eta(v,u)$ by \eqref{eta}. Equating these two expressions gives \eqref{sigma.uvw}. This completes the proof.
\end{proof}

\begin{lemma}\label{lem:tau.invar}
    $\tau$ is invariant under the actions of $\fg$ and $A$. More explicitly,
    \begin{equation*}
        \tau(xu,v)=-\tau(u,xv),\ \ \tau(au,v)=\tau(u,av),\ \ \forall\ x\in\fg,\,a\in A,\,u,v\in U.
    \end{equation*}
\end{lemma}
\begin{proof}
    Let $x\in\fg$, $a\in A$ and $u,v\in U$. Since $\tau,\sigma$ are symmetric, then by \eqref{compt.uv.5} and \eqref{sigma.beta}, we have
    \begin{gather*}
        \tau(xu,v)+\tau(u,xv)=\tau(xu,v)+\tau(xv,u)=\kappa\big(x,\eta(u,v)+\eta(v,u)\big)=0,\\
        \tau(au,v)=\frac{2}{3}\beta\big(a,\sigma(u,v)\big)=\frac{2}{3}\beta\big(a,\sigma(v,u)\big)=\tau(av,u)=\tau(u,av). \qedhere
    \end{gather*}
\end{proof}

Combine the structures obtained above into the quadruple
\begin{equation}\label{quadruple}
    \cQ(V)=\Big((\fg,[\cdot,\cdot],\kappa),(A,\circ,\beta),\,(\rho,\cdot,\theta),(U,\tau,\eta,\sigma)\Big),
\end{equation}
where the first three components form the restricted even datum \eqref{datum.even.A}, and the fourth component is $U$ equipped with the $(\fg,A)$-module structure established above. The maps $\tau$, $\eta$ and $\sigma$ take values in $\bC$, $\fg$ and $A$, respectively. The following three tables summarize the bilinear maps,  module structures, and the compatibility conditions in this quadruple. Throughout, $x,y,z\in\fg$, $a,b,c\in A$, and $u,v\in U$.

\begingroup
\renewcommand{\arraystretch}{1.40}
\setlength{\tabcolsep}{4pt}
\setlength{\abovecaptionskip}{6pt}
\setlength{\belowcaptionskip}{6pt}
\setlength{\aboverulesep}{3pt}
\setlength{\belowrulesep}{3pt}

\makeatletter
\setlength{\@fptop}{0pt}
\setlength{\@fpsep}{18pt}
\setlength{\@fpbot}{0pt plus 1fil}
\makeatother

\begin{table}[H]
\centering
\small
\linespread{1}\selectfont

\begin{threeparttable}
\caption{Product and module structures in the quadruple $\cQ(V)$.}
\label{tab:product.actions}

\begin{tabularx}{\textwidth}{
    @{}c
    >{\centering\arraybackslash}X
    >{\centering\arraybackslash}X@{}
}
\toprule
Product or action & Definition & Algebraic structure \\
\midrule

Product on $A$
& $a\circ b=\frac12a_1b$
& Commutative associative algebra with identity $\omega$\tnote{1} \\

$\fg$-action on $A$
& $\rho(x)a=x_0a$
& Lie algebra homomorphism $\rho:\fg\to\Der(A)$ \\

$A$-actions on $\fg$
& $a\cdot x=x\cdot a=a_1x=x_1a$\tnote{2}
& Unital symmetric $A$-bimodule \\

$\fg$-action on $U$
& $xu=x_0u$
& $\fg$-module \\

$A$-action on $U$
& $au=\frac23a_1u$
& Unital $A$-module \\

\bottomrule
\end{tabularx}

\begin{tablenotes}
\footnotesize
\item[1] The identity element $\omega\in A$ satisfies
$\omega\circ a=a$, $\omega\cdot x=x$ and $\omega u=u$.
\item[2] Also write $ax=xa$ for $a\cdot x=x\cdot a$.
\end{tablenotes}
\end{threeparttable}
\end{table}

\begin{table}[H]
\centering
\small
\linespread{1}\selectfont

\begin{threeparttable}
\caption{Bilinear maps in the quadruple $\cQ(V)$.}
\label{tab:bilinear.maps}

\begin{tabularx}{\textwidth}{
    @{}c
    >{\centering\arraybackslash}X
    c
    >{\centering\arraybackslash}X@{}
}
\toprule
Map & Definition & Symmetry & Invariance \\
\midrule

$\kappa:\fg\times\fg\to\bC$
& $x_1y=\kappa(x,y)\bone$
& symmetric
& $\kappa([x,y],z)=\kappa(x,[y,z])$ \\

$\beta:A\times A\to\bC$
& $a_3b=\beta(a,b)\bone$
& symmetric
& $\beta(a\circ b,c)=\beta(a,b\circ c)$ \\

$\theta:\fg\times A\to\bC$
& $x_2a=\theta(x,a)\bone$\tnote{1}
& ---
& --- \\

$\tau:U\times U\to\bC$
& $u_2v=\tau(u,v)\bone$
& symmetric
& $\begin{gathered}
    \tau(xu,v)=-\tau(u,xv)\\
    \tau(au,v)=\tau(u,av)
  \end{gathered}$ \\

$\eta:U\times U\to\fg$
& $\eta(u,v)=u_1v$
& skew-symmetric
& --- \\

$\sigma:U\times U\to A$
& $\sigma(u,v)=u_0v-\frac12L(-1)(u_1v)$
& symmetric
& --- \\

\bottomrule
\end{tabularx}

\begin{tablenotes}[flushleft]
\footnotesize
\item[1] $a_2x=-\theta(x,a)\bone$.
\end{tablenotes}
\end{threeparttable}
\end{table}

\begin{table}[H]
\centering
\small
\linespread{1}\selectfont

\caption{Compatibility conditions in the quadruple $\cQ(V)$.}
\label{tab:compatibility}

\begin{tabularx}{\textwidth}{
    @{}c
    >{\centering\arraybackslash}X
    c@{}
}
\toprule
Structures involved & Content & Equations \\
\midrule

$\fg$ and $A$
& Compatibility of the mutual actions with the Lie bracket
  and the product on $A$
& \eqref{compt.1}, \eqref{compt.2} \\

$\kappa,\beta,\theta$
& Compatibility of the even bilinear maps with the mutual actions
& \eqref{theta.1}--\eqref{theta.4} \\

Actions on $U$
& Compatibility of the $\fg$-action and the $A$-action
& \eqref{compt.xav} \\

$\tau,\eta$
& Compatibility with the actions, and relations to the even bilinear maps
& \eqref{compt.uv.3}--\eqref{compt.uv.6} \\

$\sigma$ and the even datum
& Equivariance and compatibility with the actions and bilinear maps
& \eqref{sigma.g}--\eqref{sigma.beta} \\

$\sigma,\eta$ and $U$
& Compatibility involving three odd elements
& \eqref{sigma.uvw} \\

\bottomrule
\end{tabularx}
\end{table}

\endgroup

\begin{proposition}\label{prop:mode.relations}
    Let $V=\bigoplus_{q\in\frac12\bZ}V_{(q)}$ be a conformal vertex superalgebra strongly generated by $\fg+A+U$, where $\fg=V_{(1)}$, $U=V_{(\frac{3}{2})}$ and $A\subset V_{(2)}$ is defined by \eqref{A}. Suppose that $u_0v\in A+L(-1)\fg$ for all $u,v\in U$. Then
    \allowdisplaybreaks[1]
    \begin{align}
        &[x_m,y_n]=[x,y]_{m+n}+m\kappa(x,y)\delta_{m+n,0},\label{mode.gg}\\
        &[a_m,b_n]=(m-n)(a\circ b)_{m+n-1}+\binom{m}{3}\beta(a,b)\delta_{m+n-2,0},\label{mode.AA}\\
        &[x_m,a_n]=\big(\rho(x)a\big)_{m+n}+m(ax)_{m+n-1}+\binom{m}{2}\theta(x,a)\delta_{m+n-1,0},\label{mode.gA}\\
        &[x_m,u_n]=(xu)_{m+n},\label{mode.gU}\\
        &[a_m,u_n]=\frac12(m-2n)(au)_{m+n-1},\label{mode.AU}\\
        &[u_m,v_n]_+=\sigma(u,v)_{m+n}+\frac12(m-n)\eta(u,v)_{m+n-1}+\binom{m}{2}\tau(u,v)\delta_{m+n-1,0},\label{mode.UU}
    \end{align}
    for $x,y\in\fg$, $u,v\in U$, $a,b\in A$ and $m,n\in\bZ$, where $[u_m,v_n]_+=u_mv_n+v_nu_m$. Here the operations are those of the quadruple \eqref{quadruple}, satisfying \eqref{compt.1}, \eqref{compt.2}, \eqref{theta.1}--\eqref{theta.4}, \eqref{compt.xav}, \eqref{compt.uv.3}--\eqref{compt.uv.6} and \eqref{sigma.g}--\eqref{sigma.uvw} \textup{(}or refer to Tables~\ref{tab:product.actions}--\ref{tab:compatibility}\textup{)}.
\end{proposition}
\begin{proof}
    Using the restricted datum \eqref{datum.even.A} and the operator identities proved in Subsection~\ref{sect:even.generators}, the first three relations are from the restrictions of \eqref{even.mode.gg}--\eqref{even.mode.gA} to $a,b\in A$. It remains to verify the relations involving $U$.

    Let $x,y\in\fg$, $u,v\in U$ and $a,b\in A$. Since $x_iu\in V_{(\frac32-i)}=0$ for $i\ge1$, we have $[x_m,u_n]=(x_0u)_{m+n}=(xu)_{m+n}$. Similarly, due to $a_iu=0$ for $i\ge2$, using \eqref{der} and \eqref{a_1v}, we obtain
    \begin{align*}
        [a_m,u_n]&=(a_0u)_{m+n}+m(a_1u)_{m+n-1}=\frac23\big(L(-1)a_1u\big)_{m+n}+m(a_1u)_{m+n-1}\\
        &=\frac13(m-2n)(a_1u)_{m+n-1}=\frac{1}{2}(m-2n)(au)_{m+n-1}.
    \end{align*}
    Finally, combining \eqref{trun.uv}, \eqref{eta} and \eqref{decom.uv} gives
    \begin{align*}
        [u_m,v_n]_+&=(u_0v)_{m+n}+m(u_1v)_{m+n-1}+\binom{m}{2}(u_2v)_{m+n-2}\\
        &=\sigma(u,v)_{m+n}+\frac12\big(L(-1)\eta(u,v)\big)_{m+n}+m\eta(u,v)_{m+n-1}+\binom{m}{2}\tau(u,v)\delta_{m+n-1,0}\\
        &=\sigma(u,v)_{m+n}+\frac12(m-n)\eta(u,v)_{m+n-1}+\binom{m}{2}\tau(u,v)\delta_{m+n-1,0}.
    \end{align*}
    The compatibility relations have been established above. This completes the proof.
\end{proof}

When $\fg=0$, strong generation and weight counting imply $A=V_{(2)}$. Moreover, $\eta(u,v)=0$ and $\sigma(u,v)=u_0v$ for all $u,v\in U$. Consequently, combining \eqref{mode.AA}, \eqref{mode.AU} and \eqref{mode.UU} with the relations $a\circ b=\frac12a_1b$ and $au=\frac23a_1u$, we recover the commutation relations in Proposition~\ref{prop:LY}. In particular,
\begin{equation*}
    \frac12(m-2n)(au)_{m+n-1}
    =\frac13(m-2n)(a_1u)_{m+n-1}.
\end{equation*}

\subsection{Lie superalgebra $\cL(\fg,A,U)$}\label{sect:Lie.salg}
\hspace{1.5em}We now construct a Lie superalgebra $\cL(\fg,A,U)$ from an abstract quadruple. Let $Z$ be an even vector space, and let
\begin{equation*}
    \cQ_Z=\Big((\fg,[\cdot,\cdot],\kappa),(A,\circ,\beta),
    (\rho,\cdot,\theta),(U,\tau,\eta,\sigma)\Big)
\end{equation*}
be an abstract quadruple equipped with the algebraic structures, module structures and bilinear maps described in Tables~\ref{tab:product.actions} and \ref{tab:bilinear.maps},
except that the bilinear maps
\begin{equation*}
    \kappa:\fg\times\fg\to Z,\ \ \beta:A\times A\to Z,\ \ \theta:\fg\times A\to Z,\ \ \tau:U\times U\to Z,
\end{equation*}
take values in $Z$. All identities involving these maps are understood as equalities in $Z$. More explicitly, $(\fg,[\cdot,\cdot])$ is a Lie algebra, $(A,\circ)$ is a commutative associative algebra, $\rho:\fg\to\Der(A)$ is a Lie algebra homomorphism, $\fg$ is a symmetric $A$-bimodule, and $U$ is both a $\fg$-module and an $A$-module. The bilinear maps $\kappa$, $\beta$, $\tau$ and $\sigma$ are symmetric, and $\eta$ is skew-symmetric. Moreover,
\begin{gather*}
    \kappa([x,y],z)=\kappa(x,[y,z]),\ \ \beta(a\circ b,c)=\beta(a,b\circ c),
\end{gather*}
for all $x,y,z\in\fg$ and $a,b,c\in A$.

\begin{theorem}\label{thm:construction.L}
Let $\cQ_Z$ be a quadruple as above. Set
\begin{equation*}
    \cL(\fg,A,U)=\big((\fg\oplus A\oplus U)\otimes\bC[t,t^{-1}]\big)\oplus Z,
\end{equation*}
View $\cL(\fg,A,U)$ as a superspace with
\begin{align*}
    \cL(\fg,A,U)^{\ep}=\big((\fg\oplus A)\otimes\bC[t,t^{-1}]\big)\oplus Z,\ \ \cL(\fg,A,U)^{\op}=U\otimes\bC[t,t^{-1}].
\end{align*}
Write $\xi(n)=\xi\otimes t^n$ for $\xi\in\fg\oplus A\oplus U$ and $n\in\bZ$. Define a bilinear operation $[\cdot,\cdot]$ on $\cL(\fg,A,U)$ by
\begin{align}
    &[x(m),y(n)]=[x,y](m+n)+m\kappa(x,y)\delta_{m+n,0},\label{L:gg}\\
    &[a(m),b(n)]=(m-n)(a\circ b)(m+n-1)+\binom{m}{3}\beta(a,b)\delta_{m+n-2,0},\label{L:AA}\\
    &[x(m),a(n)]=\big(\rho(x)a\big)(m+n)+m(ax)(m+n-1)\notag\\
    &\qquad\qquad\qquad+\binom{m}{2}\theta(x,a)\delta_{m+n-1,0},\label{L:gA}\\
    &[x(m),u(n)]=(xu)(m+n),\label{L:gU}\\
    &[a(m),u(n)]=\frac12(m-2n)(au)(m+n-1),\label{L:AU}\\
    &[u(m),v(n)]=\sigma(u,v)(m+n)+\frac12(m-n)\eta(u,v)(m+n-1)\notag\\
    &\qquad\qquad\qquad+\binom{m}{2}\tau(u,v)\delta_{m+n-1,0},\label{L:UU}\\
    &[Z,\cL(\fg,A,U)]=0,\label{L:K}
\end{align}
for $x,y\in\fg$, $a,b\in A$, $u,v\in U$ and $m,n\in\bZ$. All remaining brackets are determined by skew-symmetry. Then $\cL(\fg,A,U)$ is a Lie superalgebra if and only if all compatibility conditions in Table~\ref{tab:compatibility}, namely \eqref{compt.1}, \eqref{compt.2}, \eqref{theta.1}--\eqref{theta.4}, \eqref{compt.xav}, \eqref{compt.uv.3}--\eqref{compt.uv.6} and \eqref{sigma.g}--\eqref{sigma.uvw} hold.
\end{theorem}
The proof of this theorem is given in Appendix~\ref{app:Jacobi}.

\section{Vertex operator superalgebra associated with $\cL(\fg,A,U)$}\label{sect:VOSA}
\subsection{Construction of $V_\cL(\ell,0)$}
\hspace{1.5em}Let $Z$ be a vector space. For every $\ell\in Z^*=\mathrm{Hom}_{\bC}(Z,\bC)$, we associate a vertex superalgebra $V_\cL(\ell,0)$ to the Lie superalgebra $\cL(\fg,A,U)$.

Fix a quadruple $\cQ_Z=\Big((\fg,[\cdot,\cdot],\kappa),(A,\circ,\beta),(\rho,\cdot,\theta),(U,\tau,\eta,\sigma)\Big)$ satisfying all assumptions and compatibility conditions in Theorem~\ref{thm:construction.L}. We denote the identity of $(A,\circ)$ by $e$, corresponding to $\omega$ in Table~\ref{tab:product.actions}, while 
$\omega$ continues to denote the conformal vector in a vertex operator superalgebra. Thus
\begin{equation}\label{unital.e}
    e\circ a=a,\ \ ex=x,\ \ eu=u,\ \ \forall\ a\in A,\ x\in\fg,\ u\in U.
\end{equation}

By Theorem~\ref{thm:construction.L}, the quadruple $\cQ_Z$ determines a Lie superalgebra $\cL(\fg,A,U)$, which is a $\frac{1}{2}\bZ$-graded algebra with
\begin{equation*}
    \deg x(n)=-n,\ \ \deg u(n)=-n+\frac{1}{2},\ \ \deg a(n)=-n+1,\ \ \deg c=0,
\end{equation*}
for $x\in\fg$, $u\in U$, $a\in A$, $c\in Z$ and $n\in\bZ$. For $\xi\in\fg+A+U$, form a generating function
\begin{equation*}
    \xi(z)=\sum_{n\in\bZ}\xi(n)z^{-n-1}\in\cL(\fg,A,U)[[z,z^{-1}]].
\end{equation*}

Set 
\begin{equation*}
    \cL^+=(\fg\oplus A\oplus U)\otimes\bC[t]\oplus Z,\ \ \cL^-=(\fg\oplus A\oplus U)\otimes t^{-1}\bC[t^{-1}],
\end{equation*}
which are subalgebras with $\cL(\fg,A,U)=\cL^+\oplus\cL^-$ as a vector space.

Fix $\ell\in Z^*$, and write
\begin{equation*}
    \kappa_\ell=\ell\circ\kappa,\ \ \beta_\ell=\ell\circ\beta,\ \ \theta_\ell=\ell\circ\theta,\ \ \tau_\ell=\ell\circ\tau,
\end{equation*}
where $\circ$ denotes composition of maps. Let $\bC_\ell$ be the one-dimensional $\cL^+$-module on which $(\fg\oplus A\oplus U)\otimes\bC[t]$ acts trivially and each $c\in Z$ acts as the scalar $\ell(c)$. Note that $\cL^+$ is a direct sum of algebras. We can define an induced module
\begin{equation*}
    V_\cL(\ell,0)=U\big(\cL(\fg,A,U)\big)\otimes_{U(\cL^+)}\bC_\ell,
\end{equation*}
where $U(\cdot)$ denotes the universal enveloping algebra of a Lie superalgebra. By the Poincar\'e--Birkhoff--Witt (PBW) theorem, we have $V_\cL(\ell,0)\cong U(\cL^-)$ as a $U(\cL^-)$-module. Furthermore, we have
\begin{equation*}
    V_\cL(\ell,0)\cong S\big((\fg\oplus A)\otimes t^{-1}\bC[t^{-1}]\big)\otimes \Lambda(U\otimes t^{-1}\bC[t^{-1}])
\end{equation*}
as a vector superspace, where $S(\cdot)$ and $\Lambda(\cdot)$ denote the symmetric algebra and the exterior algebra, respectively.

Set 
\begin{equation*}
    \bone=1\otimes 1\in V_\cL(\ell,0).
\end{equation*}
Since $(\fg\oplus A\oplus U)\otimes\mathbb{C}[t]$ acts trivially on $\mathbb{C}_\ell$, we have $\xi(n)\bone=0$ for $\xi\in \fg+A+U$ and $n\ge0$. Meanwhile, identify $\fg\oplus A\oplus U$ as a subspace of $V_\cL(\ell,0)$ through the linear map
\begin{equation}\label{creation}
    x+a+u\mapsto x(-1)\bone+a(-1)\bone+u(-1)\bone,
\end{equation}
for $x\in\fg$, $a\in A$ and $u\in U$. For $\xi$ in one of the summands $\fg,A,U$, write $h_\xi=1,2,\frac{3}{2}$, respectively. Then $\deg \xi(n)=h_\xi-n-1$. Due to the fact that $\cL^+$ is a $\frac{1}{2}\bZ$-graded subalgebra, $V_\cL(\ell,0)$ is a $\frac{1}{2}\bZ$-graded $\cL(\fg,A,U)$-module by defining $\deg \bone=0$, that is,
\begin{equation*}
    V_\cL(\ell,0)=\bigoplus_{q\in\frac{1}{2}\bZ} V_\cL(\ell,0)_{(q)}.
\end{equation*}
By the PBW theorem, $V_\cL(\ell,0)_{(q)}=0$ for $q<0$. Since $\xi(n)V_\cL(\ell,0)_{(q)}\subset V_\cL(\ell,0)_{(q+h_\xi-n-1)}$, we have $\xi(n)V_\cL(\ell,0)_{(q)}=0$ for $n>q+h_\xi-1$.

Define a linear operator $D$ on $\cL(\fg,A,U)$ by
\begin{equation*}
    D(c)=0,\ \ D\big(\xi(n)\big)=-n\xi(n-1),
\end{equation*}
for $c\in Z$, $\xi\in\fg+A+U$ and $n\in\bZ$. Then we have $D\big(\xi(z)\big)=\frac{d}{dz}\xi(z)$. It is easy to see that $D$ is a derivation of $\cL(\fg,A,U)$. Furthermore, $D$ naturally acts as a derivation on $U\big(\cL(\fg,A,U)\big)$. Since $D$ preserves the subalgebra $\cL^+$ and acts trivially on $\bC_\ell$, $D$ induces an operator on $V_\cL(\ell,0)$ by $D(X\otimes 1)=D(X)\otimes 1$ for $X\in U\big(\cL(\fg,A,U)\big)$ such that
\begin{equation*}
    D\bone=0,\ \ [D,\xi(z)]=\frac{d}{dz}\xi(z),\ \ \text{for}\ \xi\in\fg+A+U.
\end{equation*}

$V_\cL(\ell,0)$ is $\bZz$-graded:
\begin{equation*}
    V_\cL(\ell,0)^\ep=\bigoplus_{n\in\bZ}V_\cL(\ell,0)_{(n)},\ \  V_\cL(\ell,0)^\op=\bigoplus_{n\in\frac{1}{2}+\bZ}V_\cL(\ell,0)_{(n)}.
\end{equation*}

We now show that all generating functions are mutually local to each other, by operator product expansions (OPEs). One can refer to \cite{K98} for the details on OPEs. 

\begin{lemma}
    For $V_\cL(\ell,0)$, the generating functions satisfy the following OPEs\textup{:}
    \allowdisplaybreaks[4]
    \begin{align*}
        x(z)y(w) &\sim \frac{\kappa_\ell(x,y)}{(z-w)^2}+\frac{[x,y](w)}{z-w}, \\
        a(z)b(w) &\sim \frac{\beta_\ell(a,b)}{(z-w)^4}+\frac{2(a \circ b)(w)}{(z-w)^2}+\frac{\partial_w(a \circ b)(w)}{z-w}, \\
        x(z)a(w) &\sim \frac{\theta_\ell(x,a)}{(z-w)^3}+\frac{(ax)(w)}{(z-w)^2}+\frac{\big(\rho(x)a\big)(w)}{z-w},\\
        x(z)u(w) &\sim \frac{(xu)(w)}{z-w},\\
        a(z)u(w) &\sim \frac{\frac{3}{2}(au)(w)}{(z-w)^2}+\frac{\partial_w(au)(w)}{z-w}, \\
        u(z)v(w) &\sim \frac{\tau_\ell(u,v)}{(z-w)^3}+\frac{\eta(u,v)(w)}{(z-w)^2}+\frac{\sigma(u,v)(w)+\frac{1}{2}\partial_w\eta(u,v)(w)}{z-w}.
    \end{align*}
\end{lemma}

Then by Li's local system theory for vertex superalgebras developed in \cite{Li}, the following result is standard.

\begin{proposition}\label{prop:VSA}
     For $\ell\in Z^*$, $V_\cL(\ell,0)$ admits a vertex superalgebra structure, which is uniquely determined by the condition that $\bone$ is the vacuum vector and
     \begin{equation*}
         Y(\xi,z)=\xi(z)=\sum_{n\in\bZ}\xi(n)z^{-n-1},\ \ \text{for}\ \xi\in\fg+A+U.
     \end{equation*}
     Note that in the standard vertex operator notation $Y(\xi,z)=\sum_{n\in\bZ}\xi_nz^{-n-1}$, we have
     \begin{equation*}
         \xi_n=\xi(n),\ \ \text{for}\ \xi\in\fg+A+U,\,n\in\bZ.
     \end{equation*}
\end{proposition}

Furthermore, we have the following theorem:
\begin{theorem}\label{thm:VOSA}
    The vector $\omega=e(-1)\bone$ is a conformal vector of $V_\cL(\ell,0)$, with central charge $c_V=2\beta_\ell(e,e)$. If $\fg$, $A$ and $U$ are finite-dimensional, then $V_\cL(\ell,0)$ is a vertex operator superalgebra of CFT-type.
\end{theorem}
\begin{proof}
    Set $L(m)=e(m+1)$, so that $Y(\omega,z)=\sum_{m\in\bZ}L(m)z^{-m-2}$. Then by \eqref{L:AA} and $e\circ e=e$, we obtain
    \begin{equation*}
        [L(m),L(n)]=(m-n)L(m+n)+\frac{\beta_\ell(e,e)}{6}(m^3-m)\delta_{m+n,0}.
    \end{equation*}
    This is the Virasoro relation with central charge $2\beta_\ell(e,e)$. Since $\rho(x)$ is a derivation of the unital algebra $(A,\circ)$, then
    \begin{equation}\label{rho(x)e}
        \rho(x)e=\rho(x)(e\circ e)=\big(\rho(x)e\big)\circ e+e\circ\big(\rho(x)e\big)=2\rho(x)e,
    \end{equation}
    namely $\rho(x)e=0$. Using \eqref{unital.e} and the brackets in Theorem~\ref{thm:construction.L}, we get
\begin{align}
    [L(m),x(n)]&=-n\,x(m+n)-\frac{\theta_\ell(x,e)}{2}(m^2+m)\delta_{m+n,0},\label{L.x}\\
    [L(m),a(n)]&=(m+1-n)a(m+n)+\frac{\beta_\ell(e,a)}{6}(m^3-m)\delta_{m+n-1,0},\label{L.a}\\
    [L(m),u(n)]&=\frac{1}{2}(m+1-2n)u(m+n).\label{L.u}
\end{align} 
We can obtain $L(-1)=D$ by taking $m=-1$ in these formulae. Taking $m=0$ in \eqref{L.x}--\eqref{L.u} and using $L(0)\bone=0$, we can deduce that
\begin{align*}
    L(0)x=L(0)x(-1)\bone=[L(0),x(-1)]\bone=x(-1)\bone=x,
\end{align*}
$L(0)a=2a$ and $L(0)u=\frac{3}{2}u$ for $x\in\fg$, $a\in A$ and $u\in U$. Then the conformal weight grading condition \eqref{L(0)} follows. Furthermore,
\begin{equation}\label{V.CFT}
    V_{\cL}(\ell,0)_{(0)}=\bC\bone,\ \ V_{\cL}(\ell,0)_{(\frac{1}{2})}=0,\ \ V_{\cL}(\ell,0)_{(1)}=\fg,\ \ V_{\cL}(\ell,0)_{(\frac{3}{2})}=U,
\end{equation}
and $V_{\cL}(\ell,0)_{(q)}=0$ for $q\in-\frac{1}{2}\bZ_+$. 

Finally, if the generating spaces are finite-dimensional, only finitely many ordered PBW monomials have any fixed degree. Hence all conformal weight spaces are finite-dimensional. This together with \eqref{V.CFT} shows that $V_\cL(\ell,0)$ is a vertex operator superalgebra of CFT-type.
\end{proof}

\subsection{Simple quotients and reconstruction}
\hspace{1.5em}In this subsection, we assume that $\fg$, $A$ and $U$ are finite-dimensional and $\fg$ is perfect, that is, $[\fg,\fg]=\fg$. Thus $V_\cL(\ell,0)$ is a vertex operator superalgebra of CFT-type.

First, we recall the following definitions and facts (cf. \cite{DH,FHL,Li94,Sch,Xu}).
\begin{definition}
    Let $V$ be a vertex operator superalgebra. A bilinear form $\langle\cdot,\cdot\rangle$ on $V$ is said to be \emph{invariant} if 
    \begin{equation}\label{invariant}
        \langle Y(u,z)v,w\rangle=(-1)^{|u||v|}\langle v,Y(e^{zL(1)}z^{-2L(0)}e^{\pi iL(0)}u,z^{-1})w\rangle,
    \end{equation}
    for $u,v,w\in V$, where $i=\sqrt{-1}$ and by definition $e^{\pi i L(0)}v=e^{q\pi i}v$ for $v\in V_{(q)}$ with $q\in\frac{1}{2}\bZ$.
\end{definition}

\begin{proposition}\label{prop:symm.invar.bilin.form}
    Let $V=\bigoplus_{q\in\frac{1}{2}\bZ}V_{(q)}$ be a vertex operator superalgebra. Any invariant bilinear form $\langle\cdot,\cdot\rangle$ on $V$ is automatically symmetric. Moreover,
    \begin{equation}\label{involution}
        \langle L(n)u,v\rangle=\langle u,L(-n)v\rangle\ \ \text{for}\ n\in\bZ,\,u,v\in V.
    \end{equation}
    In particular, $\langle L(0)u,v\rangle=\langle u,L(0)v\rangle$. It follows that $\langle V_{(p)},V_{(q)}\rangle=0$ for $p,q\in\frac{1}{2}\bZ$ with $p\neq q$.
\end{proposition}

The existence and uniqueness of the normalized invariant bilinear form follow from \cite{Li94,Sch,Xu}. 
\begin{proposition}\label{prop:bilinear.form}
    Let $V=\bigoplus_{q\in\frac12\bZ}V_{(q)}$ be a vertex operator superalgebra of CFT-type such that $L(1)V_{(1)}=0$. Then $V$ admits a symmetric invariant bilinear form $\langle\cdot,\cdot\rangle$ uniquely determined by $\langle\bone,\bone\rangle=1$. 
\end{proposition}

Furthermore, we have the following corollary of Proposition~\ref{prop:bilinear.form}, or one can refer to \cite{LY} for a proof.
\begin{corollary}\label{cor:simple}
    The kernel of the normalized invariant bilinear form $\langle\cdot,\cdot\rangle$ is the unique maximal proper ideal of $V$. $V$ is simple if and only if $\langle\cdot,\cdot\rangle$ is non-degenerate.
\end{corollary}

We now return to a vertex operator superalgebra $V=\bigoplus_{q\in\frac{1}{2}\bZ}V_{(q)}$ strongly generated by $\fg+A+U$, where $\fg=V_{(1)}$ is perfect, $U=V_{(\frac32)}$ and $A$ is the subspace defined by \eqref{A}. Assume that $u_0v\in A+L(-1)\fg$ for $u,v\in U$, and use the associated quadruple $\cQ(V)$ from Section~\ref{sect:data}. 
\begin{lemma}\label{lem:L(1)g=0}
    For any $x\in\fg$ and $a\in A$, we have $\theta(x,a)=0$. In particular, $L(1)\fg=0$.
\end{lemma}
\begin{proof}
    By \eqref{theta.2} and \eqref{rho(x)e}, we obtain
    \begin{equation*}
        \theta([x,y],\omega)=\theta\big(x,\rho(y)\omega\big)=0,
    \end{equation*}
    for every $x,y\in\fg$, where $\omega$ is the identity of $A$. Since $\fg$ is perfect, we have $\theta(x,\omega)=0$ for any $x\in\fg$. Furthermore, by \eqref{theta.1}, we can deduce that
    \begin{equation*}
        \theta(x,a)=\theta(x,a\circ \omega)=\theta(ax,\omega)=0,\ \ \forall\ x\in\fg,\,a\in A.
    \end{equation*}
    By \eqref{a_2x}, we have $L(1)x=\omega_2x=-\theta(x,\omega)\bone=0$. Hence $L(1)\fg=0$.
\end{proof}
Then by Lemma~\ref{lem:CFT}, Proposition~\ref{prop:bilinear.form} and Lemma \ref{lem:L(1)g=0}, $V$ has a unique normalized symmetric invariant bilinear form. Moreover, $L(1)A=0$ by \eqref{A} and $L(1)U\subset V_{(\frac12)}=0$. For $x\in\fg$, $a\in A$ and $u\in U$, by \eqref{invariant}, we have
\begin{align*}
    Y(e^{zL(1)}z^{-2L(0)}e^{\pi iL(0)}x,z^{-1})&=-Y(e^{zL(1)}z^{-2L(0)}x,z^{-1})=-z^{-2}Y(e^{zL(1)}x,z^{-1})\\
    &=-z^{-2}Y(x,z^{-1})=-\sum_{m\in\bZ}x_mz^{m-1}=-\sum_{n\in\bZ}x_{-n}z^{-n-1},\\
    Y(e^{zL(1)}z^{-2L(0)}e^{\pi iL(0)}a,z^{-1})&=\sum_{n\in\bZ}a_{2-n}z^{-n-1},\\
    Y(e^{zL(1)}z^{-2L(0)}e^{\pi iL(0)}u,z^{-1})&=(-i)\sum_{n\in\bZ}u_{1-n}z^{-n-1}.
\end{align*}
Then for $x,y\in\fg$, $a,b\in A$ and $u,v\in U$,
\begin{align*}
    \langle x,y\rangle&=\langle x_{-1}\bone,y\rangle=-\langle \bone,x_1y\rangle=-\kappa(x,y)\langle\bone,\bone\rangle=-\kappa(x,y),\\
    \langle a,b\rangle&=\beta(a,b),\ \ \langle u,v\rangle=-i\tau(u,v).
\end{align*}
Furthermore, we have the following proposition.
\begin{proposition}\label{prop:simple->}
    Let $V=\bigoplus_{q\in\frac12\bZ}V_{(q)}$ be a simple vertex operator superalgebra strongly generated by $\fg+A+U$, where $\fg=V_{(1)}$ is perfect, $U=V_{(\frac{3}{2})}$ and $A\subset V_{(2)}$ is defined by \eqref{A}. Suppose that $u_0v\in A+L(-1)\fg$ for all $u,v\in U$. Then the associated quadruple
    \begin{equation*}
        \cQ(V)=\Big((\fg,[\cdot,\cdot],\kappa),(A,\circ,\beta),(\rho,\cdot,\theta),(U,\tau,\eta,\sigma)\Big)
    \end{equation*}
    has the structures and compatibility conditions listed in Tables~\ref{tab:product.actions}--\ref{tab:compatibility}. Moreover, $\theta=0$, the forms $\kappa$ and $\tau$ are non-degenerate, and
    \begin{equation}\label{simple.con}
        \Rad\, \beta\cap \Ann_A(\fg)=0,
    \end{equation}
    where $\Rad\beta=\{a\in A\mid \beta(a,b)=0,\ \forall\ b\in A\}$ and $\Ann_A(\fg)
    =\{a\in A\mid ax=0,\ \forall\ x\in\fg\}$.
\end{proposition}
\begin{proof}
     Lemma~\ref{lem:L(1)g=0} shows $\theta=0$. Since $V$ is simple, by Proposition~\ref{prop:bilinear.form} and Corollary~\ref{cor:simple}, $V$ admits a non-degenerate symmetric invariant bilinear form uniquely determined by $\langle\bone,\bone\rangle=1$. Then by $V_{(1)}=\fg$ and $V_{(\frac{3}{2})}=U$, we know that $\kappa(x,y)=-\langle x,y\rangle$ and $\tau(u,v)=i\langle u,v\rangle$ are non-degenerate, where $x,y\in\fg$, $u,v\in U$.

    Suppose $a\in \Rad\,\beta\cap \Ann_A(\fg)$. Let $b\in A$ and $x,y\in\fg$. Then
    \begin{equation}\label{beta.1}
        \langle a,b\rangle=\beta(a,b)=0.
    \end{equation} 
    By \eqref{involution} and $L(1)a=0$, we have
    \begin{equation}\label{beta.2}
        \langle a,L(-1)x\rangle=\langle L(1)a,x\rangle=0.
    \end{equation}
    By symmetry and invariance of $\langle\cdot,\cdot\rangle$, we obtain
    \begin{equation}\label{beta.3}
        \langle a,x_{-1}y\rangle=\langle x_{-1}y,a\rangle=-\langle y,x_1a\rangle=\kappa(y,ax)=0.
    \end{equation}
    On the other hand, strong generation and weight counting show that
    \begin{equation}\label{V_(2)}
        V_{(2)}=A+L(-1)\fg+\Span\{x_{-1}y\mid x,y\in\fg\}.
    \end{equation}
    Combining \eqref{beta.1}--\eqref{V_(2)} implies that $a$ is orthogonal to $V_{(2)}$. Since $\langle\cdot,\cdot\rangle$ is non-degenerate, we obtain $a=0$. This completes the proof.
\end{proof}

\begin{remark}\label{rmk:deg.form}
    Unlike the setting in \cite[Corollary~3.22]{LY}, $A$ is in general a proper subspace of $V_{(2)}$. Thus the non-degeneracy of the invariant bilinear form on $V$ does not directly imply the non-degeneracy of its restriction $\beta$ to $A$. Indeed, $\langle a,x_{-1}y\rangle=\kappa(y,ax)$, for $a\in A$ and $x,y\in\fg$. So an element of $\Rad\,\beta$ need not be orthogonal to $V_{(2)}$, which motivates the additional condition \eqref{simple.con} and allows $\beta$ to be degenerate. When $\fg=0$, strong generation gives $V_{(2)}=A$, and this condition reduces to the non-degeneracy of $\beta$, as in \cite{LY}.
\end{remark}

We now turn to the converse construction, starting from an abstract quadruple $\cQ_Z$ and the Lie superalgebra $\cL(\fg,A,U)$ associated to it, where $\fg$, $A$ and $U$ are finite-dimensional and $\fg$ is perfect. The argument proving $\theta=0$ in Lemma~\ref{lem:L(1)g=0} applies equally to the $Z$-valued map $\theta$. Hence $\theta=0$ and $L(1)\fg=0$ in $V_\cL(\ell,0)$.

We next give sufficient conditions for the simple quotient of $V_\cL(\ell,0)$ to retain the generating subspaces $\fg$, $A$ and $U$.

By \eqref{V.CFT} and Lemma~\ref{lem:L(1)g=0}, we have $V_\cL(\ell,0)_{(0)}=\bC\bone$ and $L(1)V_\cL(\ell,0)_{(1)}=L(1)\fg=0$. Then by Proposition~\ref{prop:bilinear.form}, there exists a unique symmetric invariant bilinear form $\langle\cdot,\cdot\rangle$ on $V_\cL(\ell,0)$ with $\langle\bone,\bone\rangle=1$. For $x,y\in\fg$, by \eqref{L:gg}, we have
\begin{align*}
    x_1y=x_1y_{-1}\bone=[x_1,y_{-1}]\bone+y_{-1}x_1\bone=[x,y]_0\bone+\kappa(x,y)\bone=\ell\big(\kappa(x,y)\big)\bone=\kappa_\ell(x,y)\bone.
\end{align*}
Similarly, for $a,b\in A$ and $u,v\in U$, by \eqref{L:AA} and \eqref{L:UU}, we can get
\begin{equation*}
    a_3b=\beta_\ell(a,b)\bone,\ \ u_2v=\tau_\ell(u,v)\bone.
\end{equation*}
By invariance \eqref{invariant}, we obtain
\begin{equation*}
    \langle x,y\rangle=-\kappa_\ell(x,y),\ \ \langle a,b\rangle=\beta_\ell(a,b),\ \ \langle u,v\rangle=-i\tau_\ell(u,v).
\end{equation*}
The radical of $\langle\cdot,\cdot\rangle$
\begin{equation*}
    J=\{v\in V_{\cL}(\ell,0)\mid \langle v,V_{\cL}(\ell,0)\rangle=0\}
\end{equation*}
is the unique maximal proper ideal of $V_{\cL}(\ell,0)$. Set
\begin{equation*}
    L_{\cL}(\ell,0)=V_{\cL}(\ell,0)/J.
\end{equation*}
We have the following result.
\begin{proposition}\label{prop:simple<-}
    Assume that $\kappa_\ell$ and $\tau_\ell$ are non-degenerate and $\Rad\,\beta_\ell\cap\Ann_A(\fg)=0$. Then the quotient map $V_\cL(\ell,0)\to L_\cL(\ell,0)$ is injective on $\fg\oplus A\oplus U$. Consequently, the simple vertex operator superalgebra $L_\cL(\ell,0)$ satisfies 
    \begin{equation*}
        L_\cL(\ell,0)_{(1)}=\mathfrak g,\ \ L_\cL(\ell,0)_{(\frac32)}=U,\ \ A\subset L_\cL(\ell,0)_{(2)},
    \end{equation*}
    and $L_\cL(\ell,0)$ is strongly generated by $\fg+A+U$.
\end{proposition}
\begin{proof}
    Since $\kappa_\ell,\tau_\ell$ are non-degenerate, the restrictions of $\langle\cdot,\cdot\rangle$ to $\fg$ and $U$ are non-degenerate. Hence $J\cap \fg=0$ and $J\cap U=0$. Let $a\in J\cap A$. For every $b\in A$, $0=\langle a,b\rangle=\beta_\ell(a,b)$, so $a\in\Rad\,\beta_\ell$. Moreover, since $J$ is an ideal, $ax=x_1a\in J\cap\mathfrak g=0$, for any $x\in\fg$. Thus $a\in \Rad\,\beta_\ell\cap \Ann_A(\fg)=0$, which implies that $J\cap A=0$.  Since $J$ is graded and the three generating subspaces have distinct weights, the quotient map is injective on $\fg\oplus A\oplus U$. The remaining assertions follow.
\end{proof}

Finally, we can obtain the following isomorphism.
\begin{theorem}\label{thm:reconstruction}
    Take $Z=\bC$ and $\ell_0=\Id\in Z^*$. Let $V=\bigoplus_{q\in\frac12\bZ}V_{(q)}$ be a simple vertex operator superalgebra strongly generated by $\fg+A+U$, where $\fg=V_{(1)}$ is perfect, $U=V_{(\frac{3}{2})}$ and $A\subset V_{(2)}$ is defined by \eqref{A}. Suppose that $u_0v\in A+L(-1)\fg$ for all $u,v\in U$. Let $\cL(\fg,A,U)$ be the Lie superalgebra constructed from the associated quadruple $\cQ(V)$. Then there exists an isomorphism of vertex operator superalgebras
    \begin{equation*}
        L_{\cL}(\ell_0,0)\cong V,
    \end{equation*}
    where $L_\cL(\ell_0,0)$ is the unique simple quotient of $V_\cL(\ell_0,0)$.
\end{theorem}
\begin{proof}
    The universal property of the induced module gives a homomorphism
    \begin{equation*}
        \Phi:V_\cL(\ell_0,0)\longrightarrow V,
    \end{equation*}
    such that $\Phi(\xi)=\xi$ for $\xi\in\fg+A+U$, and thus
    \begin{equation*}
        \Phi\big(Y(\xi,z)w\big)=Y_V\big(\Phi(\xi),z\big)\Phi(w),\ \ \forall\ w\in V_{\cL}(\ell_0,0).
    \end{equation*}
    Since the identity $e$ corresponds to the conformal vector $\omega$ of $V$, we have $\Phi\big(e(-1)\bone\big)=\omega$. Then $\Phi$ is a homomorphism of vertex operator superalgebras. Strong generation implies that it is surjective. Since $V$ is simple, $\ker\Phi$ is the unique maximal proper ideal $J$ of $V_\cL(\ell_0,0)$. Hence $V\cong V_\cL(\ell_0,0)/J=L_\cL(\ell_0,0)$.
\end{proof}

\section{Tensor product decomposition and examples}\label{sect:examples}
\hspace{1.5em}In this section, we first establish a tensor product decomposition of the universal vertex operator superalgebras constructed in Section~\ref{sect:VOSA}, and then present several examples.

Notably, for the quadruple $\cQ_Z=\Big((\fg,[\cdot,\cdot],\kappa),(A,\circ,\beta),(\rho,\cdot,\theta),(U,\tau,\eta,\sigma)\Big)$ and the Lie superalgebra $\cL(\fg,A,U)$ associated to it constructed in Section~\ref{sect:data}, if $U=0$, the additional condition $a_0U\subset L(-1)U$ in \eqref{A} is trivial, so that $A=\cA$. In this case, $\cL(\fg,A,0)$ is an ordinary Lie algebra, denoted by
\begin{equation*}
    \cL(\fg,\cA)=\big((\fg\oplus\cA)\otimes\bC[t,t^{-1}]\big)\oplus Z.
\end{equation*}
By the construction in Section~\ref{sect:VOSA}, each of the Lie (super)algebras described in this section gives rise to a family of conformal vertex (super)algebras.

For convenience, we also write $ab$ for $a\circ b$ and use the corresponding notation for powers and ideals, where $a,b\in A$. 

\subsection{Tensor product decomposition}\label{sect:tensor.decomposition}
\hspace{1.5em}In this subsection, we assume that $\fg$, $A$ and $U$ are finite-dimensional. We first recall some definitions on the commutative algebra. One can refer to \cite{AM} for more details. 

For $a\in A$, $a$ is said to be \emph{idempotent} if $a^2=a\neq0$. The \emph{Jacobson radical} $\fJ(A)$ of $A$ is defined to be the intersection of all its maximal ideals. The algebra $A$ is called \emph{local} if it has a unique maximal ideal. A commutative associative algebra $A$ is called \emph{Artinian} if, regarded as an $A$-module, it admits a finite composition series
\begin{equation*}
    A=I_0\supset I_1\supset\cdots\supset I_n=0,
\end{equation*}
where each $I_j$ is an ideal of $A$ and each quotient $I_j/I_{j+1}$ is a simple $A$-module.

Since $A$ is finite-dimensional, it is Artinian. Then by the \emph{structure theorem for Artin rings}, we can obtain the following result.
\begin{proposition}\label{prop:Artin}
    Let $A$ be a finite-dimensional unital commutative associative algebra over $\bC$ with identity $e$. Then there exist a positive integer $r$ and idempotents $e_1,\dots,e_r$ such that
    \begin{gather*}
        e=e_1+\cdots+e_r,\ \ e_i  e_j=\delta_{i,j}e_i,\ \ \forall\ i,j=1,\dots,r;\\ A=A_1\oplus\cdots\oplus A_r,\ \ A_i=e_i  A.
    \end{gather*}
    Each $A_i$ is a local Artinian algebra with identity $e_i$. Moreover,
    \begin{equation}\label{Artin}
        A_i/\fJ(A_i)\cong \bC,\ \ \fJ(A_i)^{N_i}=0,
    \end{equation}
    for some positive integer $N_i$.
\end{proposition}
\begin{proof}
    From \cite[Theorem~8.7]{AM}, there exists an isomorphism
    \begin{equation*}
        \varphi:A\longrightarrow B_1\times\cdots\times B_r,
    \end{equation*}
    where each $B_i$ is a local Artinian algebra. Denote the identity of $B_i$ by $1_i$ and define $e_i=\varphi^{-1}(0,\ldots,1_i,\ldots,0)$. Then
    \begin{equation*}
        e_i^2=e_i,\ \ e_i  e_j=0\quad(i\ne j).
    \end{equation*}
    It follows that $\sum_{i=1}^{r}e_i=\varphi^{-1}(1_1,\ldots,1_r)=e$. Setting $A_i=e_i  A$, due to the fact that 
    \begin{equation*}
        \varphi(A_i)=0\times\cdots\times B_i\times\cdots\times0,
    \end{equation*}
    we have $A_i\cong B_i$ and $A=A_1\oplus\cdots\oplus A_r$. In particular, each $A_i$ is local and Artinian. Moreover, for $b=e_i  a\in A_i$, $e_i  b=e_i^2  a=e_i  a=b$. Thus $e_i$ is the identity of $A_i$.

    Since $A_i$ is local, $\fJ(A_i)$ is its unique maximal ideal. Note that $A_i$ is finite-dimensional over $\bC$. Then by \cite[Corollary~7.10]{AM}, we have $A_i/\fJ(A_i)\cong \bC$. 

    Finally, \cite[Corollary~8.2 and Proposition~8.4]{AM} imply that $\fJ(A_i)^{N_i}=0$ for some positive integer $N_i$. This completes the proof.
\end{proof}

For the rest of this subsection, we use the notation in Proposition~\ref{prop:Artin}.
\begin{lemma}\label{lem:M_i}
    Let $M$ be a module over $A$ and denote $M_i=e_iM$ for $i=1,\cdots,r$. Then we have $M=\bigoplus_{i=1}^rM_i$.
\end{lemma}
\begin{proof}
    For any $m\in M$, we can obtain $m=em=\sum_{i=1}^re_im$. So $M=M_1+\cdots+M_r$. Now we show that it is a direct sum. Assume that $m_1+\cdots+m_r=0$ for $m_i\in M_i$. Then by $e_ie_j=\delta_{i,j}e_i$, we have 
    \begin{equation*}
        0=e_j\big(\sum_{i=1}^rm_i\big)=m_j,\ \ \forall\ j=1,\dots,r.
    \end{equation*}
    Hence $M=\bigoplus_{i=1}^re_iM=\bigoplus_{i=1}^rM_i$.
\end{proof}

Write $\fg_i=e_i\fg$ and $U_i=e_iU$ for $i=1,\dots,r$. Then by Lemma~\ref{lem:M_i} and the fact that $\fg$ and $U$ are modules over $A$, we have
\begin{equation}\label{g_iU_i}
    \fg=\bigoplus_{i=1}^r\fg_i,\ \ U=\bigoplus_{i=1}^r U_i.
\end{equation}

\begin{proposition}\label{prop:decom.BMS}
    Let $\pi:\cL(\fg,A,U)\rightarrow\cL(\fg,A,U)/Z$ be the canonical projection, and define
    \begin{equation*}
        \overline{\cL}_i=\pi\big((\fg_i\oplus A_i\oplus U_i)\otimes\bC[t,t^{-1}]\big).
    \end{equation*}
    Then we have
    \begin{equation*}
        \cL(\fg,A,U)/Z\cong\bigoplus_{i=1}^r\overline{\cL}_i
    \end{equation*}
    as a direct sum of Lie superalgebras.
\end{proposition}
\begin{proof}
    By Proposition~\ref{prop:Artin} and \eqref{g_iU_i}, the decompositions of $A$, $\fg$ and $U$ lead to
    \begin{equation*}
        \cL(\fg,A,U)/Z=\bigoplus_{i=1}^{r}\overline{\cL}_i\ \ \text{as vector spaces.}
    \end{equation*}
    
    In the following, it suffices to verify that each $\overline{\cL}_i$ is a Lie superalgebra and $[\overline{\cL}_i,\overline{\cL}_j]=0$ whenever $i\neq j$. For any $D\in\Der(A)$, since $e_i^2=e_i$, a similar equation as in \eqref{rho(x)e} shows that $D(e_i)=0$. Therefore,
    \begin{equation}\label{rho(x)e_i}
        \rho(x)e_i=0,\ \ \forall\ x\in\fg.
    \end{equation}
    Then by \eqref{compt.2}, we have $\rho(e_ix)a-e_i\rho(x)a=\rho(ax)e_i-a\rho(x)e_i=0$, which implies that
    \begin{equation}\label{rho.linear}
        \rho(e_ix)a=e_i\rho(x)a,\ \ \forall\ x\in\fg,\, a\in A.
    \end{equation}
    
    Let $x,y\in\fg$, $a,b\in A$ and $u,v\in U$. It is easy to see that
    \begin{equation}\label{ab_i}
        (e_i  a) (e_j  b)=\delta_{i,j}e_i  (a  b),\ \ (e_i  a)(e_jx)=\delta_{i,j}e_iax.
    \end{equation} 
    By \eqref{rho(x)e_i} and \eqref{rho.linear},
    \begin{align}\label{rho_i}
        \rho(e_ix)(e_ja)=e_i\rho(x)(e_ja)=e_i\big(\rho(x)e_j\big)a+e_ie_j\big(\rho(x)a\big)=\delta_{i,j}e_i\rho(x)a.
    \end{align}
    Moreover, by \eqref{compt.1} and \eqref{rho(x)e_i}, we can obtain
    \begin{align}
        [e_ix,e_jy]&=e_j[e_ix,y]+\big(\rho(e_ix)e_j\big)y=e_j[e_ix,y]=-e_j[y,e_ix]\notag\\
        &=-e_je_i[y,x]-e_j\big(\rho(y)e_i\big)x=\delta_{i,j}e_i[x,y]\label{[xy]_i}.
    \end{align}
    By \eqref{compt.xav} and \eqref{rho(x)e_i}, we have $(e_ix)u=x(e_iu)=e_i(xu)$. It follows that
    \begin{align}\label{[xu]_i}
        (e_ix)(e_ju)=e_ix(e_ju)=e_ie_jxu=\delta_{i,j}e_ixu.
    \end{align}
    Since $U$ is an $A$-module and $A$ is commutative and associative, 
    \begin{equation}\label{[au]_i}
        (e_ia)(e_ju)=(e_iae_j)u=e_ie_j(au)=\delta_{i,j}e_iau.
    \end{equation}
    By \eqref{sigma.rho} and \eqref{rho(x)e_i}, we have $\sigma(e_iu,v)=\sigma(u,e_iv)$. This together with \eqref{sigma.A} shows that 
    \begin{equation}\label{sigma_i}
        e_i \sigma(u,v)=\sigma(e_iu,v)=\sigma(u,e_iv).
    \end{equation}
    It follows from \eqref{compt.uv.4} and \eqref{sigma_i} that
    \begin{equation}\label{[uu]_i}
        \eta(e_iu,e_jv)=\delta_{i,j}e_i\eta(u,v),\ \ \sigma(e_iu,e_jv)=\delta_{i,j}e_i\sigma(u,v).
    \end{equation}
    Since all central terms vanish, by the relations in Theorem~\ref{thm:construction.L} and the equalities \eqref{ab_i}--\eqref{[au]_i} and \eqref{[uu]_i}, we can conclude that
    \begin{equation*}
        [\overline{\cL}_i,\overline{\cL}_i]\subset \overline{\cL}_i,\ \ [\overline{\cL}_i,\overline{\cL}_j]=0\ \ (i\neq j).
    \end{equation*}
    Consequently, this proposition holds.
\end{proof}

Furthermore, the central terms between distinct components also vanish. More precisely, for $i \neq j$, since $e_i   e_j = 0$, the invariance of $\beta$ and $\tau$ together with \eqref{theta.1} yields
\begin{equation}\label{central_i}
\beta(e_i  a, e_j  b)=0,\ \ \tau(e_iu, e_jv)=0, \ \ \theta(e_ix,e_j  a)=0,
\end{equation}
for all $x\in\fg$, $a,b\in A$ and $u,v\in U$. Moreover, for any $x,y\in\fg$ and $i\neq j$, taking $a=e_j$ in \eqref{theta.4} and using \eqref{[xy]_i}, we obtain
\begin{equation}\label{central_i.2}
\kappa(e_ix,e_jy)=\kappa\big(e_ix,e_j(e_jy)\big)-\kappa\big(e_j(e_ix),e_jy\big)=\theta([e_ix,e_jy],e_j)=0.
\end{equation}

For each $i= 1,\dots,r$, let $\cQ_{Z,i}$ be the restriction of the quadruple $\cQ_Z$ to $(\fg_i,A_i,U_i)$, where $e_i$ is the identity of $A_i$. Let
\begin{equation*}
    \cL_i:=\cL(\fg_i,A_i,U_i)=\widetilde{\cL}_i\oplus Z
\end{equation*}
be the corresponding Lie superalgebra constructed in Theorem~\ref{thm:construction.L}, where $\widetilde{\cL}_i=(\fg_i\oplus A_i\oplus U_i)\otimes\mathbb{C}[t,t^{-1}]$. By Proposition~\ref{prop:decom.BMS}, \eqref{central_i} and \eqref{central_i.2}, the subspaces $\widetilde{\cL}_i$ and $\widetilde{\cL}_j$ satisfy
\begin{equation}\label{L_i.commute}
[\widetilde{\cL}_i,\widetilde{\cL}_j]=0 \ \ (i\neq j).
\end{equation}

To characterize the corresponding decomposition at the level of vertex operator superalgebras, we recall the tensor product construction (cf. \cite{FHL,K98,Li}).

\begin{definition}
    Let $(V_i,Y_i,\bone_i,\omega_i)$, $i=1,2$, be vertex operator superalgebras. Their tensor product is the superspace $V_1\otimes V_2$, equipped with the usual tensor product weight and parity gradings, vacuum vector and conformal vector
\begin{equation*}
    \bone=\bone_1\otimes\bone_2,\ \ \omega=\omega_1\otimes\bone_2+\mathbf1_1\otimes\omega_2,
\end{equation*}
and vertex operator map
\begin{equation*}
    Y(a_1\otimes a_2,z)=Y_1(a_1,z)\otimes Y_2(a_2,z).
\end{equation*}
\end{definition}

Here and throughout, tensor products of operators are understood in the graded sense: for homogeneous operators $S\in\End(V_1)$, $T\in\End(V_2)$ and homogeneous vectors $v_i\in V_i$,
\begin{equation*}
    (S\otimes T)(v_1\otimes v_2)=(-1)^{|T||v_1|}Sv_1\otimes Tv_2.
\end{equation*}

\begin{theorem}\label{thm:tensor_decomp}
For any $\ell\in Z^*$, let $V_{\cL}(\ell,0)$ and $V_{\cL_i}(\ell,0)$ $(1\le i\le r)$ be the vertex operator superalgebras associated with $\cL(\fg,A,U)$ and $\cL_i=\cL(\fg_i,A_i,U_i)$, respectively. Then there is an isomorphism of vertex operator superalgebras
\begin{equation*}
    V_{\cL}(\ell,0)\cong V_{\cL_1}(\ell,0)\otimes\cdots\otimes V_{\cL_r}(\ell,0).
\end{equation*}
\end{theorem}
\begin{proof}
Let $W=V_{\cL_1}(\ell,0)\otimes\cdots\otimes V_{\cL_r}(\ell,0)$ be the tensor product vertex operator superalgebra, with vacuum vector $\bone_W=\bone_1\otimes\cdots\otimes\bone_r$ and conformal vector
\begin{equation*}
    \omega_W=\sum_{i=1}^r\bone_1\otimes\cdots\otimes\omega^{(i)}\otimes\cdots\otimes\bone_r,
\end{equation*}
where $\bone_i$ and $\omega^{(i)}=e_i(-1)\bone_i$ denote the vacuum vector and conformal vector of $V_{\cL_i}(\ell,0)$, respectively. For each $i$ and $\xi_i\in\fg_i+A_i+U_i$, let $\iota_i(\xi_i)=\bone_1\otimes\cdots \otimes\xi_i(-1)\bone_i\otimes\cdots\otimes\bone_r\in W$. Then
\begin{equation*}
    Y_W\big(\iota_i(\xi_i),z\big)=\Id\otimes\cdots\otimes Y_i(\xi_i,z)\otimes\cdots\otimes\Id=\sum_{n\in\bZ}\big(\Id\otimes\cdots\otimes\xi_i(n)\otimes\cdots\otimes\Id\big)z^{-n-1}.
\end{equation*}
For distinct tensor factors, we have
\begin{equation*}
    [Y_W\big(\iota_i(\xi_i),z_1\big),Y_W\big(\iota_j(\xi_j),z_2\big)]=0\ \ (i\ne j).
\end{equation*}
Combining this with \eqref{L_i.commute}, we can deduce that $W$ becomes an $\cL(\fg,A,U)$-module by letting
\begin{equation*}
    \xi_i(n)\longmapsto\Id\otimes\cdots\otimes\xi_i(n)\otimes\cdots\otimes\Id,\ \ c\longmapsto\ell(c)\Id_W\ \ (c \in Z).
\end{equation*}
Since $(\fg\oplus A\oplus U)\otimes\bC[t]$ annihilates $\bone_W$ and $c\in Z$ acts on $\bone_W$ as $\ell(c)$, the universal property of the induced module $V_{\cL}(\ell,0)=U\big(\cL(\fg,A,U)\big)\otimes_{U(\cL^+)}\bC_\ell$ yields a unique $\cL(\fg,A,U)$-module homomorphism
\begin{equation*}
    \Psi:V_{\cL}(\ell,0)\longrightarrow W;\ \ \bone\longmapsto\bone_W.
\end{equation*}
By \eqref{L_i.commute}, the subalgebra $\cL^-=(\fg\oplus A\oplus U)\otimes t^{-1}\bC[t^{-1}]$ decomposes into a direct sum of Lie superalgebras $\cL^-=\bigoplus_{i=1}^r\cL_i^-$, where $\cL_i^-=(\fg_i\oplus A_i\oplus U_i)\otimes t^{-1}\bC[t^{-1}]$. Hence $U(\cL^-)\cong U(\cL_1^-)\otimes\cdots\otimes U(\cL_r^-)$, and the PBW theorem implies that $\Psi$ is a linear isomorphism.

Notice that $\Psi\big(\xi(-1)\bone\big)=\sum_{i=1}^r\iota_i(\xi_i)$ for $\xi=\sum_{i=1}^r\xi_i\in\fg+A+U$ with $\xi_i\in\fg_i+A_i+U_i$. In particular, since $e=\sum_{i=1}^r e_i$, we have $\Psi(\omega)=\Psi\big(e(-1)\bone\big)=\omega_W$. Moreover,
\begin{equation*}
    \Psi\big(Y(\xi, z)v\big)=\Psi\big(\xi(z)v\big)=\sum_{i=1}^r Y_W\big(\iota_i(\xi_i),z\big)\Psi(v)=Y_W\big(\Psi(\xi),z\big)\Psi(v)
\end{equation*}
for all $\xi\in\fg+A+U$ and $v\in V_{\cL}(\ell,0)$. By the uniqueness assertion in Proposition~\ref{prop:VSA}, $\Psi$ is therefore an isomorphism of vertex superalgebras. Since $\Psi(\omega)=\omega_W$, it is an isomorphism of vertex operator superalgebras. This completes the proof.
\end{proof}

Consequently, in the finite-dimensional setting, the study of the universal vertex operator superalgebras $V_{\cL}(\ell,0)$ reduces to the case where $A$ is a local Artinian algebra.

\subsection{Affine-Virasoro algebra}
\hspace{1.5em}The affine-Virasoro algebra is the semi-direct sum of the Virasoro algebra $\mathrm{Vir}$ and the untwisted affine Lie algebra. It is also named the conformal current algebra (\cite{K86}) and the entire gauge algebra (\cite{EFK}). 

Let $\fg$ be a simple Lie algebra of finite dimension, and let
\begin{gather*}
    \cA=\bC e,\ \ e\circ e=e,\ \ ex=x,\ \ \rho=\theta=0,\\
    Z=\bC K\oplus\bC C,\ \ \kappa(x,y)=(x,y)K,\ \ \beta(e,e)=\frac{1}{2}C,
\end{gather*}
where $(\cdot,\cdot)$ is the Killing form on $\fg$. Set $L_m=e(m+1)$ and $x_m=x(m)$. Then $\mathrm{AVir}:=\cL(\fg,\cA)$ is the \emph{affine-Virasoro algebra}, which is a Lie algebra generated by $\{ L_m,x_m,K,C\mid m\in\bZ\}$ satisfying the following non-vanishing relations:
\begin{align*}
    &[x_m,y_n]=[x,y]_{m+n}+m(x,y)\delta_{m+n,0}K,\\
    &[L_m,L_n]=(m-n)L_{m+n}+\frac{m^3-m}{12}\delta_{m+n,0}C,\\
    &[L_m,x_n]=-nx_{m+n},\ \ \forall\ x,y\in\fg,\,m,n\in\bZ.
\end{align*}

For $\ell\in Z^*$, write $k=\ell(K)$ and $c=\ell(C)$. By Theorem~\ref{thm:VOSA},
the universal vacuum module $V_{\cL}(\ell,0)$ is a vertex operator algebra of CFT-type with conformal vector $\omega=e(-1)\bone$ and central charge $c$. Since $\fg$ is perfect, $V_{\cL}(\ell,0)$ admits a unique simple quotient, denoted by $L_{\cL}(\ell,0)$. We now take $k\neq0$ and $c=0$. Then $\kappa_\ell(x,y)=k(x,y)$ and $\beta_\ell=0$. So $\kappa_\ell$ is non-degenerate. Moreover, since $A=\bC e$ and $e$ acts as the identity on $\fg$, we have $\Ann_A(\fg)=0$, and hence
\begin{equation*}
    \Rad\,\beta_\ell\cap\Ann_A(\fg)=0.
\end{equation*}
As $U=0$, Proposition~\ref{prop:simple<-} implies that the quotient map
$V_{\cL}(\ell,0)\to L_{\cL}(\ell,0)$ is injective on $\fg\oplus A$. Thus, this example shows that nontrivial simple quotients can be obtained even when $\beta_\ell$ is degenerate.

\subsection{Twisted Heisenberg--Virasoro algebra}
\hspace{1.5em}In this example, we take $\fg=\bC h$ as an abelian Lie algebra and
\begin{gather*}
    \cA=\bC e,\ \ e\circ e=e,\ \ eh=h,\ \ \rho=0,\ \ Z=\bC c_1\oplus\bC c_2\oplus \bC c_3,\\ 
    \beta(e,e)=\frac{1}{2}c_1,\ \ \theta(h,e)=2c_2,\ \ \kappa(h,h)=c_3.
\end{gather*}
Set $L_m=e(m+1)$ and $I_n=h(n)$. Then $\mathcal{HV}:=\cL(\fg,\cA)$ is the \emph{twisted Heisenberg--Virasoro algebra} with a basis $\{L_m,I_m,c_1,c_2,c_3\mid m\in\bZ\}$ satisfying the following non-vanishing relations:
\begin{align*}
    &[L_m,L_n]=(m-n)L_{m+n}+\frac{m^3-m}{12}\delta_{m+n,0}c_1,\\
    &[L_m,I_n]=-nI_{m+n}-(m^2+m)\delta_{m+n,0}c_2,\\
    &[I_m,I_n]=m\delta_{m+n,0}c_3,
\end{align*}
for any $m,n\in\bZ$, which was first introduced in \cite{ACKP}. The twisted Heisenberg--Virasoro vertex operator algebra was studied in \cite{GW}. Note that $\fg=\bC h$ is non-perfect. One can refer to \cite{GW} for more details on the simple quotient of the associated vertex operator algebra.

\subsection{Superconformal algebras}
\hspace{1.5em}The notion of a superconformal algebra was introduced by Kac (cf. \cite{K98}). A basic example is the \emph{Neveu--Schwarz algebra} $\mathrm{NS}$, also known as the $N=1$ \emph{superconformal algebra}. The construction in \cite{LY} realizes this algebra as a special case of $\cL(0,A,U)$ in our notation. In this subsection, we first show that the $N=2$ \emph{superconformal algebra} (of Neveu--Schwarz type) is a special case of $\cL(\fg,A,U)$.

Suppose that $A=\bC e$ and $Z=\bC C$ such that $\beta(e,e)=\frac{1}{2}C$. Then
\begin{align*}
    [e(m+1),e(n+1)]&=(m-n)(e\circ e)(m+n+1)+\binom{m+1}{3}\delta_{m+n,0}\beta(e,e),\\
    &=(m-n)e(m+n+1)+\frac{m^3-m}{12}\delta_{m+n,0}C,
\end{align*}
for $m,n\in\bZ$. Hence, by setting $L_m=e(m+1)$ for $m\in\bZ$, we get a copy of the Virasoro algebra
\begin{equation}\label{N2.L}
    [L_m,L_n]=(m-n)L_{m+n}+\frac{m^3-m}{12}\delta_{m+n,0}C,
\end{equation}
for $m,n\in\bZ$. Moreover, suppose that $\fg=\bC h$ with $[h,h]=0$ and $U=\bC u^+\oplus\bC u^-$ with $\sigma(u^+,u^-)=2e\in A$. Then we have $eh=h$ and $eu^\pm=u^\pm$ since $\fg$ and $U$ are unital $A$-modules. Equip $U$ with the $\fg$-module structure determined by $hu^\pm=\pm u^\pm$. By \eqref{rho(x)e}, we can get $\rho=0$. Then it follows from \eqref{sigma.g} that
\begin{equation*}
    0=\sigma(hu^\pm,u^\pm)+\sigma(u^\pm,hu^\pm)=\pm2\sigma(u^\pm,u^\pm),
\end{equation*}
namely $\sigma(u^\pm,u^\pm)=0$. The compatibility conditions \eqref{sigma.x}, \eqref{sigma.theta} and \eqref{sigma.beta} lead to
\begin{gather*}
    4h=2\sigma(u^+,u^-)h=\eta(hu^+,u^-)-\eta(u^+,hu^-)=2\eta(u^+,u^-),\\
    0=\theta\big(h,\sigma(u^+,u^-)\big)=2\theta(h,e),\\
    \frac{3}{2}\tau(u^+,u^-)=\beta\big(e,\sigma(u^+,u^-)\big)=\beta(e,2e)=C,\\
    \frac{3}{2}\tau(u^\pm,u^\pm)=\beta\big(e,\sigma(u^\pm,u^\pm)\big)=0,
\end{gather*}
which respectively means that
\begin{equation*}
    \eta(u^+,u^-)=2h,\ \ \theta=0,\ \ \tau(u^+,u^-)=\frac{2}{3}C\ \ \text{and}\ \ \tau(u^\pm,u^\pm)=0.
\end{equation*}
Since $\eta$ is skew-symmetric, then $\eta(u^\pm,u^\pm)=0.$ Finally, by \eqref{compt.uv.5}, we have
\begin{equation*}
    \kappa(h,h)=\frac{1}{2}\kappa\big(h,\eta(u^+,u^-)\big)=\frac{1}{2}\tau(hu^+,u^-)=\frac{1}{2}\tau(u^+,u^-)=\frac{1}{3}C.
\end{equation*}
Setting $J_m=h(m)$ for $m\in\bZ$, we obtain a copy of the Heisenberg algebra:
\begin{equation}\label{N2.J}
    [J_m,J_n]=\frac{m}{3}\delta_{m+n,0}C.
\end{equation}
We now write $G^\pm_r=u^\pm(r+\frac{1}{2})$ for $r\in\bZ+\frac{1}{2}$. It follows from \eqref{L:UU} that
\begin{align}
    [G^+_r,G^-_s]&=2e(r+s+1)+(r-s)h(r+s)+\frac{1}{3}(r+\frac{1}{2})(r-\frac{1}{2})\delta_{r+s,0}C\notag\\
    &=2L_{r+s}+(r-s)J_{r+s}+\frac{1}{3}(r^2-\frac{1}{4})\delta_{r+s,0}C,\label{N2.UU}
\end{align}
for $r,s\in\bZ+\frac{1}{2}$. One can easily check the remaining relations as follows:
\begin{align}\label{N2.remain}
    [L_m,J_n]=-nJ_{m+n},\  [L_m,G^\pm_r]=(\frac{m}{2}-r)G^\pm_{m+r},\  [J_m,G^\pm_r]=\pm G^\pm_{m+r},\  [G^\pm_r,G^\pm_s]=0,
\end{align}
for $m,n\in\bZ$, $r,s\in\bZ+\frac{1}{2}$. Combining \eqref{N2.L}--\eqref{N2.remain} and the fact that $C$ is central, we get all relations of the $N=2$ superconformal algebra. In other words, under the assumption 
\begin{gather*}
    A=\bC e,\ \ \fg=\bC h,\ \ U=\bC u^+\oplus\bC u^-,\ \ Z=\bC C,\\  
    [h,h]=0,\ \ \beta(e,e)=\frac{1}{2}C,\ \ \sigma(u^+,u^-)=2e,\ \ hu^\pm=\pm u^\pm,
\end{gather*}
the Lie superalgebra $\cL(\fg,A,U)$ is isomorphic to the $N=2$ superconformal algebra.

Another example is the \emph{small} $N=4$ \emph{superconformal algebra}. Kac and Wakimoto \cite{KW} constructed the corresponding vertex superalgebra via quantum Drinfeld--Sokolov reduction. We now realize this superconformal algebra as a special case of $\cL(\fg,A,U)$.

Take
\begin{gather}
    \fg=\mathfrak{sl}_2=\Span\{E,F,H\},\ \ [H,E]=2E,\ \ [H,F]=-2F,\ \ [E,F]=H,\label{sl_2}\\
    A=\bC e,\ \ Z=\bC C,\ \ \beta(e,e)=\frac{1}{2}C. \notag
\end{gather}
Let $V=\mathbb C^2=\mathbb Cu^+\oplus\mathbb Cu^-$ be the two-dimensional standard $\mathfrak{sl}_2$-module, where
\begin{equation*}
    u^+=\begin{pmatrix}1\\0\end{pmatrix},\ \ u^-=\begin{pmatrix}0\\1\end{pmatrix}.
\end{equation*}
Thus, $Hu^\pm=\pm u^\pm$, $Eu^-=u^+$, $Fu^+=u^-$ and $Eu^+=Fu^-=0$.
Set \(U=V\oplus V^*\), where \(V^*=\operatorname{Hom}_{\mathbb C}(V,\mathbb C)\) carries the dual representation given by
\begin{equation*}
    (x\varphi)(u)=-\varphi(xu),\ \ \forall\ x\in\mathfrak g,\ u\in V,\ \varphi\in V^*.
\end{equation*}
Let $\{(u^+)^*,(u^-)^*\}$ be the basis dual to $\{u^+,u^-\}$, and write $v^-=(u^+)^*$, $v^+=(u^-)^*$. Actually, the action on $U$ is given by
\[
\begin{array}{c|cccc}
  & u^+ & u^- & v^+ & v^- \\ \hline
H & u^+ & -u^- & v^+ & -v^- \\
E & 0 & u^+ & 0 & -v^+ \\
F & u^- & 0 & -v^- & 0
\end{array}
\]
We further assume that
\begin{equation}\label{N4.assum.sigma}
    \sigma(u^+,v^-)=2e\in A.
\end{equation}
Now we can construct the small $N=4$ superconformal algebra as follows.

By \eqref{rho(x)e}, we have $\rho=0$. Plugging this into \eqref{sigma.g} shows that $\sigma$ is $\fg$-invariant:
\begin{equation}\label{sigma.invar}
    0=\sigma(xu,v)+\sigma(u,xv),\ \ \forall\ x\in\fg,\,u,v\in U.
\end{equation}
By the symmetry of $\sigma$ and \eqref{sigma.invar}, we can obtain
\begin{equation}\label{N4.sigma}
    \begin{split}
        &0=\sigma(Eu^-,v^-)+\sigma(u^-,Ev^-)=\sigma(u^+,v^-)-\sigma(u^-,v^+),\\
    &0=\sigma(Hu^\pm,u^\pm)+\sigma(u^\pm,Hu^\pm)=\pm2\sigma(u^\pm,u^\pm),\\
    &0=\sigma(Hv^\pm,v^\pm)+\sigma(v^\pm,Hv^\pm)=\pm2\sigma(v^\pm,v^\pm),\\
    &0=\sigma(Hu^\pm,v^\pm)+\sigma(u^\pm,Hv^\pm)=\pm2\sigma(u^\pm,v^\pm),\\
    &0=\sigma(Eu^-,u^-)+\sigma(u^-,Eu^-)=2\sigma(u^+,u^-),\\
    &0=\sigma(Ev^-,v^-)+\sigma(v^-,Ev^-)=-2\sigma(v^+,v^-).
    \end{split}
\end{equation}
Combining \eqref{N4.assum.sigma} and \eqref{N4.sigma}, we can deduce that
\begin{equation}\label{N4.sigma.1}
    \sigma(u^+,v^-)=\sigma(u^-,v^+)=2e,\ \ \sigma(V,V)=\sigma(V^*,V^*)=\sigma(u^\pm,v^\pm)=0.
\end{equation}
By \eqref{sigma.theta}, $0=\theta\big(x,\sigma(u^+,v^-)\big)=2\theta(x,e)$ for all $x\in\fg$. It follows that $\theta=0$. Plugging \eqref{N4.sigma.1} into \eqref{sigma.beta}, we have
\begin{equation*}
    \tau(u^+,v^-)=\tau(u^-,v^+)=\frac{2}{3}C,\ \ \tau(V,V)=\tau(V^*,V^*)=\tau(u^\pm,v^\pm)=0.
\end{equation*}
By \eqref{sigma.x} and \eqref{N4.sigma.1}, we obtain
\begin{align*}
    4H=\eta(Hu^+,v^-)-\eta(u^+,Hv^-)=\eta(u^+,v^-)-\eta(u^+,-v^-)=2\eta(u^+,v^-),
\end{align*}
which implies that $\eta(u^+,v^-)=2H$. Similarly, using the actions of $H,E,F$ on $U$, we can get
\begin{gather*}
    \eta(u^-,v^+)=-2H,\ \ \eta(u^+,v^+)=4E,\ \ \eta(u^-,v^-)=4F.
\end{gather*}
By the skew-symmetry of $\eta$, we also have $\eta(V,V)=\eta(V^*,V^*)=0$. So far, we have determined $\theta$, $\sigma$, $\tau$ and $\eta$. This together with the compatibility condition \eqref{compt.uv.5} and the symmetry of $\kappa$ shows
\begin{equation*}
    \kappa(H,H)=\frac{1}{3}C,\ \ \kappa(E,F)=\frac{1}{6}C,
\end{equation*}
and the others vanish. Finally, set
\begin{gather*}
    L_m=e(m+1),\ \ J_m^0=H(m),\ \ J_m^+=E(m),\ \ J_m^-=F(m),\\
    G_r^\pm=u^\pm(r+\frac{1}{2}),\ \ \overline{G}_r^\pm=v^\pm(r+\frac{1}{2}),
\end{gather*}
for $m\in\bZ$, $r\in\bZ+\frac{1}{2}$. By Theorem~\ref{thm:construction.L}, the brackets between the even generators are as follows:
\begin{equation}\label{N4.rel.1}
    \begin{split}
        &[L_m,L_n]=(m-n)L_{m+n}+\frac{m^3-m}{12}\delta_{m+n,0}C,\\
    &[L_m,J^a_n]=-nJ^a_{m+n},\ \ a\in\{0,+,-\},\\
    &[J^0_m,J^0_n]=\frac{m}{3}\delta_{m+n,0}C,\ \ [J^0_m,J^\pm_n]=\pm2J^\pm_{m+n},\\
    &[J^+_m,J^-_n]=J^0_{m+n}+\frac{m}{6}\delta_{m+n,0}C,\ \ [J^\pm_m,J^\pm_n]=0,
    \end{split}
\end{equation}
which contributes to an affine-Virasoro algebra. The brackets between the even and odd generators are
\begin{equation}\label{N4.rel.2}
    \begin{split}
        &[L_m,G^\pm_r]=(\frac m2-r)G^\pm_{m+r},\ \ [L_m,\overline G^\pm_r]=(\frac m2-r)\overline G^\pm_{m+r},\\
    &[J^0_m,G^\pm_r]=\pm G^\pm_{m+r},\ \ [J^0_m,\overline G^\pm_r]=\pm\overline G^\pm_{m+r},\\
    &[J^\pm_m,G^\mp_r]=G^\pm_{m+r},\ \ [J^\pm_m,\overline G^\mp_r]=-\overline G^\pm_{m+r},\\
    &[J^\pm_m,G^\pm_r]=[J^\pm_m,\overline G^\pm_r]=0.
    \end{split}
\end{equation}
The brackets between the odd generators are
\begin{equation}\label{N4.rel.3}
    \begin{split}
        &[G^+_r,\overline G^-_s]=2L_{r+s}+(r-s)J^0_{r+s}+\frac13\left(r^2-\frac14\right)\delta_{r+s,0}C,\\
    &[G^-_r,\overline G^+_s]=2L_{r+s}-(r-s)J^0_{r+s}+\frac13\left(r^2-\frac14\right)\delta_{r+s,0}C,\\
    &[G^+_r,\overline G^+_s]=2(r-s)J^+_{r+s},\ \ [G^-_r,\overline G^-_s]=2(r-s)J^-_{r+s},\\
    &[G^\epsilon_r,G^{\epsilon'}_s]=[\overline G^\epsilon_r,\overline G^{\epsilon'}_s]=0,
      \ \ \epsilon,\epsilon'\in\{+,-\}.
    \end{split}
\end{equation}
The relations in \eqref{N4.rel.1}--\eqref{N4.rel.3} hold for any $m,n\in\mathbb Z$, $r,s\in\bZ+\frac12$, and $C$ is central.

Thus, we obtain the small $N=4$ superconformal algebra. By Theorem~\ref{thm:VOSA}, the associated vacuum module carries a vertex operator superalgebra structure.

\subsection{Generalized $W(2,2)/\mathrm{BMS}_3$ (super)algebras}
\hspace{1.5em}In this subsection, we assume that $\fg$, $A$ and $U$ are finite-dimensional. In view of Theorem~\ref{thm:tensor_decomp}, we focus on examples arising from local Artinian algebras. 

If $A$ is local, then by \eqref{Artin}, we have $A=\bC e\oplus \fJ(A)$ as vector spaces, where $e$ is the identity of $A$. If $\fJ(A)\neq0$, then we can take nonzero $a\in \fJ(A)$. From \eqref{Artin}, there exists $N\in\bZ_+$ such that $a^N=0$ and $a^{N-1}\neq0$. Therefore, we get a subalgebra $B$ of $A$ generated by $e$ and $a$:
\begin{gather*}
    B=\bC e\oplus \bC a\oplus\cdots \oplus\bC a^{N-1},\ \ \text{satisfying}\ \ a^i  a^j=\begin{cases}
        a^{i+j},\ \ &i+j<N;\\
        0,\ \ & i+j\ge N.
    \end{cases}
\end{gather*}
It is easy to check that $B\cong\bC[\varepsilon]/(\varepsilon^N)$, where $\bC[\varepsilon]/(\varepsilon^N)$ is a \emph{truncated polynomial algebra} with variable $\varepsilon$ if $N\ge 2$.

In what follows, we let $N\ge 2$. We first construct the truncated Virasoro algebra and truncated Neveu--Schwarz algebra. Define
\begin{equation}\label{trun.1}
    B=\bigoplus_{i=0}^{N-1}\bC a^i,\ \ Z=\bigoplus_{i=0}^{N-1}\bC C_i,
\end{equation}
such that $a^0=e$ is the identity and $a^k=C_k=0$ for all $k\ge N$.

Take $\fg=U=0$ and let 
\begin{equation}\label{trun.2}
    \beta(a^i,a^j)=\frac{1}{2}C_{i+j},\ \ 0\le i,j<N.
\end{equation}
Then $\beta$ is symmetric and invariant. Set $L_m^{(i)}=a^i(m+1)$ for $m\in\bZ$, $0\le i<N$. Then we get a Lie algebra $\cL(0,B,0)$ generated by $\{L_m^{(i)},C_i\mid m\in\bZ,\,0\le i<N\}$, satisfying 
\begin{equation}\label{W(2,2)}
    [L_m^{(i)},L_n^{(j)}]=(m-n)L_{m+n}^{(i+j)}+\frac{m^3-m}{12}\delta_{m+n,0}C_{i+j},
\end{equation}
for $m,n\in\bZ$, where $C_i$ is central and $L_m^{(k)}=C_k=0$ for $k\ge N$. Note that $\cL(0,B,0)$ is isomorphic to the \emph{truncated Virasoro algebra} (cf. \cite{JZ,Wilson}) $\mathrm{Vir}\otimes\bC[\varepsilon]/(\varepsilon^N)$ defined by 
\begin{equation*}
    [x\otimes \varepsilon^i,y\otimes \varepsilon^j]=[x,y]\otimes \varepsilon^{i+j},\ \ \forall\ x,y\in\mathrm{Vir},\ i,j=0,\dots,N-1.
\end{equation*}
In particular, if $N=2$, then $\mathrm{Vir}\otimes\bC[\varepsilon]/(\varepsilon^2)$ is the $W(2,2)$ algebra, or called the $\mathrm{BMS}_3$ algebra (cf. \cite{O,ZD}).

Retain $\fg=0$ and the definitions of $B,Z,\beta$ in \eqref{trun.1} and \eqref{trun.2}. Let $U=Bu$ be a free $B$-module of rank one with basis $u$. Then 
\begin{equation}
    U=\bigoplus_{i=0}^{N-1}\bC u_i,\ \ u_i=a^iu.
\end{equation}
So $a^iu_j=u_{i+j}$ for $i+j<N$. We also define 
\begin{equation}\label{trun.sigma}
    \sigma(u_i,u_j)=\begin{cases}
        2a^{i+j},\ \ &i+j<N;\\
        0,\ \ &i+j\ge N.
    \end{cases}
\end{equation}
Since $\fg=0$, we have $\rho=\theta=\eta=0$. By \eqref{sigma.beta} and \eqref{trun.sigma}, $\tau(u_i,u_j)=\frac{2}{3}C_{i+j}$ for $0\le i,j<N$. Set $G_r^{(i)}=u_i(r+\frac{1}{2})$ for $r\in\bZ+\frac{1}{2}$, $0\le i<N$. It follows that $\cL(0,B,Bu)$ is a Lie superalgebra generated by $\{L_m^{(i)}, G_r^{(i)}, C_i\mid m\in\bZ,\,r\in\bZ+\frac{1}{2},\,0\le i<N\}$, satisfying \eqref{W(2,2)} and
\begin{align*}
    &[L_m^{(i)},G_r^{(j)}]=\left(\frac m2-r\right)G_{m+r}^{(i+j)},\\
    &[G_r^{(i)},G_s^{(j)}]=2L_{r+s}^{(i+j)}+\frac13\left(r^2-\frac14\right)\delta_{r+s,0}C_{i+j},
\end{align*}
for $m,n\in\bZ$ and $r,s\in\bZ+\frac{1}{2}$, where $C_i$ is central and $G_r^{(k)}=L_m^{(k)}=C_k=0$ for $k\ge N$. It is actually a \emph{truncated Neveu--Schwarz algebra} (cf. \cite{WLP}) $\mathrm{NS}\otimes\bC[\varepsilon]/(\varepsilon^N)$ defined by  
\begin{equation*}
    [x\otimes \varepsilon^i,y\otimes \varepsilon^j]=[x,y]\otimes \varepsilon^{i+j},\ \ \forall\ x,y\in\mathrm{NS},\ i,j=0,\dots,N-1.
\end{equation*}
In particular, if $N=2$, then $\mathrm{NS}\otimes\bC[\varepsilon]/(\varepsilon^2)$ is the $N=2$ super-$\mathrm{BMS}_3$ algebra (cf. \cite{GPB,Man}). 

Next, we present some examples with $\fg\neq0$. We still take 
\begin{equation*}
    B=\bigoplus_{i=0}^{N-1}\bC a^i\cong\bC[\varepsilon]/(\varepsilon^N),
\end{equation*}
where $a^0=e$ is the identity of $B$. Then $\Der(B)$ is a Lie algebra. From Leibniz's rule, we have 
\begin{equation*}
    D(e)=0,\ \ D(a^j)=ja^{j-1}  D(a),\ \ \forall\ D\in\Der(B),\,1\le j<N.
\end{equation*}
One can check that
\begin{equation*}
    \Der(B)\longrightarrow a  B;\ \ D\longmapsto D(a),\ \ D\in\Der(B),
\end{equation*}
is an isomorphism of vector spaces. Thus $\dim\Der(B)=N-1$, and we denote 
\begin{equation*}
    \Der(B)=\bigoplus_{i=0}^{N-2}\bC D_i,\ \ D_i=a^{i+1}\frac{d}{da}.
\end{equation*}

Let $U=Z=0$ and $\fg=\Der(B)$. Further define 
\begin{equation*}
    \rho(D)=D,\ \ (bD)c=b  D(c),\ \ \forall\ D\in\fg,\, b,c\in B.
\end{equation*}
Then by Leibniz's rule, for any $D,E\in\fg$ and $b,c\in B$, we have
\begin{align*}
    [D,bE](c)&=D\big(bE(c)\big)-bE\big(D(c)\big)=D(b)E(c)+bD\big(E(c)\big)-bE\big(D(c)\big)\\
    &=D(b)E(c)+b[D,E](c)=\big(D(b)E+b[D,E]\big)(c),
\end{align*}
It follows that $[D,bE]=b[D,E]+D(b)E$. Moreover, $\rho(bD)c-b \rho(D)c=(bD)c-b D(c)=0$. Hence, $(\fg,B)$ is a Lie--Rinehart algebra and satisfies
\begin{equation*}
    [D_i,D_j]=(j-i)D_{i+j},\ \ a^i  a^j=a^{i+j},\ \ D_i(a^j)=ja^{i+j},\ \ a^jD_i=D_{i+j},
\end{equation*}
where $D_k=a^l=0$ for $k\ge N-1$, $l\ge N$. Set $L_m^{(i)}=a^i(m+1)$ with $0\le i<N$ and $J_m^{(j)}=D_j(m)$ with $0\le j<N-1$. Then $\mathcal{B}_N:=\cL(\Der(B),B)$ is a Lie algebra generated by $\{ L_m^{(i)},J_m^{(j)}\mid m\in\bZ,\, 0\le i<N,\, 0\le j<N-1\}$ satisfying
\begin{align*}
    [L_m^{(i)},L_n^{(j)}]&=(m-n)L_{m+n}^{(i+j)},\\
    [J_m^{(i)},J_n^{(j)}]&=(j-i)J_{m+n}^{(i+j)},\\
    [J_m^{(i)},L_n^{(j)}]&=jL_{m+n}^{(i+j)}+mJ_{m+n}^{(i+j)},
\end{align*}
for $m,n\in\bZ$, where $L_m^{(i)}=J_m^{(j)}=0$ for $i\ge N$, $j\ge N-1$. In particular, if $N=2$, then $\mathcal{B}_2$ is generated by $\{L_m,P_m,J_m\mid m\in\bZ\}$ and satisfies the following non-vanishing relations:
\begin{align*}
    [L_m,L_n]=(m-n)L_{m+n},\ [L_m,P_n]=(m-n)P_{m+n},\ [L_m,J_n]=-nJ_{m+n},\ [J_m,P_n]=P_{m+n},
\end{align*}
for $m,n\in\bZ$. It is actually a \emph{3D Weyl--BMS algebra} (cf. \cite{AST,Yang}). In view of this, we call $\mathcal{B}_N$ a \emph{generalized Weyl--BMS algebra}.

In all the previous examples, $(\fg,B)$ is a Lie--Rinehart algebra. We conclude this subsection with an example in which this condition fails.

Take $B=\bC e\oplus\bC a\oplus\bC a^2$, $\fg=\bC h$ with $[h,h]=0$ and $Z=0$. Equip $\fg$ with the $B$-module structure defined by
\begin{equation}\label{nonLR.1}
    eh=h,\ \ ah=a^2h=0.
\end{equation}
Further let 
\begin{equation}\label{nonLR.2}
    \rho(h)=D,\ \ D(e)=0,\ \ D(a)=a,\ \ D(a^2)=2a^2,
\end{equation}
where $D=a\frac{d}{da}$ is a derivation of $B$. Since $\fg$ is abelian, $\rho$ is a Lie algebra homomorphism. Then by \eqref{nonLR.1} and \eqref{nonLR.2}, we have
\begin{equation*}
    [h,bh]=0=b[h,h]+\big(\rho(h)b\big)h,\ \ \forall\ b\in B.
\end{equation*}
Moreover, by Remark~\ref{rmk:Lie-Rine}, the defect
\begin{equation*}
    D_h(a,a)=\rho(ah)a-a\big(\rho(h)a\big)=-a^2\neq0.
\end{equation*}
So under this setting, $(\fg,B)$ is not a Lie--Rinehart algebra.

Take $U=\bC u$ and equip it with the $B$-module structure determined by $eu=u$, $au=a^2u=0$, and the $\fg$-module structure determined by $hu=u$. Suppose that
\begin{equation*}
    \sigma(u,u)=2a^2,\ \ \eta=0.
\end{equation*}
Then we can obtain
\begin{gather*}
    D\sigma(u,u)=4a^2=\sigma(hu,u)+\sigma(u,hu),\\
    \sigma(bu,u)=\sigma(u,bu)=b  \sigma(u,u).
\end{gather*}
So compatibility conditions \eqref{sigma.g} and \eqref{sigma.A} hold. One can check all the remaining compatibility conditions hold easily.

Set $L_m=e(m+1)$, $P_m=a(m+1)$, $Q_m=a^2(m+1)$, $J_m=h(m)$ and $G_r=u(r+\frac{1}{2})$ for $m\in\bZ$ and $r\in\bZ+\frac{1}{2}$. Then we get a Lie superalgebra $\cL(\fg,B,U)$ generated by $\{L_m,P_m,Q_m,J_m,G_r\mid m\in\bZ,\,r\in\bZ+\frac{1}{2}\}$, satisfying the following non-vanishing relations:
\begin{align*}
    &[L_m,L_n]=(m-n)L_{m+n},\ \ [L_m,P_n]=(m-n)P_{m+n},\ \ [L_m,Q_n]=(m-n)Q_{m+n},\\
    &[P_m,P_n]=(m-n)Q_{m+n},\ \ [L_m,J_n]=-nJ_{m+n},\ \ [J_m,P_n]=P_{m+n},\ \ [J_m,Q_n]=2Q_{m+n},\\
    &[L_m,G_r]=(\frac{m}{2}-r)G_{m+r},\ \ [J_m,G_r]=G_{m+r},\ \ [G_r,G_s]=2Q_{r+s},
\end{align*}
for $m,n\in\bZ$, $r,s\in\bZ+\frac{1}{2}$. To the best of our knowledge, this Lie superalgebra is new, and its even part is isomorphic to the centerless \emph{extended Schrödinger--Virasoro algebra} \cite{U}.

\begin{remark}
    In the examples of this subsection, we have chosen $A$ to be a truncated polynomial algebra. In general, a finite-dimensional local unital commutative associative algebra over $\bC$ need not be of this form. In dimension three, there are exactly two isomorphism classes of such algebras, represented by $\bC[\varepsilon]/(\varepsilon^3)$ and $\bC[x,y]/(x^2,xy,y^2)$ (cf. \cite{Poonen}). Moreover, for each fixed integer $n\ge7$, there are infinitely many isomorphism classes of $n$-dimensional local unital commutative associative algebras over $\bC$ (cf. \cite{Poonen,S}). A complete classification of the Lie algebras $\cL(\fg,A)$ with $\dim A\le 6$, even in the case where $(\fg,A)$ is a Lie--Rinehart algebra, seems out of reach at present.
\end{remark}

\appendix

\section{Proof of Theorem~\ref{thm:construction.L}}\label{app:Jacobi}
\hspace{1.5em}Throughout this appendix, let $\cQ_Z$ be as specified before Theorem~\ref{thm:construction.L}, and assume the structural, symmetry and invariance
properties stated there. Note that the maps $\kappa,\beta,\theta,\tau$ take values in $Z$. All central expressions below are therefore $Z$-valued.

The skew-symmetry of bilinear operation $[\cdot,\cdot]$ on $(\fg,\fg)$ and $(A,A)$ follows from that $\fg$ is a Lie algebra, $A$ is a commutative algebra and the bilinear forms $\kappa,\beta$ are symmetric. Since $\sigma,\tau$ are symmetric and $(m-n)\eta(u,v)=(n-m)\eta(v,u)$, we get the skew-symmetry of $[\cdot,\cdot]$ on $(U,U)$. Therefore, it remains to prove that the super Jacobi identity is equivalent to the listed compatibility conditions.

For convenience, we collect the identities used in the proof below. Let
\[
    x,y,z\in\fg,\ \ a,b,c\in A,\ \ u,v,w\in U.
\]
The numbered identities below retain their numbering from Section~\ref{sect:data}. These identities are assumed in the sufficiency direction and recovered from the Jacobi identity in the necessity direction. 

\medskip
\noindent\textbf{Compatibility conditions for the even datum.}

The mutual actions of $\fg$ and $A$ satisfy
\begin{gather}
    [x,ay]=a[x,y]+\big(\rho(x)a\big)y,\tag{\ref{compt.1}}\\
    \rho(ax)b-a\circ\big(\rho(x)b\big)=\rho(bx)a-b\circ\big(\rho(x)a\big).\tag{\ref{compt.2}}
\end{gather}
The even bilinear maps satisfy
\begin{align}
    \theta(ax,b)&=\theta(x,a\circ b)=\theta(bx,a),\tag{\ref{theta.1}}\\
    \theta\big(x,\rho(y)a\big)&=\theta([x,y],a)=-\theta\big(y,\rho(x)a\big),\tag{\ref{theta.2}}\\
    \beta\big(\rho(x)a,b\big)&=0=\beta\big(a,\rho(x)b\big),\tag{\ref{theta.3}}\\
    \theta([x,y],a)&=\kappa(x,ay)-\kappa(ax,y).\tag{\ref{theta.4}}
\end{align}

\medskip
\noindent\textbf{Compatibility of the actions on $U$.}

The two actions on $U$ satisfy
\begin{equation}
    (ax)u=x(au)-\frac32\big(\rho(x)a\big)u=a(xu)-\frac12\big(\rho(x)a\big)u.\tag{\ref{compt.xav}}
\end{equation}
In particular,
\begin{equation}\label{compt.xav.prime}
    x(au)-a(xu)=\big(\rho(x)a\big)u.\tag{\ref{compt.xav}$'$}
\end{equation}

\medskip
\noindent\textbf{Conditions involving $\tau$ and $\eta$.}

The maps $\tau$ and $\eta$ satisfy
\begin{align}
    &[x,\eta(u,v)]=\eta(xu,v)+\eta(u,xv),\tag{\ref{compt.uv.3}}\\
    &a\eta(u,v)=\eta(au,v)=\eta(u,av),\tag{\ref{compt.uv.4}}\\
    &\kappa\big(x,\eta(u,v)\big)=\tau(xu,v),\tag{\ref{compt.uv.5}}\\
    &\theta\big(\eta(u,v),a\big)=0.\tag{\ref{compt.uv.6}}
\end{align}

\medskip
\noindent\textbf{Conditions involving $\sigma$.}

The symmetric bilinear map $\sigma$ satisfies
\allowdisplaybreaks[4]
\begin{align}
    \rho(x)\sigma(u,v)&=\sigma(xu,v)+\sigma(u,xv),\tag{\ref{sigma.g}}\\
    2\sigma(u,v)x&=\eta(xu,v)-\eta(u,xv),\tag{\ref{sigma.x}}\\
    \theta\big(x,\sigma(u,v)\big)&=0,\tag{\ref{sigma.theta}}\\
    2a\circ\sigma(u,v)&=\sigma(au,v)+\sigma(u,av),\tag{\ref{sigma.A}}\\
    \rho\big(\eta(u,v)\big)a&=\sigma(au,v)-\sigma(u,av),\tag{\ref{sigma.rho}}\\
    \beta\big(a,\sigma(u,v)\big)&=\frac32\tau(au,v),\tag{\ref{sigma.beta}}\\
    3\big(\sigma(u,v)w-\sigma(u,w)v\big)&=\eta(u,v)w-\eta(u,w)v-2\eta(v,w)u.\tag{\ref{sigma.uvw}}
\end{align}
Interchanging $u$ and $v$ in \eqref{sigma.uvw}, and using the symmetry of $\sigma$ and the skew-symmetry of $\eta$, gives
\begin{equation}\label{sigma.uvw'}
    3\big(\sigma(u,v)w-\sigma(v,w)u\big)=-\eta(u,v)w-\eta(v,w)u-2\eta(u,w)v.\tag{\ref{sigma.uvw}$'$}
\end{equation}

We now verify the Jacobi identity for the generators of $\cL(\fg,A,U)$.

\afCase{1} $(\fg,\fg,\fg)$.

The Jacobi identity follows from the Jacobi identity in $\fg$ and the symmetry and invariance of $\kappa$, as in the standard affine Lie algebra construction.

\afCase{2} $(\fg,\fg,U)$.

Let $x,y\in\fg$, $u\in U$ and $m,n,k\in\bZ$. Using the defining relations in Theorem~\ref{thm:construction.L} and the centrality of $Z$, we obtain
\begin{align*}
    &[x(m),[y(n),u(k)]]-[[x(m),y(n)],u(k)]-[y(n),[x(m),u(k)]]\\
    &\quad=\big(x(yu)-[x,y]u-y(xu)\big)(m+n+k)=0,
\end{align*}
where the last equality follows from the assumed $\fg$-module structure on $U$.

\afCase{3} $(A,A,A)$.

Define a product on $A$ by
\begin{equation*}
    a*b=2(a\circ b).
\end{equation*}
Then $(A,*)$ is a commutative associative algebra, and $\beta$ is symmetric and satisfies
\begin{equation*}
    \beta(a*b,c)=\beta(a,b*c).
\end{equation*}
Moreover, the defining bracket \eqref{L:AA} becomes
\begin{equation*}
    [a(m),b(n)]=\frac12(m-n)(a*b)(m+n-1)+\binom{m}{3}\beta(a,b)\delta_{m+n-2,0}.
\end{equation*}
Thus the Jacobi identity in this case follows directly from \cite[proof of Theorem~4.1, Case~1]{LY}, whose calculation remains valid for $Z$-valued $\beta$, since it uses only its symmetry, invariance and linear combinations of its values.

\afCase{4} $(A,A,U)$.  

With the product $*$ defined in Case~3, define an action on $U$ by
\begin{equation*}
    a\diamond u=2au.
\end{equation*}
The $A$-module property of the original action gives
\begin{equation*}
    a\diamond(b\diamond u)=4a(bu)=4(a\circ b)u=(a*b)\diamond u.
\end{equation*}
Hence $U$ is a module over $(A,*)$. Then the bracket \eqref{L:AU} is of the form
\begin{equation*}
    [a(m),u(n)]=\frac14(m-2n)(a\diamond u)(m+n-1).
\end{equation*}
Together with \eqref{L:AA} and the centrality of $Z$, these are precisely the relations considered in \cite[proof of Theorem~4.1, Case~2]{LY}, whose calculation proves the Jacobi identity in this case.

\afCase{5} $(\fg,\fg,A)$.

Let $x,y\in\fg$, $a\in A$ and $m,n,k\in\bZ$, and set $s=m+n+k$. Write
\begin{equation*}
    J_1=[x(m),[y(n),a(k)]],\ \ J_2=[y(n),[x(m),a(k)]],\ \ J_3=[[x(m),y(n)],a(k)].
\end{equation*}
We have
\begin{align*}
    J_1&=[x(m),\big(\rho(y)a\big)(n+k)]+n[x(m),(ay)(n+k-1)]\\
    &=\big(\rho(x)\rho(y)a\big)(s)+m\Big(\big(\rho(y)a\big)x\Big)(s-1)+\binom{m}{2}\theta\big(x,\rho(y)a\big)\delta_{s-1,0}\\
    &\quad+n[x,ay](s-1)+mn\kappa(x,ay)\delta_{s-1,0}.
\end{align*}
Using permutations, we can deduce that
\begin{align*}
    J_2&=\big(\rho(y)\rho(x)a\big)(s)+n\Big(\big(\rho(x)a\big)y\Big)(s-1)+m[y,ax](s-1)\\
    &\quad+\left(\binom{n}{2}\theta\big(y,\rho(x)a\big)+mn\kappa(y,ax)\right)\delta_{s-1,0}.
\end{align*}
Moreover,
\begin{align*}
    J_3&=[[x,y](m+n),a(k)]\\
    &=\big(\rho([x,y])a\big)(s)+(m+n)([x,y]a)(s-1)+\binom{m+n}{2}\theta([x,y],a)\delta_{s-1,0}.
\end{align*}
Set
\begin{equation*}
    J=J_1-J_2-J_3=P_A(s)+P_\fg(s-1)+P_Z\delta_{s-1,0},
\end{equation*}
where $P_A\in A$, $P_\fg\in\fg$ and $P_Z\in Z$. The $A$-component is determined by 
\begin{align*}
    P_A=\rho(x)\rho(y)a-\rho(y)\rho(x)a-\rho([x,y])a=0,
\end{align*}
since $\rho$ is a Lie algebra homomorphism. The $\fg$-component is 
\begin{equation}\label{J.ggU}
    P_\fg=n[x,ay]-m[y,ax]+m\big(\rho(y)a\big)x-n\big(\rho(x)a\big)y-(m+n)a[x,y].
\end{equation}
By \eqref{compt.1}, we obtain
\begin{equation*}
    [x,ay]=a[x,y]+\big(\rho(x)a\big)y,\ \ [y,ax]=-a[x,y]+\big(\rho(y)a\big)x.
\end{equation*}
Substituting these two identities into \eqref{J.ggU} shows that $P_\fg$ also vanishes. Finally, the central component is $P_Z\delta_{s-1,0}$, where $P_Z\in Z$ is given by
\begin{align*}
    P_Z&=\binom{m}{2}\theta\big(x,\rho(y)a\big)-\binom{n}{2}\theta\big(y,\rho(x)a\big)-\binom{m+n}{2}\theta([x,y],a)\\
    &\quad+mn\big(\kappa(x,ay)-\kappa(y,ax)\big).
\end{align*}
It follows from \eqref{theta.2} that
\begin{equation*}
    \theta\big(x,\rho(y)a\big)=\theta([x,y],a),\ \ \theta\big(y,\rho(x)a\big)=-\theta([x,y],a).
\end{equation*}
By the symmetry of $\kappa$ and \eqref{theta.4}, we also have
\begin{equation*}
    \kappa(x,ay)-\kappa(y,ax)=\theta([x,y],a).
\end{equation*}
Consequently,
\begin{equation*}
    P_Z=\left(\binom{m}{2}+\binom{n}{2}-\binom{m+n}{2}+mn\right)\theta([x,y],a)=0.
\end{equation*}
Thus $P_A=P_\fg=P_Z=0$, then the Jacobi identity holds in this case.

Conversely, if the Jacobi identity holds, then take $(m,n,k)=(0,1,0)$ in \eqref{J.ggU}, which implies that
\begin{equation*}
    [x,ay]-a[x,y]-\big(\rho(x)a\big)y=0.
\end{equation*}
So we obtain \eqref{compt.1}. For the central component, write
\begin{align*}
    T_1=\theta\big(x,\rho(y)a\big),\ \ T_2=\theta\big(y,\rho(x)a\big),\ \ T_3=\theta([x,y],a),\ \ T_4=\kappa(x,ay)-\kappa(y,ax).
\end{align*}
Using $\binom{m+n}{2}=\binom{m}{2}+\binom{n}{2}+mn$, we get
\begin{equation*}
    P_Z=\binom{m}{2}(T_1-T_3)-\binom{n}{2}(T_2+T_3)+mn(T_4-T_3).
\end{equation*}
Letting $s=1$, the $Z$-component of the Jacobi identity implies $P_Z=0$. The choices $(m,n,k)=(2,0,-1)$, $(0,2,-1)$, $(1,1,-1)$ show, respectively,
\begin{equation*}
    T_1=T_3,\ \ T_2=-T_3,\ \ T_4=T_3.
\end{equation*}
The first two equalities give \eqref{theta.2}; the last leads to \eqref{theta.4} by the symmetry of $\kappa$.

Hence, the Jacobi identity for $(\fg,\fg,A)$ is equivalent to $\eqref{compt.1}$, \eqref{theta.2} and \eqref{theta.4}.

\afCase{6} $(\fg,A,A)$.

Let $x\in\fg$, $a,b\in A$ and $m,n,k\in\bZ$, and set
$s=m+n+k$. Write
\begin{align*}
    J_1=[x(m),[a(n),b(k)]],\ \ J_2=[[x(m),a(n)],b(k)],\ 
    \ J_3=[a(n),[x(m),b(k)]].
\end{align*}
Using the defining brackets in Theorem~\ref{thm:construction.L} and the centrality of $Z$, we obtain
\begin{align*}
    J_1&=(n-k)[x(m),(a\circ b)(n+k-1)]\\
    &=(n-k)\big(\rho(x)(a\circ b)\big)(s-1)+m(n-k)\big((a\circ b)x\big)(s-2)\\
    &\quad+(n-k)\binom{m}{2}\theta(x,a\circ b)\delta_{s-2,0}.
\end{align*}
For the second term,
\begin{align*}
    J_2&=[\big(\rho(x)a\big)(m+n),b(k)]+m[(ax)(m+n-1),b(k)]\\
    &=\big((m+n-k)\big(b\circ\rho(x)a\big)+m\rho(ax)b\big)(s-1)\\
    &\quad+m(m+n-1)\big(b(ax)\big)(s-2)\\
    &\quad+\left(\binom{m+n}{3}\beta\big(\rho(x)a,b\big)+m\binom{m+n-1}{2}\theta(ax,b)\right)\delta_{s-2,0}.
\end{align*}
Similarly,
\begin{align*}
    J_3&=[a(n),\big(\rho(x)b\big)(m+k)]+m[a(n),(bx)(m+k-1)]\\
    &=\big((n-m-k)\big(a\circ\rho(x)b\big)-m\rho(bx)a\big)(s-1)\\
    &\quad-m(m+k-1)\big(a(bx)\big)(s-2)\\
    &\quad+\left(\binom{n}{3}\beta\big(a,\rho(x)b\big)-m\binom{m+k-1}{2}\theta(bx,a)\right)\delta_{s-2,0}.
\end{align*}

Set
\begin{align*}
    J=J_1-J_2-J_3=P_A(s-1)+P_{\fg}(s-2)+P_Z\delta_{s-2,0},
\end{align*}
where $P_A\in A$, $P_{\fg}\in\fg$ and $P_Z\in Z$. The $A$-component is determined by
\begin{align*}
    P_A&=(n-k)\rho(x)(a\circ b)-(m+n-k)\big(b\circ\rho(x)a\big)-m\rho(ax)b\\
    &\quad-(n-m-k)\big(a\circ\rho(x)b\big)+m\rho(bx)a\\
    &=(n-k)\big(\rho(x)(a\circ b)-a\circ\rho(x)b-b\circ\rho(x)a\big)\\
    &\quad+m\big(a\circ\rho(x)b+\rho(bx)a-b\circ\rho(x)a-\rho(ax)b\big)\\
    &=0.
\end{align*}
Here the first term vanishes because $\rho(x)$ is a derivation, and the second vanishes by \eqref{compt.2}. For the $\fg$-component, the symmetric $A$-bimodule structure gives
\begin{align*}
    P_{\fg}&=m(n-k)\big((a\circ b)x\big)-m(m+n-1)\big(b(ax)\big)+m(m+k-1)\big(a(bx)\big)\\
    &=m\big((n-k)-(m+n-1)+(m+k-1)\big)\big((a\circ b)x\big)=0.
\end{align*}
Finally, the $Z$-component is $P_Z\delta_{s-2,0}$, where $P_Z\in Z$ is given by
\begin{align*}
    P_Z&=(n-k)\binom{m}{2}\theta(x,a\circ b)-\binom{m+n}{3}\beta\big(\rho(x)a,b\big)-m\binom{m+n-1}{2}\theta(ax,b)\\
    &\quad-\binom{n}{3}\beta\big(a,\rho(x)b\big)+m\binom{m+k-1}{2}\theta(bx,a).
\end{align*}
By \eqref{theta.1} and \eqref{theta.3}, this becomes
\begin{align*}
    P_Z&=\left((n-k)\binom{m}{2}-m\binom{m+n-1}{2}+m\binom{m+k-1}{2}\right)\theta(x,a\circ b)\\
    &=\frac12m(n-k)(2-s)\theta(x,a\circ b).
\end{align*}
Consequently, $P_Z\delta_{s-2,0}=0$ for all $m,n,k\in\bZ$. Thus $J=0$, proving the Jacobi identity in this case.

Conversely, suppose that $J=0$ for all $m,n,k\in\bZ$. Comparing components gives $P_A=P_\fg=P_Z\delta_{s-2,0}=0$. For $P_A=0$, we take $(m,n,k)=(1,0,0)$ and obtain \eqref{compt.2}. For the central coefficient $P_Z$, write
\begin{gather*}
    B_1=\beta\big(\rho(x)a,b\big),\ \ B_2=\beta\big(a,\rho(x)b\big),\\
    T_1=\theta(x,a\circ b),\ \ T_2=\theta(ax,b),\ \ T_3=\theta(bx,a).
\end{gather*}
Then by $P_Z\delta_{s-2,0}=0$, we can set $s=2$ and obtain the following equation
\begin{equation}\label{gAA.C}
    0=\binom{m+n}{3}B_1+\binom{n}{3}B_2-(m+2n-2)\binom{m}{2}T_1+m\binom{m+n-1}{2}T_2-m\binom{1-n}{2}T_3.
\end{equation}
Taking $(m,n)=(0,3)$, $(1,2)$, $(2,1)$ and $(3,0)$ in \eqref{gAA.C} gives, respectively, 
\begin{equation*}
    \begin{cases}
        B_1+B_2=0,\\
        B_1+T_2-T_3=0,\\
        B_1-2T_1+2T_2=0,\\
        B_1-3T_1+3T_2=0.
    \end{cases}
\end{equation*}
It implies that
\begin{equation*}
    B_1=B_2=0,\ \ T_1=T_2=T_3.
\end{equation*}
The equalities $B_1=B_2=0$ give \eqref{theta.3}, whereas $T_1=T_2=T_3$ lead to \eqref{theta.1}.

Hence, the Jacobi identity for $(\fg,A,A)$ is equivalent to $\eqref{compt.2}$, \eqref{theta.1} and \eqref{theta.3}.

\afCase{7} $(\fg,A,U)$.

Let $x\in\fg$, $a\in A$, $u\in U$ and $m,n,k\in\bZ$, and set $s=m+n+k$. Using the defining brackets in Theorem~\ref{thm:construction.L} and the centrality of $Z$, we obtain
\begin{align*}
    [x(m),[a(n),u(k)]]&=\frac{1}{2}(n-2k)[x(m),(au)(n+k-1)]\\
    &=\frac{1}{2}(n-2k)\big(x(au)\big)(s-1),\\
    [[x(m),a(n)],u(k)]&=[\big(\rho(x)a\big)(m+n),u(k)]+m[(ax)(m+n-1),u(k)]\\
    &=\Big(\frac{1}{2}(m+n-2k)\big(\rho(x)a\big)u+m(ax)u\Big)(s-1),\\
    [a(n),[x(m),u(k)]]&=[a(n),(xu)(m+k)]\\
    &=\frac{1}{2}(n-2m-2k)\big(a(xu)\big)(s-1).
\end{align*}
Then
\begin{align}
    [x(m),[&a(n),u(k)]]-[[x(m),a(n)],u(k)]-[a(n),[x(m),u(k)]]\notag\\
    &=\frac{1}{2}(n-2k)\Big(x(au)-\big(\rho(x)a\big)u-a(xu)\Big)(s-1)\notag\\
    &\quad+m\Big(a(xu)-(ax)u-\frac{1}{2}\big(\rho(x)a\big)u\Big)(s-1).\label{J.xau}
\end{align}
In \eqref{J.xau}, the first term vanishes by \eqref{compt.xav.prime} and the second vanishes by \eqref{compt.xav}. Then the Jacobi identity holds in this case.

Conversely, if the Jacobi identity holds, then taking $(m,n,k)=(0,1,0)$ and $(1,0,0)$ in \eqref{J.xau}, respectively, gives 
\begin{align*}
    x(au)=\big(\rho(x)a\big)u+a(xu),\ \ a(xu)=(ax)u+\frac{1}{2}\big(\rho(x)a\big)u.
\end{align*}
It follows that $(ax)u=x(au)-\frac{3}{2}\big(\rho(x)a\big)u$. Combining these equations, we have \eqref{compt.xav}. 

So the Jacobi identity for $(\fg,A,U)$ is equivalent to \eqref{compt.xav}.

\afCase{8} $(\fg,U,U)$.

Let $x\in\fg$, $u,v\in U$ and $m,n,k\in\bZ$, and set $s=m+n+k$. Write
\begin{equation*}
    J_1=[x(m),[u(n),v(k)]],\ \ J_2=[[x(m),u(n)],v(k)],\ \ J_3=[u(n),[x(m),v(k)]].
\end{equation*}
Using the defining brackets in Theorem~\ref{thm:construction.L} and the centrality of $Z$, we obtain
\allowdisplaybreaks[4]
\begin{align*}
    J_1&=[x(m),\sigma(u,v)(n+k)]+\frac{1}{2}(n-k)[x(m),\eta(u,v)(n+k-1)]\\
    &=\big(\rho(x)\sigma(u,v)\big)(s)+m\big(\sigma(u,v)x\big)(s-1)+\binom{m}{2}\theta\big(x,\sigma(u,v)\big)\delta_{s-1,0}\\
    &\quad+\frac{1}{2}(n-k)\left([x,\eta(u,v)](s-1)+m\kappa\big(x,\eta(u,v)\big)\delta_{s-1,0}\right)\\
    &=\big(\rho(x)\sigma(u,v)\big)(s)\\
    &\quad+\big(m\sigma(u,v)x+\frac{1}{2}(n-k)[x,\eta(u,v)]\big)(s-1)\\
    &\quad+\left(\binom{m}{2}\theta\big(x,\sigma(u,v)\big)+\frac{1}{2}m(n-k)\kappa\big(x,\eta(u,v)\big)\right)\delta_{s-1,0},\\
    J_2&=[(xu)(m+n),v(k)]\\
    &=\sigma(xu,v)(s)+\frac{1}{2}(m+n-k)\eta(xu,v)(s-1)+\binom{m+n}{2}\tau(xu,v)\delta_{s-1,0},\\
    J_3&=[u(n),(xv)(m+k)]\\
    &=\sigma(u,xv)(s)+\frac{1}{2}(n-m-k)\eta(u,xv)(s-1)+\binom{n}{2}\tau(u,xv)\delta_{s-1,0}.
\end{align*}

Set
\begin{equation*}
    J=J_1-J_2-J_3=P_A(s)+P_\fg(s-1)+P_Z\delta_{s-1,0},
\end{equation*}
where $P_A\in A$, $P_\fg\in\fg$ and $P_Z\in Z$. The $A$-component is determined by
\begin{align*}
    P_A=\rho(x)\sigma(u,v)-\sigma(xu,v)-\sigma(u,xv).
\end{align*}
Then $P_A=0$ if and only if \eqref{sigma.g} holds. The $\fg$-component is 
\begin{align*}
    P_\fg&=m\big(\sigma(u,v)x-\frac{1}{2}\eta(xu,v)+\frac{1}{2}\eta(u,xv)\big)\\
    &\qquad+\frac{1}{2}(n-k)\big([x,\eta(u,v)]-\eta(xu,v)-\eta(u,xv)\big).
\end{align*}
If \eqref{sigma.x} and \eqref{compt.uv.3} hold, then $P_\fg=0$. Conversely, letting $P_\fg=0$ and taking $(m,n,k)=(1,0,0)$, $(0,1,0)$, respectively, we arrive at \eqref{compt.uv.3} and \eqref{sigma.x}.

For the $Z$-component, write
\begin{equation*}
    T_1=\theta\big(x,\sigma(u,v)\big),\ \ T_2=\kappa\big(x,\eta(u,v)\big),\ \ T_3=\tau(xu,v),\ \ T_4=\tau(u,xv).
\end{equation*}
Then substituting them into the central component with $k=1-m-n$, we have
\begin{align*}
    P_Z&=\frac{1}{2}m(m-1)T_1+\frac{1}{2}m(m+2n-1)T_2-\frac{1}{2}(m+n)(m+n-1)T_3-\frac{1}{2}n(n-1)T_4\\
    &=\frac{1}{2}m(m-1)(T_1+T_2-T_3)+mn(T_2-T_3)-\frac{1}{2}n(n-1)(T_3+T_4)\\
    &=\frac{1}{2}m(m-1)(T_1+T_2-T_3)+mn(T_2-T_3).
\end{align*}
The third term vanishes by the $\fg$-invariance of $\tau$. Moreover, if \eqref{compt.uv.5} and \eqref{sigma.theta} hold, then $P_Z=0$. Conversely, letting $P_Z=0$ and taking $(m,n)=(1,1)$ and $(2,0)$, respectively, we can get \eqref{compt.uv.5} and \eqref{sigma.theta}. Thus, $P_Z\delta_{s-1,0}=0$ for all $m,n,k\in\bZ$ if and only if \eqref{compt.uv.5} and \eqref{sigma.theta} hold.

Hence the Jacobi identity for $(\fg,U,U)$ is equivalent to \eqref{compt.uv.3}, \eqref{compt.uv.5} and \eqref{sigma.g}--\eqref{sigma.theta}.

\afCase{9} $(A,U,U)$.

Let $a\in A$, $u,v\in U$ and $m,n,k\in\bZ$, and set $s=m+n+k$. Write
\begin{equation*}
    J_1=[a(m),[u(n),v(k)]],\ \ J_2=[[a(m),u(n)],v(k)],\ \ J_3=[u(n),[a(m),v(k)]].
\end{equation*}
Using the defining brackets in Theorem~\ref{thm:construction.L} and the centrality of $Z$, we obtain
\allowdisplaybreaks[1]
\begin{align*}
    J_1&=[a(m),\sigma(u,v)(n+k)]+\frac12(n-k)[a(m),\eta(u,v)(n+k-1)]\\
    &=(m-n-k)\big(a\circ\sigma(u,v)\big)(s-1)+\binom{m}{3}\beta\big(a,\sigma(u,v)\big)\delta_{s-2,0}\\
    &\quad-\frac{1}{2}(n-k)\Big(\big(\rho(\eta(u,v))a\big)(s-1)+(n+k-1)\big(a\eta(u,v)\big)(s-2)\\
    &\qquad+\binom{n+k-1}{2}\theta\big(\eta(u,v),a\big)\delta_{s-2,0}\Big)\\
    &=\Big((m-n-k)a\circ\sigma(u,v)-\frac12(n-k)\rho\big(\eta(u,v)\big)a\Big)(s-1)\\
    &\quad-\frac12(n-k)(n+k-1)\big(a\eta(u,v)\big)(s-2)\\
    &\quad+\left(\binom m3\beta\big(a,\sigma(u,v)\big)-\frac12(n-k)\binom{n+k-1}{2}\theta\big(\eta(u,v),a\big)\right)\delta_{s-2,0}.
\end{align*}
For the second term,
\begin{align*}
    J_2&=\frac12(m-2n)[(au)(m+n-1),v(k)]\\
    &=\frac12(m-2n)\sigma(au,v)(s-1)+\frac14(m-2n)(m+n-k-1)\eta(au,v)(s-2)\\
    &\quad+\frac12(m-2n)\binom{m+n-1}{2}\tau(au,v)\delta_{s-2,0}.
\end{align*}
Similarly,
\begin{align*}
    J_3&=\frac12(m-2k)[u(n),(av)(m+k-1)]\\
    &=\frac12(m-2k)\sigma(u,av)(s-1)+\frac14(m-2k)(n-m-k+1)\eta(u,av)(s-2)\\
    &\quad+\frac12(m-2k)\binom n2\tau(u,av)\delta_{s-2,0}.
\end{align*}

Set
\begin{equation*}
    J=J_1-J_2-J_3=P_A(s-1)+P_\fg(s-2)+P_Z\delta_{s-2,0},
\end{equation*}
where $P_A\in A$, $P_\fg\in\fg$ and $P_Z\in Z$. For convenience, write
\begin{gather*}
    S_1=a\circ\sigma(u,v),\ \ S_2=\sigma(au,v),\ \ S_3=\sigma(u,av),\ \ R=\rho(\eta(u,v))a,\\
    E_1=a\eta(u,v),\ \ E_2=\eta(au,v),\ \ E_3=\eta(u,av).
\end{gather*}
Then
\begin{align}
    P_A&=(m-n-k)S_1-\frac12(n-k)R-\frac12(m-2n)S_2-\frac12(m-2k)S_3\notag\\
    &=\frac12(m-n-k)(2S_1-S_2-S_3)+\frac12(n-k)(S_2-S_3-R).\label{J.auu.A}
\end{align}
By \eqref{sigma.A} and \eqref{sigma.rho}, we have $P_A=0$. For the $\fg$-component,
\begin{align}
    P_\fg&=-\frac12(n-k)(n+k-1)E_1-\frac14(m-2n)(m+n-k-1)E_2\notag\\
    &\quad-\frac14(m-2k)(n-m-k+1)E_3.\label{J.auu.g}
\end{align}
By \eqref{compt.uv.4}, we get $E_1=E_2=E_3$. Consequently, $P_\fg=0$, since
\begin{align*}
    &\quad2(n-k)(n+k-1)+(m-2n)(m+n-k-1)\\
    &\qquad+(m-2k)(n-m-k+1)\\
    &=2(n-k)(n+k-1)+(m-2n)\big((n-k)+(m-1)\big)\\
    &\qquad+(m-2k)\big((n-k)-(m-1)\big)\\
    &=2(n-k)(n+k-1)+(2m-2n-2k)(n-k)-2(n-k)(m-1)\\
    &=2(n-k)\big((n+k-1)+(m-n-k)-(m-1)\big)\\
    &=0.
\end{align*}
For the $Z$-component, write
\begin{gather*}
    B=\beta\big(a,\sigma(u,v)\big),\ \ H=\theta\big(\eta(u,v),a\big),\ \ T_1=\tau(au,v),\ \ T_2=\tau(u,av).
\end{gather*}
Then
\begin{align}
    P_Z&=\binom m3 B-\frac12(n-k)\binom{n+k-1}{2}H-\frac12(m-2n)\binom{m+n-1}{2}T_1
      \notag\\
      &\quad-\frac12(m-2k)\binom n2 T_2.\label{J.auu.C}
\end{align}
By \eqref{compt.uv.6} and \eqref{sigma.beta},
\begin{equation*}
    H=0,\ \ B=\frac32T_1.
\end{equation*}
Moreover, the $A$-invariance of $\tau$ gives $T_1=T_2$. Therefore,
\begin{align*}
    P_Z&=\frac32\binom m3 T_1-\frac12(m-2n)\binom{m+n-1}{2}T_1-\frac12(m-2k)\binom n2
      T_1\\
    &=\binom n2(s-2)T_1.
\end{align*}
Thus $P_Z\delta_{s-2,0}=0$, and hence $J=0$.

Conversely, suppose that $J=0$ for all $m,n,k\in\bZ$. Comparing components gives
\begin{equation*}
    P_A=P_\fg=P_Z\delta_{s-2,0}=0.
\end{equation*}
Taking $(m,n,k)=(1,0,0)$ and $(1,1,0)$ in \eqref{J.auu.A}, respectively, yields
\begin{equation*}
    2S_1-S_2-S_3=0,\ \  S_2-S_3-R=0,
\end{equation*}
which are \eqref{sigma.A} and \eqref{sigma.rho}.
For the $\fg$-component, the choices $(m,n,k)=(0,2,0)$ and $(2,0,0)$ in \eqref{J.auu.g}
give, respectively,
\begin{equation*}
    E_2-E_1=0,\ \ \frac12(E_3-E_2)=0.
\end{equation*}
Thus $E_1=E_2=E_3$, proving \eqref{compt.uv.4}. Finally, restrict to $s=2$, so that $P_Z=0$. The assumed $A$-invariance of $\tau$ gives $T_1=T_2$. Taking $(m,n,k)=(2,1,-1)$ and $(3,0,-1)$ in \eqref{J.auu.C} gives, respectively,
\begin{equation}
    H=0,\ \ B-\frac32H-\frac32T_1=0.
\end{equation}
Hence $H=0$ and $B=\frac32T_1$, proving \eqref{compt.uv.6} and \eqref{sigma.beta}, respectively.

Therefore, the Jacobi identity for $(A,U,U)$ is equivalent to \eqref{compt.uv.4}, \eqref{compt.uv.6} and \eqref{sigma.A}--\eqref{sigma.beta}.

\afCase{10} $(U,U,U)$.

Let $u,v,w\in U$ and $m,n,k\in\bZ$, and set $s=m+n+k$. Since all three generators are odd, the Jacobi identity is
\begin{equation*}
    [u(m),[v(n),w(k)]]=[[u(m),v(n)],w(k)]-[v(n),[u(m),w(k)]].
\end{equation*}
Write
\begin{equation*}
    J_1=[u(m),[v(n),w(k)]],\ \ J_2=[[u(m),v(n)],w(k)],\ \ J_3=[v(n),[u(m),w(k)]].
\end{equation*}
For convenience, set
\begin{gather*}
    S_1=\sigma(u,v)w,\ \ S_2=\sigma(u,w)v,\ \ S_3=\sigma(v,w)u,\\
    E_1=\eta(u,v)w, \ \ E_2=\eta(u,w)v,\ \ E_3=\eta(v,w)u.
\end{gather*}
Using the defining brackets in Theorem~\ref{thm:construction.L}, we obtain
\begin{align*}
    J_1&=[u(m),\sigma(v,w)(n+k)]+\frac12(n-k)[u(m),\eta(v,w)(n+k-1)]\\
    &=\big(\frac12(2m-n-k)S_3-\frac12(n-k)E_3\big)(s-1).
\end{align*}
Similarly,
\begin{align*}
    J_2&=[\sigma(u,v)(m+n),w(k)]+\frac12(m-n)[\eta(u,v)(m+n-1),w(k)]\\
    &=\big(\frac12(m+n-2k)S_1+\frac12(m-n)E_1\big)(s-1),\\
    J_3&=[v(n),\sigma(u,w)(m+k)]+\frac12(m-k)[v(n),\eta(u,w)(m+k-1)]\\
    &=\big(\frac12(2n-m-k)S_2-\frac12(m-k)E_2\big)(s-1).
\end{align*}

Set $J=J_1-J_2+J_3$. Collecting coefficients gives
\begin{align}
    2J&=\Big(m(2S_3-S_1-S_2-E_1-E_2)+n(2S_2-S_1-S_3+E_1-E_3)\notag\\
    &\qquad+k(2S_1-S_2-S_3+E_2+E_3)
      \Big)(s-1).\label{J.uuu}
\end{align}
Define
\begin{align*}
    R_1&=3(S_1-S_2)-E_1+E_2+2E_3,\\
    R_2&=3(S_1-S_3)+E_1+2E_2+E_3.
\end{align*}
Then \eqref{J.uuu} can be rewritten as
\begin{equation}\label{J.uuu.1}
    J=\frac16\Big((m-2n+k)R_1+(-2m+n+k)R_2\Big)(s-1).
\end{equation}

Suppose that \eqref{sigma.uvw} or \eqref{sigma.uvw'} holds. Then we have $R_1=R_2=0$. Consequently, $J=0$ for all $m,n,k\in\bZ$.

Conversely, suppose that the super Jacobi identity holds for all $u,v,w\in U$ and $m,n,k\in\bZ$. Taking $(m,n,k)=(0,-1,1)$ in \eqref{J.uuu.1} and comparing the coefficients, we obtain $R_1=0$, which is \eqref{sigma.uvw}.

Hence, the Jacobi identity for $(U,U,U)$ is equivalent to \eqref{sigma.uvw}.

The above ten cases exhaust the verification of the Jacobi identity. This completes the proof of Theorem~\ref{thm:construction.L}.

\end{document}